\documentclass[a4paper,11pt,reqno]{article}
\usepackage{amsmath,amssymb,amsthm} 
\usepackage{graphics} 
\usepackage{enumerate}
\usepackage{txfonts}

\usepackage[colorlinks,citecolor=blue]{hyperref}

\usepackage[capitalize]{cleveref}
\crefname{equation}{}{}

\usepackage{epsfig}
\usepackage{color}
\usepackage[english]{babel}

\numberwithin{equation}{section}

\newtheorem{theorem}{Theorem}[section]
\newtheorem{lemma}[theorem]{Lemma}
\newtheorem{proposition}[theorem]{Proposition}

\newtheorem{corollary}[theorem]{Corollary}
\newtheorem{assumption}[theorem]{Assumption}

\theoremstyle{definition}

\newtheorem{remark}[theorem]{Remark}

\newcommand{\hel} {
\hskip2.5pt{\vrule height7pt width.5pt depth0pt}
\hskip-.2pt\vbox{\hrule height.5pt width7pt depth0pt}
\, }
\newcommand{\restr}{\hel}

\newcommand{\Var}{\operatorname{Var}}

\renewcommand{\Cap}{\operatorname{Cap}}
\newcommand{\R}{\mathbb{R}}
\newcommand{\N}{\mathbb{N}}

\newcommand{\eps}{\varepsilon}

\newcommand{\be}{\begin{equation*}}
\newcommand{\ee}{\end{equation*}}

\DeclareMathOperator{\supp}{supp}

\newcommand\lt{\left}
\newcommand{\cY}{\mathcal{Y}}
\newcommand\rt{\right}

\def\les{\lesssim}
\def\ges{\gtrsim}

\renewcommand{\supp}{\operatorname{supp}}

\renewcommand{\Cap}{\operatorname{Cap}}

\def\EE{\mathbb{E}}
\def\PP{\mathbb{P}}
\def\diam{\operatorname{diam}}

\newcommand{\osc}{\operatorname{osc}}

\newcommand{\cS}{\mathcal{S}}

\newcommand{\bra}[1]{\left( #1 \right)}
\newcommand{\sqa}[1]{\left[ #1 \right]}
\newcommand{\cur}[1]{\left\{ #1 \right\}}

\newcommand{\abs}[1]{\left| #1 \right|}
\newcommand{\nor}[1]{\left\| #1 \right\|}

\newcommand{\B}{\mathrm{B}}

\newcommand{\WN}{\mathsf{W}}
\newcommand{\WD}{\mathsf{Wb}}
\newcommand{\Wdist}{\mathsf{Wd}}

\newcommand{\dist}{\mathsf{d}}

\newcommand{\cost}{\mathsf{w}}
\newcommand{\costb}{\mathsf{wb}}

\newcommand{\Lip}{\operatorname{Lip}}

\newcommand{\cQ}{\mathcal{Q}}

\newcommand{\cX}{\mathcal{X}}

\renewcommand{\PP}{\mathbb{P}}

\newcommand{\costm}{\mathsf{c}}

\renewcommand{\div}{\operatorname{div}}

\newcommand{\bist}{\mathsf{b}}
\newcommand{\sfQ}{\mathsf{Q}}
\newcommand{\err}{\mathsf{err}}

\usepackage[hyperref,doi,url=false,doi=false, isbn=false, giveninits=true, backref,  style=numeric,maxbibnames=99, backend=biber]{biblatex}
\bibliography{OT.bib}

\title{Asymptotics for Random Quadratic Transportation Costs \\
}
\author{Martin Huesmann \\
\emph{Universität Münster}\\
\and
Michael Goldman \\
  \emph{CNRS and \'Ecole polytechnique}\\
 \and
 Dario Trevisan\\
 \emph{Università di Pisa}\\
}

\usepackage{fancyhdr}
\begin{document}

\maketitle

 \newcounter{step}

\begin{abstract}
We establish the validity of asymptotic limits for the general transportation problem between random i.i.d.\  points and their common distribution, with respect to the squared Euclidean distance cost, in any dimension larger than three. Previous results were essentially limited to the two (or one) dimensional case, or to distributions whose absolutely continuous part is uniform.

The proof relies upon recent advances in the stability theory of optimal transportation, combined with functional analytic techniques and some ideas from quantitative stochastic homogenization. The key tool we develop is a quantitative upper bound for the usual quadratic optimal transportation problem in terms of its boundary variant, where points can be freely transported along the boundary. The methods we use are applicable to more general random measures, including occupation measure of Brownian paths, and may open the door to further progress on challenging problems at the interface of analysis, probability, and discrete mathematics.

\vspace{1em}
\end{abstract}



\section{Introduction}

The assignment problem (or bipartite matching) is a classic optimization problem that involves finding an optimal correspondence between two sets of objects, such that a total matching cost is minimized. It counts innumerate applications and is subject of intense research by several communities, from operation research and algorithm theory, combinatorics and graph theory \parencite{lovasz2009matching}, probability and statistics \parencite{talagrand2022upper} and even theoretical physics \parencite{mezard1987spin, mezard2009information}. In the Euclidean formulation of the problem, the two sets of objects are families of points $\bra{x_i}_{i=1}^n$, $\bra{y_i}_{i=1}^n \subseteq \R^d$,  and  the matching cost is defined, for a given parameter $p>0$, as
\begin{equation}\label{eq:assignment} \min_{\sigma \in \mathcal{S}_{n}} \sum_{i=1}^n | x_i - y_{\sigma(i)} |^p,\end{equation}
where $\mathcal{S}_{n}$ denotes the set of permutations and $|\cdot|$ the Euclidean norm.

By the Birkhoff-Von Neumann Theorem,  a natural linear programming reformulation of \Cref{eq:assignment} is given by the optimal transport problem \parencite{VilOandN}, where $\sigma$ is replaced by a transport plan or coupling, and the cost \Cref{eq:assignment} becomes a special instance of the Wasserstein cost of order $p$, denoted by $\WN^p(\mu, \lambda)$, between the two empirical measures $\mu:= \sum_{i=1}^n \delta_{x_i}$ and $\lambda:= \sum_{i=1}^n \delta_{y_i}$ (see \Cref{sec:wass} for precise notation). Such ``relaxed'' formulation allows also for non-discrete measures, and opens the way to more general applications. For example, by letting $\lambda$ be absolutely continuous and keeping $\mu$ discrete over a set of $n$ points, one can interpret $\WN^p(\mu, \lambda)$ as a quantization cost \parencite{graf2007foundations, dereich2013constructive}. In this context it is then quite natural to assume e.g.\ that $x_i = X_i$ are randomly sampled from a given probability distribution, e.g., $\lambda$ normalized to be a probability.

There is in fact a rich literature studying random instances of combinatorial optimization problems in Euclidean spaces \parencite{steele1997probability, yukich2006probability}. For bipartite problems, like the assignment problem, it is well-known that classical methods encounter limitations due to local fluctuations in the number of samples between the two sets. This can lead to unexpected scaling behaviors of the cost as $n$ grows, as first observed in \parencite{ajtai1984optimal} for  uniform i.i.d.\ samples on the unit square $(0,1)^2$ and $p=1$. Subsequent contributions by several authors \parencite{talagrand1992matching, dobric1995asymptotics, boutet2002almost, talagrand2022upper, barthe2013combinatorial, dereich2013constructive, BoLe16, Le17, BobLe19, bobkov2020transport} lead to a more complete understanding of the whole picture, for general dimensions $d$, cost exponents $p$ and sample distributions $\lambda$, providing asymptotic upper and lower bounds. We report here the combination  of \parencite{dereich2013constructive,barthe2013combinatorial,goldman2022optimal} (see also \parencite{pieroni-gaussian}).

\begin{theorem}\label{thm:dereich}
 For any $d\ge 3$ and  $p\ge 1$, there exist constants $0<\underline{\mathsf{c}}(p,d) \le \overline{\mathsf{c}}(p,d)<\infty$ such that the following holds for every probability measure  $\lambda$ on $\R^d$.  Assume that either
 \begin{itemize}
  \item[(i)]  $ p<d/2$  and $\lambda$ has  finite $q$-th moment, for some $q > dp/(d-p)$;
  \item[(ii)] or $\lambda$ is absolutely continuous with respect to the Lebesgue measure on a bounded connected $C^2$-smooth or convex domain, with H\"older continuous density uniformly bounded from above and below by a strictly positive constant.
 \end{itemize}
 Then denoting by $\lambda_a$ the absolutely continuous part of $\lambda$ with respect to the Lebesgue measure,
 \begin{equation}\begin{split}\label{eq:dereich}
   \underline{\mathsf{c}}(p,d) \int_{\R^d} \lambda_a^{1-p/d} & \le \liminf_{n \to \infty} \EE\sqa{ \WN^p\bra{ \frac 1 n \sum_{i=1}^n\delta_{X_i} , \lambda}}n^{p/d} \\
   & \le \limsup_{n \to \infty} \EE\sqa{ \WN^p\bra{ \frac 1 n \sum_{i=1}^n\delta_{X_i} , \lambda}}n^{p/d}\le \overline{\mathsf{c}}(p,d)\int_{\R^d} \lambda_a^{1-p/d},
   \end{split}
 \end{equation}
where $(X_i)_{i=1}^n$ are i.i.d.\ random variables with common distribution $\lambda$. Moreover, the constants $0<\underline{\mathsf{c}}(p,d) \le \overline{\mathsf{c}}(p,d)<\infty$ are given by
\begin{equation*}\begin{split}
\overline{\mathsf{c}}(p,d)   = \lim_{L \to \infty}  \frac{1}{L^d} & \EE\sqa{ \WN_{(0,L)^d}^p\bra{ \mu, \frac{\mu((0,L)^d)}{L^d} }} \\ &  \text{and} \quad  \underline{\mathsf{c}}(p,d)= \lim_{L \to \infty} \frac{1}{L^d} \EE\sqa{ \WD_{(0,L)^d}^p\bra{ \mu, \frac{ \mu((0,L)^d)}{L^d} }},
\end{split}
\end{equation*}
where $\mu$ is a Poisson point process with unit intensity on $\R^d$, $\WN_{(0,L)^d}^p$ denotes the Wasserstein cost of order $p$ between the measures restricted on the cube $(0,L)^d$, while $\WD_{(0,L)^d}^p$ denotes instead its ``boundary'' variant (see \Cref{sec:boundary} for precise definitions).
\end{theorem}

 The boundary transport problem $\WD_{(0,L)^d}^p$ was first introduced in the optimal transport literature in \parencite{figalli2010new} to study the heat semigroup with Dirichlet boundary conditions as a gradient flow of the entropy. In the context of random Euclidean bipartite matching,  it is proposed in \parencite{barthe2013combinatorial} as a natural bipartite counterpart of the ``dual'' boundary functionals studied in \parencite{steele1997probability, yukich2006probability}.

Let us point out that under the assumptions (i), the statement of \Cref{thm:dereich} does not appear  in this precise form in the literature.  It can be however  easily obtained by combining the proofs of  \parencite{dereich2013constructive,barthe2013combinatorial}. Indeed, while \parencite{barthe2013combinatorial} treats the case of the random Euclidean assignment problem (or bipartite matching), the proofs extend straightforwardly to the case of the matching to the reference measure. By \parencite[Theorem 5 \& Lemma 4]{barthe2013combinatorial} this yields the upper bound in \Cref{eq:dereich} for measures with bounded support. Arguing then as in \parencite[Theorem 4]{dereich2013constructive} one can conclude the proof of the upper bound in \Cref{eq:dereich} for measures with finite $q-$moment. Notice that  this is a stronger statement than \parencite[Theorem 6]{barthe2013combinatorial}, where more integrability is required, or \parencite[Theorem 4]{dereich2013constructive}, where a Riemann integrability condition is imposed on $\lambda_a$. Regarding the lower bound in \Cref{eq:dereich}, it can be proven exactly as in \parencite[Theorem 8]{barthe2013combinatorial}, avoiding once more the Riemann integrability issue from \parencite[Theorem 4]{dereich2013constructive}.

In the statement  of  \Cref{thm:dereich}, the cases $d \in \cur{1,2}$ are excluded because of a different  asymptotic scaling, see \parencite{BoLe16,caracciolo2014scaling, ambrosio2019pde,ambrosio2022quadratic}.  In the range $p\ge d/2$, the upper bound in \Cref{eq:dereich} cannot hold without structural  assumptions on $\lambda$, see e.g. \parencite{FoGu15}. See however \parencite{pieroni-gaussian} for some results for Gaussian or Gaussian-like densities (as well as a  more general lower bound).
%




\subsection{Main result}

The natural question raised by \Cref{thm:dereich} is the existence of the limit in \Cref{eq:dereich}. It would be of course  a consequence of the equality $\underline{\mathsf{c}}(p,d)= \overline{\mathsf{c}}(p,d)$. This  is
 explicitly conjectured in \parencite[Remark 2]{dereich2013constructive},  albeit for a smaller constant, see \Cref{sec:open} -- but also in the bipartite case in \parencite{barthe2013combinatorial}. In this paper we settle this conjecture in the case $p=2$.
 \begin{theorem}\label{thm:main-poisson}
 For every $d\ge 3$, it holds
 \begin{equation}\label{eq:limit-poisson-intro}
  \underline{\mathsf{c}}(2,d)= \overline{\mathsf{c}}(2,d).
 \end{equation}
\end{theorem}

So far the only results in this direction were obtained in the case $d=p=2$ \parencite{caracciolo2014scaling, ambrosio2019pde,  benedetto2020euclidean, ambrosio2022quadratic}, where the constant $\overline{\costm}(2,2)  = 1/(4\pi)$ can be explicitly computed, and the concave one-dimensional case $d=1$, $p \in (0,1/2)$, \parencite{trevisan2022concave} where the special structure of the solution on the line can be exploited.

\begin{remark}
This problem is more subtle than it may first appear: for example, in \parencite{dobric1995asymptotics}, it is claimed that the $\liminf$ and $\limsup$ in \eqref{eq:dereich} always coincide and the limit equals the right hand side, for $p=1$ in any dimension $d \ge 3$. However, an error is contained in \parencite[Lemma 4.2]{dobric1995asymptotics} which cannot be fixed, also according to \parencite{barthe2013combinatorial}. Moreover, the validity of \Cref{eq:limit-poisson-intro} strongly relies on the stochastic properties of the Poisson point process. Indeed, the analog statement would not be true for the random grid, i.e.\ replacing $\mu$ by $\sum_{x\in \mathbb{Z}^d+Z}\delta_x$ where $Z$ is a uniformly distributed random variable on $(0,1)^d$. In light of  \Cref{ass:mu}, we see that \Cref{eq:limit-poisson-intro} fails for the random grid because it is not hyperuniform when tested with  with cubes (but it is when tested w.r.t.\ balls).
 \end{remark}

We also point out that if one is only interested in the equality between the $\liminf$ and the $\limsup$ in \Cref{eq:dereich}, when $\lambda_a$ is constant (and possibly $\lambda$ has a singular part with respect to the Lebesgue measure) it is still possible to conclude without knowing that $\underline{\mathsf{c}}(p,d)= \overline{\mathsf{c}}(p,d)$ by  \parencite[Theorem 2]{barthe2013combinatorial}. See also \parencite{goldman2022optimal}.

\subsection{Comments on the proof technique}\label{sec:pde-ansatz}

The proof of \Cref{thm:main-poisson} relies on novel tools from the stability theory of optimal transportation, combined with functional analytic tools and some ideas from quantitative stochastic homogenization. Since some fundamental features may be obscured by technical aspects, in this section we provide an overview of the argument in the simplified setting provided by the PDE ansatz from \parencite{caracciolo2014scaling}, which has been a stimulus for much of the recent progress in the field.

\subsubsection*{Deterministic bound} The starting point of the heuristics consists in linearizing the Monge-Amp\`ere equation for the optimal transport map, into a Poisson equation, essentially replacing the Wasserstein distance with a Sobolev norm of negative order. Given two measures $\mu$, $\lambda$ on a bounded domain $\Omega \subseteq \R^d$ (e.g., a cube), we approximate (writing here $\approx$ for a heuristic equivalence, $\les$ for an inequality)
\begin{equation*}
 \WN^2_{\Omega} (\mu, \lambda) \approx \int_{\Omega} |\nabla f|^2,
\end{equation*}
where $f$ solves the following Poisson equation, with null Neumann boundary conditions,
\begin{equation*}
 \begin{cases} \Delta f = \mu -\lambda & \text{in $\Omega$,}\\
  \nabla f \cdot \nu_{\Omega} = 0 & \text{on $\partial \Omega$,}
 \end{cases}
\end{equation*}
where $\nu_{\Omega}$ denotes the normal to the boundary of $\Omega$. A similar ansatz applied to $\WD^2_{\Omega}$ yields the approximation
\begin{equation*}
 \WN^2_{\Omega} (\mu, \lambda) \approx \int_{\Omega} |\nabla u|^2,
\end{equation*}
where $u$ solves the Poisson equation with same datum, but null Dirichlet boundary conditions:
\begin{equation}\label{eq:poisson-dirichlet}
 \begin{cases} \Delta u = \mu -\lambda & \text{in $\Omega$,}\\
  u  = 0 & \text{on $\partial \Omega$.}
 \end{cases}
\end{equation}
Despite being a rather drastic approximation, the main problem remains: how to compare the ``Dirichlet'' energy, $\int_{\Omega} |\nabla u|^2$, with the ``Neumann'' one, $\int_{\Omega} |\nabla f|^2$? To this aim, we first recall the following  representation (akin to the Benamou-Brenier formula) for the Neumann energy:
\begin{equation}\label{dualH1}
 \int_{\Omega} |\nabla f|^2 = \min\cur{ \int_{\Omega}|b|^2 \, : \, \div b = \mu -\lambda, \, \,  b\cdot \nu_{\Omega} = 0}.
\end{equation}
Similarly, for the Dirichlet energy we have:
\begin{equation}\label{dualH1Dir}
 \int_{\Omega} |\nabla u|^2 = \min\cur{ \int_{\Omega}|b|^2 \, : \, \div b = \mu -\lambda }.
\end{equation}
In particular, since $b=\nabla f$ is an admissible vector field for \Cref{dualH1Dir},
\begin{equation*}
 \int_{\Omega} |\nabla u|^2 \le \int_{\Omega} |\nabla f|^2.
\end{equation*}
This reflects a similar property for the Wasserstein distance. Hence one only needs to provide a converse inequality, up to an error term that is negligible in the limit when the difference $\mu-\lambda$ is small, in a suitable weak sense. \\
We see immediately that if $u$ solves \Cref{eq:poisson-dirichlet}, to obtain a competitor for \Cref{dualH1}, the gradient vector field $b_0:=\nabla u$,  only needs to be corrected for its flux at the boundary. We then introduce a smooth cut-off function $\eta: \Omega \to [0,1]$, that is identically $1$ on a neighborhood of $\partial \Omega$ and such that $\eta(x) = 0$ for any $x\in \Omega$ whose distance from $\partial \Omega$ is larger than a  parameter $r$ (to be chosen suitably small). 
In particular, we may assume that  $|\nabla \eta|\les r^{-1}$ and $|\nabla^2 \eta| \les r^{-2}$.
Setting $\xi:=1-\eta$, the corrected vector field  $b_1 := b_0 - \eta b_0 = \xi b_0$ satisfies the no-flux condition on $\partial \Omega$ -- it is actually null on $\partial \Omega$ -- but its divergence reads
\begin{equation*}
 \div b_1 -(\mu-\lambda)=  -\eta \div b_0 - \nabla \eta \cdot b_0 = -\eta (\mu-\lambda) - \nabla \eta \cdot \nabla u,
\end{equation*}
hence one needs  to introduce a further correction. At this stage, we observe the following elementary identity:
\begin{equation}\label{eq:magic-identity-intro}
 \nabla \eta \cdot \nabla u = \div(u \nabla \eta) - u \Delta \eta.
\end{equation}
This is a key observation for our analysis as it allows, via integration by parts,   to bound weak norms of our error terms by $u$  rather than $\nabla u$, see in particular \Cref{flatu} below.  Moreover, this suggests as a further correction, to use $b_2:= b_1+ u \nabla \eta$, which still has no flux at the boundary and
\begin{equation*}
 \div b_2 - (\mu-\lambda)= -\eta (\mu-\lambda) - u \Delta \xi.
\end{equation*}
To complete our construction, we solve the Poisson equation
\begin{equation}\label{eqcorrg}
\begin{cases}
 \Delta g = \eta (\mu-\lambda) + u \Delta \xi & \text{in $\Omega$}\\
 \nabla g \cdot \nu_{\Omega} = 0 & \text{on $\partial \Omega$,}
 \end{cases}
\end{equation}
and set $b_3:= b_2+\nabla g$.  Notice that using integration by parts, \Cref{eq:poisson-dirichlet} and the fact that $\xi$ vanishes on a neighborhood of $\partial \Omega$, we find
\begin{equation*}
 \int_{\Omega} \eta (\mu-\lambda)+u \Delta \xi=\int_{\Omega} \eta (\mu-\lambda)+ \xi(\mu-\lambda)=0
\end{equation*}
so that \Cref{eqcorrg} is indeed solvable.
Since $b_3$ is admissible for \Cref{dualH1}, we  have
\begin{equation*}
 \int_{\Omega}|\nabla f|^2  \le \int_{\Omega}|b_3|^2.
\end{equation*}
Using the triangle inequality  and standard energy estimates for elliptic equations we  eventually find that for any $\eps\in (0,1)$, 
\begin{equation}\begin{split}\label{eq:main-pde-simplified}
 \int_{\Omega}|\nabla f|^2  - \int_{\Omega}|\nabla u|^2 & \les \eps  \int_{\Omega}|\nabla u|^2\\
 & \quad + \frac{1}{\eps}\sqa{  \nor{ \eta(\mu - \lambda)}_{H^{-1,2}(\Omega)}^2 + \bra{r^{-2}+ \diam(\Omega)^2r^{-4}}\int_{\Omega}|u|^2 }.
 \end{split}
\end{equation}
 This is the PDE counterpart of our main deterministic bound for Wasserstein distances, \Cref{eq:main-deterministic} in \Cref{thm:main-deterministic}. Indeed, the first and last terms in the square brackets have a formal correspondence, by replacing the negative Sobolev norm with the Wasserstein distance and the solution to the Poisson problem  with the Kantorovich potential.

\subsubsection*{Application to the random setting}
Inequality \Cref{eq:main-pde-simplified} however becomes useful only if one can show that the terms in square brackets become negligible, which is indeed the case in our random setting. Let us assume therefore that $\Omega = (0,L)^d$, $\mu$  is a Poisson point process and $\lambda$ is the Lebesgue measure renormalized to have the same total mass of $\mu$ on $(0,L)^d$. Inserting the known upper (and lower) bounds for the matching problem in the PDE heuristics, we expect
\begin{equation}\label{PDEheuristics}
\frac 1 {L^d} \EE\sqa{\int_{(0,L)^d}|\nabla u|^2} \approx \frac 1 {L^d} \EE\sqa{\WD_{(0,L)^d}^2(\mu, \lambda) } \approx 1,
\end{equation}
which can be roughly summarized as $|\nabla u| \approx 1$, corresponding to fact that the transport happens typically on a distance of order $\approx 1$. 

Regarding the cut-off function $\eta$,  we choose $r = \delta L$ for some given $\delta\in (0,1)$. Then, we clearly see that a competition arises between the terms in the square brackets in \Cref{eq:main-pde-simplified}. The first term should become negligible when $\delta \to 0$, as seen again by the heuristics,
\begin{equation*}
 \frac 1 {L^d} \nor{   \eta(\mu - \lambda)}_{H^{-1,2}( (0,L)^d )}^2 \approx  \frac 1 {L^d} \WN^2_{(0,L)^d}(\eta \mu, \eta \lambda) \approx \delta.
\end{equation*}
Indeed, we expect that also the transport between the weighted measures happens at a distance of order $\approx 1$, but the volume of the support is now $\approx \delta L^d$. On the other side, the second term evidently grows as $\delta \to 0$:
\begin{equation*}
 \frac 1 {L^d} \bra{r^{-2}+ L^2r^{-4}}\int_{(0,L)^d}|u|^2 \approx \bra{ \delta^{-2} + \delta^{-4} } \cdot \frac{1}{L^{d+2}} \int_{(0,L)^d}|u|^2.
\end{equation*}
We notice  that Poincaré inequality on $(0,L)^d$ yields the upper bound
\begin{equation}\label{Ponctrivial}
 \frac{1}{L^{d+2}} \int_{(0,L)^d}|u|^2 \les \frac{L^2}{L^{d+2}} \int_{(0,L)^d} |\nabla u|^2 \approx 1,
\end{equation}
which is however not sufficient to conclude. Instead, we are able to argue that
\begin{equation}\label{flatu}
\EE\sqa{ \int_{\Omega}|u|^2} \ll L^{d+2}
\end{equation}
where $\ll$ means that the left-hand side is asymptotically of lower order. Thus, we may let first $L \to \infty$ and then $\delta \to 0$ to conclude. The intuition behind \Cref{flatu}  is that $\nabla u$ is highly oscillating, hence Poincaré inequality fails to correctly capture cancellations and  thus the order of $u$. A similar phenomenon happens in quantitative stochastic homogenization, see for instance \parencite{armstrong2019quantitative,armstrong2016quantitative}, and in the PDE heuristics one sees how to borrow ideas from that field. Namely, we introduce an intermediate scale $1\ll L_0 \ll L$ and partition $(0,L)^d$ into $m \approx (L/L_0)^d$ cubes, each of side length $L_0$, and (roughly) define a competitor $\tilde u$ on $(0,L)^d$ for the Dirichlet energy by gluing the solutions to \Cref{eq:poisson-dirichlet} on each sub-cube. Decomposing the average expected energy into the sum of  the contributions on each sub-cube, by stationarity of the Poisson point process, we find that it is  equal to the expected average energy of the solution $u_0$ to the Poisson problem \Cref{eq:poisson-dirichlet} on a single cube $(0,L_0)^d$:
\begin{equation}\label{tildeuuzero}
 \frac 1 {L^d} \EE\sqa{ \int_{(0,L)^d} |\nabla \tilde u|^2} =  \frac{1}{L_0^d} \EE\sqa{ \int_{(0,L_0)^d} |\nabla  u_0|^2}.
\end{equation}
Using  super-additivity arguments, it can be shown that the right-hand side of \Cref{tildeuuzero} converges as $L_0\to \infty$. In particular for $L_0$ large (which also implies $L$ large) it must be almost equal to the left-hand side of \Cref{PDEheuristics}.  This implies that the vector field $\nabla \tilde{u}$ is almost minimizing \Cref{dualH1Dir}, i.e.,
\begin{equation*}
 \frac 1 {L^d} \EE\sqa{  \int_{(0,L)^d} |\nabla \tilde u|^2} -   \frac 1 {L^d}\EE\sqa{  \int_{(0,L)^d} |\nabla u|^2} =  \frac{1}{L_0^d} \EE\sqa{ \int_{(0,L_0)^d} |\nabla  u_0|^2} -  \frac 1 {L^d}\EE\sqa{  \int_{(0,L)^d} |\nabla u|^2} \ll 1.
\end{equation*}
By strong convexity of the Dirichlet energy, we have the  stability (in)equality
\begin{equation}\label{eq:stability}
  \EE\sqa{ \int_{(0,L)^d} |\nabla u- \nabla \tilde u|^2 }=\EE\sqa{  \int_{(0,L)^d} |\nabla \tilde u|^2} -  \EE\sqa{  \int_{(0,L)^d} |\nabla u|^2} \ll L^d.
\end{equation}
Coupling this with  Poincaré inequality and $\EE\sqa{ \int_{(0,L)^d} |\tilde u|^2} \les L_0^2 L^d$, by \eqref{Ponctrivial} applied in each sub-cube of side length $L_0$,  we obtain
\begin{equation*}\begin{split}
 \EE\sqa{ \int_{(0,L)^d} |u|^2} & \les \EE\sqa{ \int_{(0,L)^d} |\tilde u|^2} +\EE\sqa{  \int_{(0,L)^d} |u-\tilde u|^2} \\
 & \les L_0^2 L^d + L^2 \EE\sqa{\int_{(0,L)^d} |\nabla u- \nabla \tilde u|^2 } \ll L^{d+2},
 \end{split}
\end{equation*}
as claimed in \Cref{flatu}.

\subsubsection*{Further comments} The argument above fairly describes the main structure of the proofs of our results.  However, we need to take into account that the Wasserstein distance is not actually approximated, at least up to length scales  of order $\approx 1$, by the negative Sobolev norm, hence much of technical effort is to fit these ideas, that may appear simple in the PDE heuristics, into the actual optimal transport setting. A crucial point we highlight here is for example the need for an analogue of the stability inequality \Cref{eq:stability}, which in the PDE setting uses strong convexity of the energy, while the optimal transport cost is usually only weakly convex -- it is a linear programming problem. For this, we rely in \Cref{thm:delalande}, on some recent results in the literature of stability for optimal transport problem for the standard quadratic cost \parencite{delalande2022quantitative, delalande2023quantitative} together with a localization argument.  Let us point out that in our application, the uniform monotonicity of the gradient of the Kantorovich functional as stated in \parencite{delalande2023quantitative,letrouit2024gluing,mischler2024,chizat2024sharper} is not enough and we need the seemingly stronger uniform convexity of the Kantorovich functional itself as presented in \parencite{delalande2022quantitative} (following the proof from \parencite{delalande2023quantitative}).

Let us point out that contrary to \Cref{eq:stability} which gives directly control on the gradients, the stability inequality \Cref{eq:taylor-entropic-delalande-mu-vt-main-body} in the case of the transport cost gives only control on the potentials. This forces us to rely on the div-curl type identity \Cref{eq:magic-identity-intro} which acts as a replacement for the Caccioppoli inequality.  Indeed, in the PDE heuristics it would be actually possible to avoid \Cref{eq:magic-identity-intro} thanks to this stronger stability property. As will be apparent in the proofs, see in particular \Cref{eq:main-deterministic} in \Cref{thm:main-deterministic}, it seems that the lack of \Cref{eq:magic-identity-intro} is the main obstacle to the extension of our results from $p=2$ to $p>1$ arbitrary.

Finally, we notice that our arguments apply to general stationary random measures that satisfy suitable concentration bounds. Precise  assumptions are given in  \Cref{ass:mu}. We thus actually show \Cref{thm:main-general}, which contains \Cref{thm:main-poisson} as a particular case. Another notable example covered by \Cref{thm:main-general}, is the Brownian interlacement occupation measure -- see \parencite{mariani2023wasserstein, drewitz2014introduction, sznitman2013scaling} in any dimension $d \ge 5$.

\subsection{Further questions and conjectures}\label{sec:open} Our work makes substantial progress towards the full resolution of the conjectured validity of the limit in \Cref{eq:dereich}. We list here some questions and open problems that may be further addressed:
\begin{enumerate}[i)]

 \item A result for the random Euclidean bipartite matching problem, thus settling the conjecture from \parencite[Theorem 2]{barthe2013combinatorial}, seems also quite natural, and the only missing ingredient seems to be a stability inequality akin to \Cref{thm:delalande}, when both measures are singular -- only the case when one of the two is sufficiently regular can be treated at the moment.

 \item In \parencite{dereich2013constructive}, the constant $\underline{\costm}(d,p)$ is defined in terms of an asymmetric boundary Wasserstein distance, which is smaller than the one we consider here. It would be interesting to prove that both constants coincide. In particular it would also imply that they coincide with the constant from  \parencite{huesmann2013optimal}. In light of \parencite[Theorem 1.1]{huesmann2013optimal} this would be a major step towards proving that all these problems generate the same  coupling between the Lebesgue measure and the Poisson point process in the limit $L\to \infty$.

 \item While we drew inspiration from quantitative stochastic homogenization, we only proved here  a qualitative convergence result. As its name suggests, the main and far more ambitious goal of quantitative stochastic homogenization is to obtain quantitative (and hopefully optimal) convergence rate. It would be interesting to investigate if such bounds could be obtained also for the transport problem.  It is in fact conjectured in \parencite{CaLuPaSi14} that the rate of convergence should be of the order of $L^{2-d}$, see also \parencite{goldman2023almost} for the current rate of convergence known in the case $p=d=2$.

 \item As already mentioned, our convergence result applies as well to the optimal transport of the Brownian interlacement occupation measure, which in \parencite{mariani2023wasserstein} has been shown to  provide an upper bound to the rate of convergence of the occupation measure of a single Brownian path $(B_t)_{t \ge 0}$ in the flat torus $\mathbb{T}^d = \R^d/\mathbb{Z}^d$. We conjecture that the boundary variant may provide also an asymptotic lower bound, hence settling a question left open in \parencite{mariani2023wasserstein}. We remark that by using the decomposition arguments from \parencite[Section 6]{ambrosio2022quadratic} our main result yields that similar limits exist in the case of transport of i.i.d.\ points on compact manifolds. Asymptotic upper and lower bounds for the occupation measure of general diffusion processes have been studied recently, see e.g.\ \parencite{wang2021precise, wang2020convergence, wang2022wasserstein, wang2021wasserstein, wang2023limit, wang2021convergence, wang2024sharp}, although to our knowledge limit results are available only in lower dimensions -- corresponding to the cases where the scaling behaviors of the costs are somehow exceptional.
\end{enumerate}

\subsection{Structure of the paper} The exposition is structured as follows. \Cref{sec:notation} provides the necessary background on general notation, measure, metric and Sobolev spaces. In \Cref{sec:deterministic}, after introducing the optimal transport theory for distance costs and its boundary variant, we establish the main deterministic results, in particular \Cref{thm:main-deterministic}. \Cref{sec:random} follows by providing \Cref{thm:main-general}, which is a generalized version of \Cref{thm:main-poisson} -- we also show that the Brownian interlacement occupation measure satisfies its assumptions.

\section{Notation and basic facts}\label{sec:notation}

Given a set $A$, we write $\chi_A$ for its indicator function. We write $x\cdot y$ for the standard scalar product in $\R^d$ and $\abs{x} = \sqrt{ x \cdot x}$ for the Euclidean norm. We always endow $\R^d$ with the Euclidean distance. We introduce the following notation for the ``power'' $p>0$ of a vector $x\in \R^d$, given by $x^{(p)}:= |x|^{p-1} x$, with the convention that $0^{(p)} = 0$. The notation $A\les B$ means that there exists a
constant $C>0$,  such that $A\le C B$, where $C$ depends on the dimension $d$ and  $p$.  We use the notation $\les_q$ to indicate the dependence on the parameter $q$. We write $A\sim B$ if both $A\les B$ and $B\les A$.  

\subsection{Measure spaces}


Given a finite measure space  $(E, \mathcal{E}, \lambda)$ and a function $f$ on $E$, taking values in $\R^d$, we write for $p \in [1, \infty)$, $\nor{f}_{L^p(\lambda)} := \bra{ \int_{E}|f|^p  d \lambda}^{1/p}$  for the Lebesgue norm, 
and  $\nor{f}_\infty := \sup_{x\in E} |f(x)|$ for the uniform norm. For a real valued $f \in L^2(\lambda)$, we write
\begin{equation*}
  \Var_\lambda (f) := \inf_{a \in \R} \int_{E}\abs{ f - a}^2 d\lambda, 
\end{equation*}
which is attained at $a = \fint_E f d\lambda:= \lambda(E)^{-1} \int_E f d \lambda$. Notice that  we do not require $\lambda$ to be a probability: however, it holds
\begin{equation}\label{eq:var-c-lambda}
  \Var_{c \lambda} (f)  = c \Var_{\lambda}(f) \quad \text{for every $f \in L^2(\lambda)$, $c >0$}.
\end{equation}
The function $f \mapsto \sqrt{ \Var_{\lambda}(f)}$ defines a semi-norm on $L^2(\lambda)$ that is null on constant functions. 
If $\tilde \lambda$ is a measure on $E$ such that $\tilde \lambda \le\lambda$, then it holds
\begin{equation}\label{eq:var-tilde-lambda}
  \Var_{\tilde \lambda}(f)  \le \Var_{\lambda} (f) \quad \text{for every $f \in L^2(\lambda)$.}
\end{equation}
%
%
%
When $\lambda = P$ is a probability measure, we employ the standard probabilistic notation, e.g.\ $\EE\sqa{\cdot} = \int_{E} \cdot dP$ for the expectation value and $\Var(\cdot)=\Var_{P}(\cdot)$ for the variance.

When $E \subseteq \R^d$ is measurable and bounded, and $\lambda$ is the Lebesgue measure on $E$, we  simply write $\nor{f}_{L^p(E)} :=\nor{f}_{L^p(\lambda)}$, $|E| := \lambda(E)$, $\int_E f := \int_E f(x) dx$ for $f \in L^1(E)$ and $\fint_Ef:= |E|^{-1} \int_E f$.  We also write
\begin{equation*}
 \Var_E (f) := \Var_{\lambda}(f) \quad \text{for every $f \in L^2(E)$,}
\end{equation*}
and notice that, for a measurable $\tilde E \subseteq E$ it holds
\begin{equation*}
 \Var_{\tilde E} (f) \le \Var_{E}(f) \quad \text{for every $f \in L^2(E)$.}
\end{equation*}
For simplicity, we always identify any finite measure $\tilde \lambda$ on $E$ that is absolutely continuous with respect to the Lebesgue measure with its density that we keep denoting with the same letter. For example, we write $\int_E f \lambda := \int_E f d\lambda$.

\subsection{Metric spaces} Given a  metric space $(\cX, \dist)$  and $\Omega \subseteq \cX$, we write $\overline{\Omega}$ for its closure, 
\begin{equation*}
 \dist (\Omega) =\sup_{x,x' \in \Omega} \dist(x,x')
\end{equation*}
for its diameter (also written $\diam(\Omega)$ if $\dist$ is understood). We  write
\begin{equation*}
  \dist(x, \Omega^c) := \inf_{z \in \Omega^c} \dist(x,z),
\end{equation*}
for the distance function from $\Omega^c$ and for $\cY \subseteq \cX$,
\begin{equation*}
 \dist(\Omega, \cY) = \sup_{x \in \Omega} \dist(x, \cY).
\end{equation*}
Notice that this quantity is not symmetric. It will be convenient to introduce the notation
\begin{equation}\label{defOmr}
 \Omega^{(r)} := \cur{ x \in \Omega\, : \, \dist(x, \Omega^c) > r}, \quad \text{for $r \ge 0$.}
\end{equation}
For a  function $f: \cX \to \R$, we write $\sup_{\Omega} f :=\sup_{x \in \Omega} f(x)$, $\inf_{\Omega}f := \inf_{x\in \Omega} f(x)$, $\osc_{\Omega}\bra{f} := \sup_{\Omega} f - \inf_\Omega f$ for its oscillation and
\begin{equation*}
 \Lip_{\Omega}(f) := \sup_{x\neq x' \in \Omega} \frac{|f(x)-f(x')|}{|x-x'|}
\end{equation*}
for its Lipschitz constant.  We have
\begin{equation*}
 \osc_{\Omega}(f) \le \diam(\Omega) \Lip_{\Omega}(f),
\end{equation*}
and, for every finite measure $\lambda$ supported on $\Omega$,
\begin{equation}\label{eq:var-osc}
 \Var_{\lambda} (f) \le  \frac{\lambda(\Omega)}{4} \osc_{\Omega}(f)^2.
\end{equation}

When $\cX \subseteq  \R^d$ (endowed with the Euclidean distance) we often consider $\Omega \subseteq \cX$ to be a bounded domain with Lipschitz boundary. In our main results we require -- but always specify -- that $\Omega$ is convex. Examples include the cases of a rectangle $R = \prod_{i=1}^d(a_i, b_i)$ with $a_i < b_i$ for $i=1, \ldots, d$ or a cube $Q$, i.e.\ the side lengths $b_i-a_i$ for $i=1, \ldots, d$ all coincide. In particular, for $\Omega \subseteq \R^d$,  \Cref{eq:var-osc} reads
\begin{equation}\label{eq:variance-osc}
 \Var_{\Omega} (f) \le \frac{|\Omega|}{4} \osc_{\Omega}(f)^2.
\end{equation}

%
%
%

%


\subsection{Sobolev spaces}

Given a bounded connected domain $\Omega \subseteq \R^d$ with Lipschitz boundary, the gradient of a  function $f: \Omega \to \R$ is denoted by $\nabla f$, the divergence of a vector field $b: \Omega \to \R^d$ is  denoted by $\div b$, and the outward unit normal to the boundary of $\Omega \subseteq\R^d$ is denoted by $\nu_{\Omega}$. We recall that a Lipschitz function is differentiable at Lebesgue a.e.\ $x \in \Omega$ and that
\begin{equation*}
 \nor{\nabla f}_{L^\infty(\Omega)} \le  \Lip_{\Omega}(f),
\end{equation*}
with equality if $\Omega$ is convex.


Given a function $\eta : \Omega \to [0, \infty)$, which we assume for simplicity bounded, weighted Sobolev spaces $H^{1,p}(\eta)$ are defined as the completion of smooth (or Lipschitz) functions with respect to the Sobolev norm
\begin{equation*}
 \nor{f}_{H^{1,p}(\eta)}:= \bra{ \nor{f}_{L^p(\eta)}^p + \nor{\nabla f}_{L^p(\eta)}^p}^{1/p}.
\end{equation*}
Since convergence in $H^{1,p}(\eta)$ entails convergence in $L^p(\eta)$ it follows that $H^{1,p}(\eta) \subseteq L^p(\eta)$. The subspace $H^{1,p}_0(\eta)$ of functions that have zero trace at the boundary is defined as the closure of smooth functions compactly supported in $\Omega$. When, $\eta = \chi_{\Omega}$, we write $H^{1,p}(\Omega) := H^{1,p}(\eta)$ and similarly $H^{1,p}_0(\Omega) := H^{1,p}_0(\eta)$.


A key tool in Sobolev spaces analysis is the validity of a Poincaré-Wirtinger inequality, which in its weighted form reads as follows: there exists a constant $c<\infty$ such that for every Lipschitz function $f$ and null weighted mean $\int_{\Omega} f \eta = 0$,
\begin{equation}
 \label{eq:weighted-poincare}
  \nor{ f    }_{L^p(\eta)} \le c\nor{ \nabla f}_{L^p(\eta)}.
\end{equation}
 We write throughout $c_P(\eta, p)$ for the smallest constant such that \Cref{eq:weighted-poincare} holds. Let us notice that, if $\eta'$ is a weight equivalent to $\eta$ in the sense that, for some constant $c_0 \in (0,\infty)$, $c^{-1}_0 \eta(x) \le \eta' (x) \le c_0\eta(x)$ for every $x \in \Omega$, then it holds
\begin{equation}\label{comparePoincare}
  \frac 1 {2 c_0^{1/p}} c_P(\eta,p)  \le c_P(\eta', p) \le 2 c^{1/p}_0 c_P(\eta).
\end{equation}
Indeed, for $p=2$ this follows (even without the factors $1/2$ and $2$) from the fact that the left-hand side in \Cref{eq:weighted-poincare} can be equivalently rewritten as $\Var_{\eta}(f)$, dropping the null mean requirement, and the homogeneity properties \Cref{eq:var-c-lambda} \Cref{eq:var-tilde-lambda} of the variance easily allow for a change of measure. For general $p$, one first notices that  $c_P(\eta, p)<\infty$ is actually equivalent to the validity of an inequality of the type
\begin{equation}\label{eq:poincare-equivalent}
 \inf_{a \in \R} \nor{ f  -a  }_{L^p(\eta)} \le \tilde c\nor{ \nabla f}_{L^p(\eta)} \quad \text{for every $f \in H^{1,p}(\eta)$,}
\end{equation}
for some constant $\tilde c \in (c_P(\eta, p)/2, c_P(\eta,p))$, and then argues as in the case of $p=2$ using similar properties for the ``$p$-variance'' appearing in the left-hand side. Indeed, \Cref{eq:weighted-poincare}  implies \Cref{eq:poincare-equivalent} with $\tilde c = c_P(\eta, p)$, by taking $a = \fint_{\Omega} f \eta$. Conversely, assuming  for simplicity that $\int_{\Omega} \eta = 1$, we notice that, for every $a \in \R$,
\begin{equation*}
  \abs{ \int_{\Omega} f \eta  - a} \le \nor{ f - a}_{L^1(\eta)} \le \nor{f-a}_{L^p(\eta)},
\end{equation*}
so that if  \Cref{eq:poincare-equivalent} holds we obtain, by the triangle inequality,
\begin{equation*}\begin{split}
 \nor{f- \int_{\Omega} f \eta }_{L^p(\eta)} & \le  \inf_{a \in \R} \cur{ \nor{ f  -a  }_{L^p(\eta)}  + \abs{ \int_{\Omega} f \eta  - a}} \le 2  \inf_{a \in \R} \nor{f-a}_{L^p(\eta)}\\
 & \le 2 c_P(\eta, p) \nor{\nabla f}_{L^p(\eta)}.
 \end{split}
\end{equation*}


In case of a constant weight $\eta = \chi_{\Omega}$, we write $c_P(\Omega,p):= c_P(\chi_\Omega,p)$ and we recall that since we assume that $\Omega \subseteq \R^d$ is a bounded connected domain with Lipschitz boundary, it always holds $c_P(\Omega,p ) <\infty$. In particular, in the case of $\Omega = Q_L$ a cube of side length $L$, we have $c_P(Q_L, p) \les L$ and more generally $c_P$ scales linearly with respect to dilations. In \parencite[Theorem 1.2]{ruiz} it is  also proved that, for any family of domains $\cur{\Omega_i}_{i \in I}$, the constants $c_P(\Omega_i, p)$ are uniformly bounded, provided that the following three conditions hold:
\begin{enumerate}[a)]
 \item the diameters $\cur{\diam(\Omega_i)}_{i \in I}$ are uniformly bounded,
 \item the sets $\cur{\Omega_i}_{i \in I}$ are uniformly connected, i.e., there exists $r>0$ such that all the sets $\Omega_i^{(r)}$ are connected,
 \item a uniform interior cone condition holds, i.e., there exists an open cone $\mathcal{C}$ (with finite height) such that, for every $i \in I$, $x \in \partial \Omega_i$, one can find a cone $\mathcal{C}_x$ with vertex at $x$ and congruent to $\mathcal{C}$ with $\mathcal C_x\subset \Omega_i$.
\end{enumerate}
 The uniform cone condition holds in particular if the boundaries $\cur{\partial \Omega_i}_{i\in I}$ are uniformly Lipschitz.

  Back to the weighted case, in \parencite[Theorem 1]{Rath}, it is proved that
\begin{equation} \label{eq:poincare-rath}
 c_P(\eta,p)\les \inf_{0 \le \tau \le \nor{\eta}_\infty } \bra{ \frac{ \nor{\eta}_\infty |\Omega| } {\tau |\cur{\eta >\tau}| }}^{1/p} \nor{ t\mapsto c_P( \cur{\eta>t}, p ) }_{L^\infty([0,\tau])},
\end{equation}
where the sets $\cur{\eta >t} = \cur{x \in \Omega\, : \, \eta(x) >t}$ are assumed to be connected domains with Lipschitz boundary.

\begin{remark}\label{rem:finite-weighted-poincare} As a simple application of this inequality and the uniformity discussed above, let us consider a weight  $\eta: \Omega \to [0,1]$ such that $|\cur{\eta>1/2}|>0$ and the sets $(\cur{\eta>t})_{t \in [0,1/2]}$ have uniformly Lipschitz boundaries, hence they satisfy a uniform cone property, and they are uniformly connected. Then, $\nor{ t\mapsto c_P( \cur{\eta>t}, p ) }_{L^\infty([0,1/2])}$ is finite and we obtain by \Cref{eq:poincare-rath} with $\tau=1/2$ that $c_P(\eta, p)< \infty$. The same argument extend to any family of weights $(\eta_i)_{i \in I}$ such that
\begin{equation*}
\sup_{i \in I} \nor{\eta_i}_\infty < \infty, \quad   \inf_{i \in I} |\cur{\eta_i>1/2}| >0,
\end{equation*}
and the sets $(\cur{\eta_i>t})_{t \in [0,1/2], i \in I}$   satisfy a uniform cone property and are uniformly connected. In such a case, we deduce
\begin{equation*}
 \sup_{i \in I} c_P(\eta_i, p)<\infty.
\end{equation*}
This may also follow by the continuity of the Poincar\'e constant with respect to Hausdorff convergence for sets satisfying a uniform cone condition, see  \parencite[Corollary 7.4.2]{bucur2000variational} for the case $p=2$.
\end{remark}

%


Sobolev spaces  of negative order are defined in terms of the following dual ``norm'':
\begin{equation}\label{defHmoinsp}
 \|f\|_{H^{-1,p}(\eta)}=\sup_{\nor{\nabla \phi}_{L^{p'}(\eta)}\le 1} \int_{\Omega} f \phi \eta \in [0, \infty]
\end{equation}
for $p \in (1, \infty)$ where $p' = p/(p-1)$. By definition, the supremum is upon $\phi \in H^{1,p'}(\eta)$ or equivalently  upon Lipschitz functions. We notice in particular that if $\|\cdot\|_{H^{-1,p}(\eta)}<\infty$ then $\int_\Omega f \eta=0$. In this case we may  restrict the supremum to functions $\phi$ having also average zero.




Using H\"older's inequality and \Cref{eq:weighted-poincare}, we obtain immediately that
 \begin{equation}\label{eq:poincare-negative}
  \|f\|_{H^{-1,p}(\eta)}\les c_P(\eta, p') \|f\|_{L^p(\eta)},
 \end{equation}
 provided that $\int_{\Omega}f \eta = 0$.
We define the operator
\begin{equation}\label{defdiveta}
 \div_\eta(b) = \div(b) - \eta^{-1} \nabla \eta \cdot b
\end{equation}
which is understood in the weak sense, i.e.,
\begin{equation}\label{divetaweak}
 \int_{\Omega} \nabla \phi \cdot b  \eta  = - \int_{\Omega} \div_\eta(b) \phi \eta \quad \text{for every Lipschitz function $\phi$.}
\end{equation}
It is then classical that provided the Poincaré inequality \Cref{eq:weighted-poincare} holds for $p'$, then  we have the dual representation
\begin{equation}\label{eq:dual-neg-sob}
 \nor{f }_{H^{-1,p}(\eta)} = \inf \cur{ \nor{ b}_{L^p(\eta)}\, : \, -\div_\eta(b) = f }.
\end{equation}
Indeed, inequality $\le$ is  true as a simple consequence of H\"older inequality. For the reverse inequality, thanks to \Cref{eq:weighted-poincare} we see that the supremum in the definition \Cref{defHmoinsp} is attained for some function $\phi$. Using a first variation argument we have $\phi=\lambda h$ where $\lambda=\|\nabla h\|_{L^{p'}(\eta)}^{-1}$ and $h$ solves (in the weak sense)
\begin{equation*}
 -\div_\eta (|\nabla h|^{p'-2}\nabla h)=f.
\end{equation*}
Using \Cref{divetaweak} we then see that setting $b=|\nabla h|^{p'-2}\nabla h$, we have $-\div_\eta b= f$ and
\begin{equation*}
 \|f\|_{H^{-1,p}(\eta)}=\|\nabla h\|_{L^{p'}(\eta)}^{p'-1}=\|b\|_{L^p(\eta)}
\end{equation*}
proving the other inequality.
%

As a consequence of \Cref{eq:dual-neg-sob}, we have the following sub-additivity property: given a (finite) partition $\Omega = \bigcup_{k}\Omega_k$ into domains such that $c_P(\chi_{\Omega_k} \eta, p')<\infty$ for every $k$, then it holds
\begin{equation}\label{eq:sub-add-sobolev-neg}
 \nor{f}_{H^{-1,p}(\eta) }^p \le (1+\eps) \sum_{k} \nor{ f - \int_{\Omega_k} f \eta_k}_{H^{-1,p}(\chi_{\Omega_k} \eta)}^p + \frac{C}{\eps^{p-1}} \nor{ \sum_k \chi_{\Omega_k}\int_{\Omega_k} f \eta_k }_{H^{-1,p}(\eta)}^p,
\end{equation}
 where $\eps \in (0,1)$ is arbitrary and $\eta_k:= \eta/\int_{\Omega_k}\eta$ (set equal to $0$ if $\eta$ is null on $\Omega_k$). To see this, we first use the triangle inequality for $\nor{\cdot}_{H^{-1,p}(\eta)}$,
 \begin{equation*}
  \nor{f}_{H^{-1,p}(\eta) }^p \le (1+\eps)   \nor{ f - \sum_{k}\chi_{\Omega_k}  \int_{\Omega_k} f \eta_k}_{H^{-1,p}(\eta)}^p + \frac{C}{\eps^{p-1}} \nor{ \sum_k \chi_{\Omega_k}\int_{\Omega_k} f \eta_k }_{H^{-1,p}(\eta)}^p.
 \end{equation*}
 For the first term we then use the formulation \Cref{eq:dual-neg-sob} with $b = \sum_{k} b_k \chi_{\Omega_k}$, where $b_k$ is admissible for  \Cref{eq:dual-neg-sob} with $f_k:=  \bra{ f  - \int_{\Omega_k} f \eta_k}\chi_{\Omega_k} \in L^p(\eta_k)$ and then minimize over $b_k$.

We end this section with two results that hold in the case of uniform weight $\eta = \chi_{\Omega}$. The first one provides a bound for the $H^{-1,2}$ norm of a function of the form $f = \nabla \eta \cdot \nabla u$ and is based on the identity \Cref{eq:magic-identity-intro} discussed in \Cref{sec:pde-ansatz}. Unfortunately this result seems to be very specific to the case $p=2$.

\begin{lemma}\label{lem:magic-identity}
 Let $\Omega \subseteq \R^d$ be a bounded domain with Lipschitz boundary, let $u \in H^{1,2}(\Omega)$  and $\eta: \Omega \to \R$ be $C^2$ smooth and such that $\nabla \eta \cdot \nu_{\Omega} = 0$ on $\partial \Omega$.
 Then, it holds
 \begin{equation}\label{eq:magic-bound}
  \nor{ \nabla \eta \cdot \nabla u - \fint_\Omega \nabla \eta \cdot \nabla u }_{H^{-1,2}(\Omega)}\les  \bra{  \nor{\nabla \eta}_{\infty}  + \diam(\Omega) \nor{\Delta \eta}_\infty}  \Var_{\Omega_\eta}(u)^{1/2}.
 \end{equation}
where $\Omega_\eta := \Omega \cap \supp |\nabla \eta|$. Here, the implicit constant depends on $\Omega$ only and is invariant under dilations.
\end{lemma}


\begin{proof}
 We use \Cref{eq:magic-identity-intro} and rewrite, for any smooth function $\phi$ with zero average on $\Omega$,
 \begin{equation*}\begin{split}
  \int_{\Omega} \phi \nabla \eta \cdot \nabla u & = \int_{\Omega } \phi \sqa{ \div(u \nabla \eta) - u \Delta \eta}\\
  & = -\int_{\Omega} u \nabla \phi  \cdot \nabla \eta  - \int_{\Omega}  \phi u \Delta \eta\\
  \end{split}
 \end{equation*}
 where in the second identity there is no boundary term because $\nabla \eta \cdot \nu_{\Omega} = 0$. We bound from above
 \begin{equation*}
   \int_{ \Omega} u \nabla \phi  \cdot \nabla \eta  \le \nor{\nabla \eta}_\infty \nor{u}_{L^2( \Omega_\eta)} \nor{\nabla \phi}_{L^2(\Omega)}
 \end{equation*}
 and, using Poincaré inequality \Cref{eq:weighted-poincare}, we find
 \begin{equation*}\begin{split}
  \int_{ \Omega}  \phi u \Delta \eta  & \le \nor{ \Delta \eta}_\infty \nor{ u}_{L^2( \Omega_\eta)} \nor {\phi}_{L^2(\Omega)} \\
  & \le   \diam(\Omega) \nor{ \Delta \eta}_\infty \nor{ u}_{L^2( \Omega_\eta)} \nor{\nabla \phi}_{L^2(\Omega)}.
  \end{split}
 \end{equation*}
 Summing these two bounds and taking the supremum over $\phi$ with $\nor{\nabla \phi}_{L^2(\Omega)} \le 1$ yields \Cref{eq:magic-bound} with $\nor{u}_{L^2( \Omega_\eta)}^2$ instead of $\Var_{\Omega_\eta}(u)$. To conclude it is sufficient to notice that the left-hand side remains unchanged if one replaces $u$ with $u-a$ for any constant $a \in \R$, and  minimize upon $a$.
\end{proof}


For the second result, which is  a particular (but not totally standard) case of Gagliardo-Nirenberg interpolation inequality, we recall first the case $p>d$ of Sobolev inequality, which reads
\begin{equation}\label{eq:sobolev}
  \sup_{x\neq x'\in \Omega} \frac{\abs{ f(x)-f(x')}}{|x-x'|^{1-d/p}} \les \nor{\nabla f}_{L^p(\Omega)} \quad \text{for every $f \in H^{1,p}_0(\Omega)$,}
\end{equation}
where the implicit constant depends on $p$ and $\Omega$, but is invariant with respect to dilations of $\Omega$.

\begin{lemma}\label{lem:gagliardo}
For every $p\ge 1$ and every Lipschitz function $f:\R^d \to \R$ that is constant on $\Omega^c$, it holds
\begin{equation*}
 \osc_{\Omega}(f) \les \diam(\Omega)^{\frac{p}{p+d}} \Lip_{\Omega}(f)^{\frac{d}{d+p}} \nor{\nabla f}_{L^p(\Omega)}^{\frac{p}{p+d}},
\end{equation*}
where the implicit constant depends on $p$  and $\Omega$ but is invariant with respect to dilations of $\Omega$.
\end{lemma}

\begin{proof}
 Without loss of generality, we can assume that $f(x) = 0$ if $x \in \Omega^c$. Moreover, by scaling we may also assume that $\diam(\Omega)=1$. Write $f_\eps$ for the convolution of $f$ with a standard convolution kernel and split
\begin{equation*}
\|f\|_\infty\les \|f-f_\eps\|_\infty +\|f_\eps\|_\infty.
\end{equation*}
We have
\begin{equation*}
 \|f-f_\eps\|_\infty\les \eps \|\nabla f\|_\infty  \le \eps \Lip_{\Omega}(f)
 \end{equation*}
 and
 \begin{equation*}
 \|f_\eps\|_\infty\les \eps^{-d/p} \nor{f}_{L^p(\R^d)} \le c_P(\Omega,p) \eps^{-d/p} \nor{\nabla f}_{L^p(\Omega)},
\end{equation*}
where the second inequality follows from Poincaré inequality in $H^{1,p}_0(\Omega)$ (one does not need to subtract the mean). An optimization over $\eps$ yields the claimed interpolation inequality.
\end{proof}

\section{Deterministic results}\label{sec:deterministic}

In this section, we collect all the deterministic tools related to optimal transport theory that we need in the proof of our main result. After recalling several basic facts that can be found in any of the standard references on the subject \parencite{AGS,VilOandN, peyre2019computational}, we proceed towards the proof of our main result for this section,  \Cref{thm:main-deterministic}. It provides a non-trivial upper bound for $\WN^2$ in terms of $\WD^2$ and suitable error terms. A key tool that we use is a  stability result for the quadratic optimal transport problem,  see \Cref{eq:stability}. Although we are only ultimately able to deal with the quadratic case, whenever the proofs allow, we provide more general results that readers may find useful.

\subsection{Generalities}
\label{sec:wass}
Let $\cX$ denote a compact metric  space and $\dist: \cX \times \cX \to [0, \infty)$ be a pseudo-distance, i.e.\ like a standard distance with only the exception that $\dist(x,y) = 0$ does not necessarily entail $x=y$. We also assume that $\dist$ is (lower-)semicontinuous with respect to the topology of $\cX$. Given two finite positive Borel measures $\mu$, $\lambda$ on $\cX$, with $\mu(\cX) = \lambda(\cX)$, the  optimal transport cost of order $p\ge 1$ between  $\mu$ and $\lambda$  is defined as the quantity 
\begin{equation*} \Wdist^p(\mu, \lambda) =  \inf_{\pi\in\mathcal{C}(\mu,\lambda)}  \int_{\cX \times \cX} \dist(x,y)^p d \pi(x,y)
\end{equation*}
where $\mathcal{C}(\mu, \lambda)$ denotes the set of couplings between $\mu$ and $\lambda$ i.e., positive Borel measures $\pi$ on $\cX \times \cX$ such that the marginal measure of $x$ is  $\lambda$ and that of $y$ is $\mu$. In order to highlight the space $\cX$ we may write $\Wdist_\cX^p$ instead of $\Wdist^p$. Moreover, if $\cX \subseteq \R^d$ and $\dist(x,y) = |y-x|$ is the Euclidean distance, we simply write  $\WN_{\cX}^p$ for $\Wdist_{\cX}^p$.

The measures $\mu$, $\lambda$ are not necessarily probabilities, however for every positive constant $a\ge 0$, it holds
\begin{equation*}
\Wdist^p( a \mu, a \lambda) = a \Wdist^p(\mu, \lambda), 
\end{equation*}
which allows in principle to often reduce to the case $\mu(\cX) = \lambda(\cX)=1$, whenever $\mu(\cX) \neq 0$. For $\alpha \ge 1$, by H\"older inequality,
 \begin{equation}\label{eq:wass-different-p}
  \bra{ \Wdist^p(\mu, \lambda)}^\alpha \le \mu(\cX)^{\alpha-1}  \Wdist^{\alpha p} (\mu, \lambda).
 \end{equation}

Given two finite or countable collections of measures $(\mu_k)_k$, $(\lambda_k)_k$ such that $\mu_k(\cX) = \lambda_k(\cX)$ for every $k$, the following subadditivity inequality holds:
\begin{equation}\label{eq:convexity-wass}
  \Wdist^p\bra{ \sum_k\mu_k, \sum_k\lambda_k} \le \sum_k   \Wdist^p\bra{ \mu_k, \lambda_k}.
\end{equation}
To prove it,  given any collection of couplings $\pi_k \in \mathcal{C}_\Omega(\mu_k, \lambda_k)$, simply notice that the measure $\sum_k \pi_k$ provides a coupling between $\sum_k\mu_k$ and $\sum_k\lambda_k$. 

Clearly, if $\dist_1(x,y) \le \dist_2(x,y)$ for every $x,y\in \cX$, where $\dist_1$, $\dist_2$ denote lower-semicontinuous pseudo-distances, it holds
\begin{equation*}
 \Wdist^p_1\bra{ \mu, \lambda} \le \Wdist^p_2(\mu, \lambda).
\end{equation*}
The trivial distance $\dist(x,y) = 1_{\cur{x \neq y}}$ yields (for any $p\ge 1$) that $\Wdist^p(\mu, \lambda)$ equals the total variation distance, $\nor{\mu-\lambda}_{\operatorname{TV}}$. Writing
\begin{equation*}
 \dist(x,y)^p \le 1_{\cur{x\neq y} } \dist(\cX)^p,
\end{equation*}
where we recall that $\dist(\cX)$ is the diameter of $\cX$ with respect to $\dist$, 
we obtain the bound
%
\begin{equation}\label{eq:wass-tv}
 \Wdist^p\bra{ \mu, \lambda} \le \nor{\mu-\lambda}_{\operatorname{TV}} \dist(\cX)^p.
\end{equation}

The general theory of optimal transportation yields  the following dual representation: 
\begin{equation}\label{eq:dual-general}
 \Wdist^p(\mu, \lambda) = \sup \cur{ \int_{\cX} u d \lambda - \int_{\cX} v  d \mu \, : \, u(x) - v(y) \le \dist(x,y)^p \quad \forall x,y \in
\cX},
\end{equation}
where the $\sup$ runs among continuous and bounded functions $u$, $v$ on $\cX$.  It is slightly convenient to restrict $v$ on a closed set $\cY\subseteq \cX$, where $\mu$ is supported. Of course, one can always let $\cY = \cX$, but a smaller $\cY$ (possibly finite) may lead to more elementary considerations. 
Duality can be then formulated in a less symmetric way, in terms of $v$ (from now on defined on $\cY$) and its $\dist^p$-transform, given by
\begin{equation*}
  \sfQ_{\dist^p} (v) (x) := \inf_{y \in \cY} \cur{ v(y) + \dist(x,y)^p}.
\end{equation*}
When both $\dist$ and $p$ are clearly understood, we simply write $\sfQ( v) = \sfQ_{\dist^p}(v)$ for such transform. It holds therefore
\begin{equation}\label{eq:dual-WDd}
 \Wdist^p(\mu, \lambda) = \sup_v \cur{ \int_{\cX} \sfQ (v) d \lambda- \int_{\cY} v  d \mu }.
\end{equation}
Indeed, if $u$ is such that $u(x)-v(y) \le \dist(x,y)^p$, then $u (x) \le \sfQ(v) (x)$, hence by replacing $u$ in the right-hand side of \Cref{eq:dual-general} with $\sfQ(v)$ one increases it. Let us notice that if $\dist(x,x') = 0$ then by the triangle inequality $\dist(x,y) = \dist(x',y)$ for every $y$, hence
\begin{equation}\label{eq:degenerate-Q-equal}
\dist(x,x') = 0 \quad \Rightarrow \quad \sfQ(v)(x)  =\sfQ(v)(x').
\end{equation}
The transform is also monotone increasing:
\begin{equation*}
  v \le \tilde v \quad \Rightarrow \quad \sfQ (v)  \le  \sfQ( \tilde v ).
\end{equation*}
If $a\in \R$, then
\begin{equation}\label{eq:Qa+v}
 \sfQ( a+v) (x) = a + \sfQ(v)(x),
\end{equation}
and setting $v(y) = 0$ for every $y\in \cY$, we find
\begin{equation*}
 \sfQ(0)(x) = \inf_{y \in \cY} \dist(x,y)^p := \dist(x,\cY)^p.
\end{equation*}
By monotonicity, we  then obtain
\begin{equation}\label{eq:sup-v-transform}
 \sup_{x \in \cX} \sfQ (v) (x)  \le \sup_{y \in \cY} v(y) + \dist(\cX,\cY)^p
\end{equation}
where we recall the notation
\begin{equation*}
 \dist(\cX,\cY) := \sup_{x \in \cX} \dist(x,\cY).
\end{equation*}
This upper bound, coupled with the lower bound
\begin{equation*}
 \inf_{x \in \cX} \sfQ (v) (x) = \inf_{y \in \cY} \cur{v(y) + \inf_{x \in \cX} \dist(x,y)^p }= \inf_{y \in \cY} v(y),
\end{equation*}
yields the bound on the oscillation
\begin{equation}\label{eq:contraction-infty}
 \osc_{\cX}\bra{ \sfQ (v)   } \le \osc_{\cY} \bra{v} +  \dist(\cX,\cY)^p.
\end{equation}
Finally, if the family $(\dist(\cdot, y)^p)_{y \in \cY}$ is uniformly continuous (with respect to the topology on $\cX$), then $\sfQ(v)$ is also continuous. In particular,
\begin{equation}\label{LipXQ}
 \Lip_{\cX} \bra{ \sfQ(v) } \le \sup_{y\in \cY} \Lip_{\cX}\bra{ \dist(\cdot, y)^p }.
 \end{equation}


The roles played by $u$ and $v$ in \Cref{eq:dual-general}  can be reversed by defining a dual $\dist^p$-transform, for $u: \cX \to \R$, given by
\begin{equation*}
 \hat{\sfQ}_{\dist^p}(u)(y) := \sup_{x \in \cX} \cur{ u(x)- \dist(x,y)^p}.
\end{equation*}
Again, if $\dist$ and $p$ are understood, we simply write $\hat{\sfQ} = \hat{\sfQ}_{\dist^p}$. The dual transform is also monotone increasing and the oscillation bound holds:
\begin{equation}\label{oschatQ}
 \osc_{\cY} \bra{ \hat{\sfQ}_{\dist^p}(u)} \le \osc_{\cX} \bra{u} +  \dist(\cX,\cY)^p.
\end{equation}
Moreover, it always holds
\begin{equation*}
  \sfQ \bra{ \hat{\sfQ} (u) } \ge u \quad \text{and}  \quad  \hat \sfQ \bra{ \sfQ (v) } \le v
\end{equation*}
hence  $\sfQ \circ  \hat{\sfQ} \circ  \sfQ  = \sfQ$. Combining these two transforms, we see in particular that in \Cref{eq:dual-WDd} we can always assume that $v = \hat \sfQ(u)$ for some $u$,  and denote in what follows $ \cS(\dist^p, \cY)$  the set of such functions. For $v \in \cS(\dist, \cY)$, its $\dist^p$-sub-differential is defined as the set
\begin{equation*}
 \partial v (y) := \cur{x \in \cX\, : \, \sfQ(v)(x) - v(y) = \dist(x,y)^p}.
\end{equation*}
Similarly, we define 
 \begin{equation*}
  \partial \sfQ(v) (x) := \cur{y \in \cY\, : \, \sfQ(v)(x) - v(y) = \dist(x,y)^p},
 \end{equation*}
 so that $x \in \partial v(y)$ if and only if $y \in \partial \sfQ(v)(x)$.

Most of the above facts actually hold for a general lower semicontinuous cost function $\mathsf{c}(x,y)$ instead of $\dist^p$. However, we additionally have that the $p$-th root of $\Wdist^p$ enjoys the triangle inequality, as a consequence of the triangle inequality of $\dist$. For our purposes, it is more useful to state and prove the following inequality for $\Wdist^p$: for every $\bar{\lambda}$ with $\bar{\lambda}(\cX) = \mu(\cX)$ and $\eps\in (0,1)$, it holds
 \begin{equation}\label{eq:triangle-wass}
  \Wdist^p(\mu, \lambda)  \le (1+ \eps)\Wdist^p(\mu, \bar{\lambda})  + \frac{c}{\eps^{p-1}}\Wdist^p(\bar{\lambda}, \lambda),
 \end{equation}
 where $c = c(p)<\infty$. Although this follows at once from the triangle inequality for $\Wdist$, we give a proof that uses the dual formulation, as it entails an intermediate inequality \Cref{eq:triangle-potential}, also useful for our purposes. Keeping $\dist$ and $p$ fixed, for any $t >0$ we define $\dist_t(x,y):= t^{\frac1 p -1} \dist(x,y)$. This scaling simply reflects on the transport cost as follows:
\begin{equation*}
 \Wdist_t^p\bra{\mu, \lambda} = t^{1-p} \Wdist^p\bra{\mu, \lambda},
\end{equation*}
but yields a modified transform
\begin{equation}\label{eq:qtv}
 \sfQ_t(v)(x) := \sfQ_{ \dist_t ^p} (v) (x):= \inf_{y \in \cY} \cur{ v(y) + \frac{ \dist(x,y)^p}{t^{p-1}}},
\end{equation}
for which  the following semigroup inequality holds:
\begin{equation*}
 \sfQ_t \bra{  \sfQ_s \bra{v}  } \ge \sfQ_{s+t}\bra{v},
\end{equation*}
where $\sfQ_s\bra{v}$ is restricted to $\cY$. Indeed, it holds
\begin{equation*}\begin{split}
 \sfQ_t \bra{  \sfQ_s \bra{v}  }(x) & = \inf_{y \in \cY} \cur{  \inf_{z \in \cY} \cur{ v(z) +  \dist_s(y,z)^p}  +  \dist_t(x, y)^p  }\\
 & = \inf_{z \in \cY} \cur{ v(z) + \inf_{y \in \cY}\cur{  \dist_t(x, y)^p   +  \dist_s(y,z)^p }  } \\ 
 & \ge \inf_{z \in \cY} \cur{ v(z)  +  \dist_{s+t}(x,z)^p  },
 \end{split}
\end{equation*}
where the second inequality follows from the triangle inequality for $\dist$ and H\"older inequality:
\begin{equation*}
  \dist(x,z) \le    t^{1-1/p}  \dist_t(x,y) +  s^{1-1/p} \dist_s(y, z)  \le \bra{ t + s }^{1-1/p} \bra{ \dist_t(x,y)^p + \dist_s(y,z)^p}^{1/p}.
\end{equation*}
Given $v: \cY \to \R$, we then write
\begin{equation*}\begin{split}
 \int_{\cX} \sfQ_{s+t}(v) d \lambda - \int_{\cY} v d \mu & = \int_{\cX} \sfQ_{s+t}(v) d \lambda - \int_{\cX}\sfQ_{s}( v) d \bar{\lambda}+\int_{\cX}\sfQ_{s}( v) d \bar{\lambda} -  \int_{\cY} v d \mu\\
 & \le \sqa{ \int_{\cX} \sfQ_{t}\bra{ \sfQ_{s}(v)} d \lambda - \int_{\cX}\sfQ_{s}( v) d \bar{\lambda}} +\sqa{\int_{\cX}\sfQ_{s}( v) d \bar{\lambda} -  \int_{\cY} v d \mu}\\
 & \le \Wdist_t^p\bra{ \bar{\lambda}, \lambda} +\sqa{\int_{\cX}\sfQ_{s}( v) d \bar{\lambda} -  \int_{\cY} v d \mu}.
 \end{split}
\end{equation*}
The choices $s=1$ and $t=\eps$ lead to the inequality
\begin{equation}\label{eq:triangle-potential}
 \sqa{\int_{\cX}\sfQ_{1}( v) d \lambda -  \int_{\cY} v d \mu} \ge \sqa{ \int_{\cX} \sfQ_{1+\eps}(v) d \lambda - \int_{\cY} v d \mu } - \eps^{1-p}\Wdist^p\bra{ \bar{\lambda}, \lambda},
\end{equation}
that we collect for later use. Back to the general case, further bounding from above we find
\begin{equation*}\begin{split}
\int_{\cX} \sfQ_{s+t}(v) d \lambda - \int_{\cY} v d \mu &\le \Wdist_t^p\bra{ \bar{\lambda}, \lambda} +\sqa{\int_{\cX}\sfQ_{s}( v) d \bar{\lambda} -  \int_{\cY} v d \mu}\\
& \le  t^{1-p} \Wdist^p\bra{ \bar{\lambda}, \lambda} + s^{1-p} \Wdist^p\bra{\mu, \bar{\lambda}}.
\end{split}
\end{equation*}
Finally, taking the supremum over $v$ gives
\begin{equation*}
 (s+t)^{1-p} \Wdist^p\bra{\mu, \lambda} \le  t^{1-p} \Wdist^p\bra{ \bar{\lambda}, \lambda} + s^{1-p} \Wdist^p\bra{\mu, \bar{\lambda}}.
\end{equation*}
Choosing $s=\eps / c$ and $t = 1-s$ yields \Cref{eq:triangle-wass}, provided that
\begin{equation*}
 (1-\eps/c)^{1-p} \le 1+  \eps
\end{equation*}
which holds if $c = c(p) < \infty$ is sufficiently large.

As a consequence of \Cref{eq:convexity-wass} and \Cref{eq:triangle-wass} we obtain that,  $\Wdist^p$ is (approximately) sub-additive: for any finite or countable Borel partition $\cX  = \bigcup_{k} \Omega_k$, it holds
\begin{equation}\label{eq:Neumann-exact-sub-additive}
 \Wdist^p\bra{ \mu, \lambda} \le (1+\eps) \sum_k \Wdist^p\bra{ \chi_{\Omega_k} \mu, \frac{\mu (\Omega_k)}{\lambda(\Omega_k)} \chi_{\Omega_k} \lambda} + \frac{c}{\eps^{p-1} } \Wdist^p\bra{ \sum_k \frac{\mu (\Omega_k)}{\lambda(\Omega_k)}  \chi_{\Omega_k} \lambda, \lambda},
\end{equation}
where $c = c(p)<\infty$. Indeed, it is sufficient to apply \Cref{eq:triangle-wass} with
\begin{equation}\label{eq:nu-def}
 \bar{\lambda} :=  \sum_k  \frac{\lambda (\Omega_k)}{\mu(\Omega_k)}  \chi_{\Omega_k} \lambda,
\end{equation}
and subsequently \Cref{eq:convexity-wass} with
\begin{equation*}\mu_k = \chi_{\Omega_k} \mu \quad \text{and} \quad \lambda_k:=\frac{\mu (\Omega_k)}{\lambda(\Omega_k)} \chi_{\Omega_k} \lambda.
\end{equation*}

%

\subsection{Boundary cost}\label{sec:boundary} Let us assume that $(\cX, \dist)$ is an (actual) compact metric space. Given $p\ge 1$, we define for any  $\Omega \subseteq \cX$ the following ``boundary'' pseudo-distance, for $x$, $y \in \cX$:
\begin{equation*}
 \bist_\Omega(x,y) := \min\cur{ \dist(x,y), \bra{ \dist(x,\Omega^c)^{p} + \dist(y, \Omega^c)^{p}}^{1/p} }.
\end{equation*}
We remark that the notation $\bist_{\Omega}$ is slightly ambiguous, as it does not make explicit the dependence on $\dist$ and also $p$, which will be anyway made clear by considering more often the quantity $\bist_{\Omega}^p$. It is not difficult to check that $\bist_{\Omega}$ defines a pseudo-distance such that  $\bist_{\Omega} = \bist_{\Omega^\circ}$, where $\Omega^\circ$ denotes the interior part of $\Omega$. Moreover, $\bist_{\Omega} \le \dist$ and
\begin{equation*}
 \bist_{\Omega}(x,y) = \dist(y, \Omega^c) \quad \text{if $x \in  \overline{\Omega^c}$,}
\end{equation*}
hence $\bist_{\Omega}(x,x') = 0$ if $x$, $x' \in \overline{\Omega^c}$. Given measures $\mu$, $\lambda$ on $\cX$ with $\mu(\cX) = \lambda(\cX)$, we write accordingly to general the notation introduced above, with $\bist_{\Omega}$ in place of $\dist$,
\begin{equation*}
 \WD_{\Omega}^p\bra{\mu, \lambda} =  \inf_{\pi\in\mathcal{C}(\mu,\lambda)}  \int_{\cX \times \cX} \bist_{\Omega}(x,y)^p d \pi(x,y)
\end{equation*}
for the primal problem and
\begin{equation}\label{eq:dual-WD}
 \WD^p_{\Omega}(\mu, \lambda) = \sup_{v\in \cS(\bist_{\Omega}^p, \cY)} \cur{ \int_{\cX} \sfQ_{\bist_{\Omega}^p} (v) d \lambda- \int_{\cY} v  d \mu }
\end{equation}
for the dual formulation.
Let us notice that if $\mu(\Omega) = \lambda(\Omega)$,  it holds
\begin{equation*}
 \WD_{\Omega}^p\bra{ \mu, \lambda}  = \WD_{\Omega}^p\bra{\chi_{\Omega} \mu, \chi_{\Omega} \lambda},
\end{equation*}
for one can transport without cost the respective masses in $\Omega^c$. In all what follows, we always consider $\WD_{\Omega}^p(\mu, \lambda)$ for measures such that $\mu(\Omega) = \lambda(\Omega)$. 
On the dual side, from \Cref{eq:degenerate-Q-equal} it follows that 
$\sfQ_{\bist_{\Omega}^p} (v)$ is constant on $\overline{\Omega^c}$, and in particular
\begin{equation*}
 \sfQ_{\bist_{\Omega}^p} (v)(x) = \inf_{y \in \cY} \cur{ v(y) + \dist(y, \Omega^c)^p} \quad \text{for every $x \in \overline{\Omega^c}$.}
\end{equation*}
We denote in what follows with $\cS_0(\bist_{\Omega}^p, \cY)$ the subset of $\cS(\bist_{\Omega}^p, \cY)$ such that $ \sfQ_{\bist_{\Omega}^p} (v) =0$ on $\overline{\Omega^c}$, i.e.,
\begin{equation*}
\inf_{y \in \cY} \cur{ v(y) + \dist(y, \Omega^c)^p} = 0.
\end{equation*}
For $v \in \cS_0(\bist_{\Omega}^p, \cY)$ it also holds $v(y) = 0$ for every $y \in \cY  \cap \overline{\Omega^c}$. Indeed, we have immediately that
\begin{equation}\label{eq:v-lower-bound}
 v(y) \ge - \dist(y, \Omega^c)^p \quad \text{for every $y \in \cY$,}
\end{equation}
and in particular $v(y) \ge 0$ for $y \in \cY \cap \overline{\Omega^c}$. Moreover, for a general $x \in \cX$, we find
\begin{equation}\label{eq:upper-bound-qv-x}
\begin{split}
 \sfQ_{\bist_{\Omega}^p}(v)(x) & = \inf_{y \in \cY} \cur{ v(y) + \bist_{\Omega}(x,y)^p} \\
 & \le  \inf_{y \in \cY} \cur{ v(y) + \dist(y, \Omega^c)^p}  + \dist(x,\Omega^c)^p = \dist(x,\Omega^c)^p.
 \end{split}
\end{equation}
This inequality entails, for $y \in \cY \cap  \overline{\Omega^c}$, that
\begin{equation*}\begin{split}
 v(y)  & = \sup_{x \in \cX}\cur{ \sfQ_{\bist_{\Omega}^p}(v)(x)  - \bist_{\Omega}(x,y)^p } \\
 & \le \sup_{x \in \cX} \cur{\dist(x,\Omega^c)^p - \dist(x,\Omega^c)^p } = 0,
 \end{split}
\end{equation*}
hence  $v(y) = 0$ for $y \in \cY \cap \overline{\Omega^c}$.

Let us notice that, by \Cref{eq:Qa+v} and the analogue property for $\hat{\sfQ}$, one can always restrict the supremum in \Cref{eq:dual-WD} to $v \in \cS_0(\bist_{\Omega}^p, \cY)$.
 Furthermore, since $x\mapsto \bist_{\Omega}(x,y)$ is $1$-Lipschitz with respect to $\dist$ (for every $y$), and
\begin{equation*}
 \bist_{\Omega}(\cX) = \sup_{x, x' \in \cX} \bist_{\Omega}(x,x') \le \dist(\Omega)
\end{equation*}
we obtain by \Cref{LipXQ} and \Cref{eq:contraction-infty} that $\sfQ_{\bist_{\Omega}^p}(v)$ is Lipschitz (with respect to $\dist$) with
\begin{equation}\label{eq:osc-Qv}
 \Lip_{\cX}\bra{\sfQ_{\bist_{\Omega}^p}(v)  } \le p \dist(\Omega)^{p-1}, \quad \text{and} \quad  \osc_{\cX}\bra{\sfQ_{\bist_{\Omega}^p}(v) } \le \osc_{\cY} \bra{v} + \dist(\Omega, \cY)^p.
\end{equation}
Similarly,
\begin{equation*}
 \Lip_{\cY}\bra{\hat \sfQ_{\bist_{\Omega}^p}(u)  } \le p \dist(\Omega)^{p-1},\quad \text{and} \quad  \osc_{\cY}\bra{\hat \sfQ_{\bist_{\Omega}^p}(u) } \le \osc_{\cX} \bra{u} + \dist(\Omega, \cY)^p.
\end{equation*}
In particular, any $v \in \cS(\bist_{\Omega}^p, \cY)$ is Lipschitz with
\begin{equation}\label{eq:osc-y-v}
\Lip_{\cY}(v) \le p \dist(\Omega)^{p-1},\quad \text{hence} \quad  \osc_{\cY}\bra{v} \le  p \dist(\Omega)^{p}, \quad \osc_{\cX}\bra{ \sfQ_{\bist_{\Omega}^p} (v)} \le p \dist(\Omega)^p + \dist(\Omega, \cY)^p.
\end{equation}
If moreover $v \in \cS_0(\bist_{\Omega}^p, \cY)$, then
\begin{equation}
 \label{eq:infty-q-v-x}
 \nor{v}_\infty \le   p \dist(\Omega)^p, \quad \text{and} \quad  \nor{ \sfQ_{\bist_{\Omega}^p} (v)}_\infty \le   p \dist(\Omega)^p + \dist(\Omega, \cY)^p.
\end{equation}

The quantity $\WD_{\Omega}^p(\mu, \lambda)$ enjoys the following approximate super-additivity property: for any finite Borel partition $\Omega  = \bigcup_{k} \Omega_k$, it holds
\begin{equation}\label{eq:Dirichlet-sub-additive}
  \WD^p_{\Omega}\bra{ \mu, \lambda} \ge  (1+\eps)^{1-p} \sum_k \WD^p_{\Omega_k} \bra{  \mu, \frac{\mu (\Omega_k)}{\lambda(\Omega_k)} \lambda} - \eps^{1-p} \Wdist^p_{\Omega}\bra{ \sum_k \frac{\mu (\Omega_k)}{\lambda(\Omega_k)}  \chi_{\Omega_k} \lambda, \lambda},
\end{equation}
for any $\eps \in (0,1)$.  To prove it, we rely on the following construction.


\begin{lemma}[gluing]\label{lem:potential-gluing}
Consider a finite Borel partition $\Omega = \bigcup_{k} \Omega_k$ such that, for every $x \in \Omega_k$, $y \in \Omega_h$ with $h\neq k$, it holds
\begin{equation}\label{eq:length}
 \dist(x,y) \ge \dist(x, \Omega_k^c) + \dist(y, \Omega_h^c).
\end{equation}
For every $k$, let $v_k \in \cS_0(\bist_{\Omega_k}^p, \cY)$ and define $v: \cY \to \R$, $v(y) := v_k(y)$ for $y \in \cY \cap \overline{\Omega_k}$ and $v(y) = 0$ for $y \in \cY \cap \overline{\Omega^c}$. 
Then, $v$ is well-defined, with
\begin{equation}\label{eq:norm-v-glued}
 \nor{v}_\infty \le    p  \max_k \dist(\Omega_k)^{p}, \quad \nor{\sfQ_{\bist_{\Omega}^p} (v ) }_\infty \le  p\max_k \dist(\Omega_k)^{p} + \dist(\Omega, \cY)^p
\end{equation} and it holds, for every $k$,
\begin{equation}\label{eq:ge-glued}
 \sfQ_{\bist_{\Omega}^p} (v )  \ge \sfQ_{\bist_{\Omega_k}^p}(v_k) \quad \text{on $\Omega_k$.}
\end{equation}
\end{lemma}


\begin{remark}
 Condition \Cref{eq:length} is always satisfied if $(\cX, \dist)$ is a length space.
\end{remark}

\begin{proof}  We notice that if $y \in \cY$ belongs both to $\overline{\Omega_k}$ and $\overline{\Omega_h}$ with $h \neq k$, it follows that $y \in \overline{\Omega_k^c}$ and $y \in \overline{\Omega_h^c}$, hence $v_k(y) = v_h(y) = 0$ and thus $v$ is well-defined.
Inequalities  \Cref{eq:norm-v-glued} are  immediate from \Cref{eq:infty-q-v-x}. 
To prove \Cref{eq:ge-glued}, fix $k$ and  $x \in \overline{\Omega_k}$. If $y \in \cY \cap \overline{\Omega_k}$, then,
\begin{equation*}
 \dist(x, \Omega^c) \ge \dist(x, \Omega^c_k) \quad \text{and} \quad \dist(y, \Omega^c) \ge \dist(y, \Omega^c_k),
\end{equation*}
hence
\begin{equation*}
 \bist_{\Omega}(x,y)^p \ge \bist_{\Omega_k}(x,y)^p
\end{equation*}
and therefore
\begin{equation*}
   \inf_{y \in \cY \cap \overline{\Omega_k}} \cur{ v(y)+  \bist_{\Omega}(x,y)^p  } \ge  \inf_{y \in \cY \cap \overline{\Omega_k}} \cur{ v_k(y)+\bist_{\Omega_k}(x,y)^p } \ge \sfQ_{\bist_{\Omega_k}^p} (v_k)(x). 
\end{equation*}
Using that from \Cref{eq:upper-bound-qv-x} with $\Omega_k$ instead of $\Omega$ and $v_k$ instead of $v$,
\begin{equation*}
  \dist(x, \Omega_k^c)^p \ge \sfQ_{\bist_{\Omega_k}^p} (v_k)(x),
\end{equation*}
 in order to prove \Cref{eq:ge-glued}, it is thus enough to argue that
\begin{equation*}
  \inf_{y \in \cY \cap \overline{\Omega_k^c}} \cur{ v(y)+  \bist_{\Omega}(x,y)^p  } \ge \dist(x, \Omega_k^c)^p.
\end{equation*}
Given $y \in \cY \cap \overline{\Omega_k^c}$,  we distinguish two cases. Assume first that $y \in \overline{\Omega_h}$ for some $h \neq k$, then by \Cref{eq:v-lower-bound} with $\Omega_h$ instead of $\Omega$ and $v_h$ instead of $v$,
\begin{equation*}
  v(y) = v_h(y) \ge -\dist(y, \Omega_h^c)^p.
\end{equation*}
%
 Moreover,
\begin{equation*}\begin{split}
 \bist_{\Omega}(x,y)^p  &=  \min \cur{ \dist(x,y)^p, \dist(x, \Omega^c)^p + \dist(y, \Omega^c)^p} \\
 & \stackrel{\Cref{eq:length}}{\ge}   \min \cur{ \bra{ \dist(x, \Omega_k^c) + \dist(y, \Omega_h^c)}^p, \dist(x, \Omega^c_k)^p + \dist(y, \Omega^c_h)^p}\\
 & \ge  \dist(x, \Omega^c_k)^p +  \dist(y, \Omega_h^c)^p.
 \end{split}
\end{equation*}
Therefore, we find in this first case that
\begin{equation}\label{eq:claim-y-gluing}
 v(y) + \bist_{\Omega}(x,y)^p \ge  \dist(x, \Omega^c_k)^p.
\end{equation}
If instead $y \in \overline{\Omega^c}$, so that $v(y) = 0$, and $\bist_{\Omega}(x,y)  = \dist(x, \Omega^c) \ge \dist(x, \Omega_k^c)$, we conclude that \Cref{eq:claim-y-gluing} holds as well.
%
%
%
\end{proof}

In order to  establish \Cref{eq:Dirichlet-sub-additive}, we need to apply the result above to the modified cost $\bist_{\Omega, t}(x,y) = t^{\frac{1}{p}-1}\bist_{\Omega}(x,y)$, for $t>0$. We  notice that the properties and definitions introduced above easily extend from the case $t=1$ to the general case, with minor modifications in the bounds -- to see this, it is sufficient to replace the original $\dist$ with $\dist_t$.

Consider then $v_k \in \cS_0(\dist_{\Omega_k, (1+\eps) }^p, \cY)$ that is optimal in the dual formulation of  the cost, so that
\begin{equation*}
\int_{\Omega_k} \sfQ_{\dist_{\Omega,(1+\eps)}^p}(v_k) d \lambda - \int_{\cY\cap \Omega_k} v_k d \mu   = (1+\eps)^{p-1} \WD_{\Omega_k}^p\bra{\mu,   \frac{\mu(\Omega_k)}{\lambda(\Omega_k)} \lambda}.
\end{equation*} 
By \Cref{eq:triangle-potential} with $\bar{\lambda}$ as in the approximate sub-additivity \Cref{eq:nu-def}, $\sfQ = \sfQ_{\bist_{\Omega}^p}$ and $v$ as provided by \Cref{lem:potential-gluing}, we find
\begin{equation*}\begin{split}
  \sqa{\int_{\cX}\sfQ_{\bist_{\Omega}^p}( v) d \lambda -  \int_{\cY} v d \mu} & \ge \sqa{ \int_{\cX} \sfQ_{\bist_{\Omega}^p, (1+\eps) }(v) d \lambda - \int_{\cY} v d \mu } - \eps^{1-p}\Wdist^p\bra{ \bar{\lambda}, \lambda}\\
  &   \ge \sum_k \sqa{ \int_{\Omega_k} \sfQ_{\bist_{\Omega}^p, (1+\eps) }(v) d \lambda - \int_{\cY\cap \Omega_k} v d \mu } - \eps^{1-p}\Wdist^p\bra{ \bar{\lambda}, \lambda}\\
  & \stackrel{\Cref{eq:ge-glued}}{\ge} \sum_k \sqa{ \int_{\Omega_k} \sfQ_{\dist_{\Omega_k}^p, (1+\eps) }(v_k) d \lambda - \int_{\cY\cap \Omega_k} v_k d \mu } - \eps^{1-p}\Wdist^p\bra{ \bar{\lambda}, \lambda}\\
  & = (1+\eps)^{1-p} \sum_k  \WD_{\Omega_k}^p\bra{\mu,   \frac{\mu(\Omega_k)}{\lambda(\Omega_k)} \lambda} - \eps^{1-p}\Wdist^p\bra{ \bar{\lambda}, \lambda}.
  \end{split}
\end{equation*}
If we further bound from above the left-hand side
\begin{equation*}
 \int_{\cX}\sfQ_{\bist_{\Omega}^p}( v) d \lambda -  \int_{\cY} v d \mu \le  \WD_\Omega^p(\mu, \lambda),
\end{equation*}
inequality \Cref{eq:Dirichlet-sub-additive} then follows. For later use, let us notice that we may also bound from above the left-hand side replacing $v$ with
\begin{equation}\label{eq:tilde-v-superadd}
 \tilde v := \hat \sfQ_{\bist_{\Omega}^p}\bra{ \sfQ_{\bist_{\Omega}^p}( v) } \in \cS\bra{ \bist_{\Omega}^p, \cY}
\end{equation}
for which $\sfQ_{\bist_{\Omega}^p}( v) = \sfQ_{\bist_{\Omega}^p}( \tilde v)$ (recall that $\sfQ\circ \hat \sfQ \circ \sfQ=\sfQ$). We find
\begin{equation}\label{eq:potential-subadditive-new}
  \int_{\cX}\sfQ_{\bist_{\Omega}^p}( \tilde v) d \lambda -  \int_{\cY} \tilde v d \mu \ge  (1+\eps)^{1-p} \sum_k  \WD_{\Omega_k}^p\bra{\mu,   \frac{\mu(\Omega_k)}{\lambda(\Omega_k)} \lambda} - \eps^{1-p}\Wdist^p\bra{ \sum_k\frac{\mu(\Omega_k)}{\lambda(\Omega_k)} \chi_{\Omega_k}\lambda , \lambda}.
\end{equation}
Let us also notice that by \Cref{eq:norm-v-glued}
\begin{equation}\label{eq:norm-tilde-v}
 \nor{ \sfQ_{\bist_{\Omega}^p} (\tilde v) }_\infty \le p\max_k \dist(\Omega_k)^{p} + \dist(\Omega, \cY)^p .
\end{equation}

We end this section with a general result providing uniform upper bounds for the optimal transport plan, extending in particular \parencite[Lemma 4.4]{AGT19} (see also \parencite{bouchitte2007new,GO}). As always in this section, we consider $\mu$, $\lambda$ to be measures on $\cX$ such that $\mu(\Omega) = \lambda(\Omega)$ and $\mu$ supported on $\cY$.

\begin{lemma}\label{lem:uniform-upper-bounds}
  Let  $v \in \mathcal{S}(\bist_{\Omega}^p, \cY)$  and set
 \begin{equation}\label{eq:def-ell}
  \ell := \bra{\osc_{\cY}\bra{v} + \dist(\cX, \cY)^p }^{1/p}.
 \end{equation}
 \begin{enumerate}
 \item For every $x \in \partial_{\bist_{\Omega}^p} v(y)$, it holds
 \begin{equation}\label{eq:l-infty-bound-dirichlet}
   \bist_{\Omega}(x,y) \le \ell 
 \end{equation}
 and therefore (recall the definition \Cref{defOmr} of $\Omega^{(\ell)}$),
 \begin{equation}\label{eq:dist-omega-dist-p-l-infty-bound}
 x \in \Omega^{(\ell)} \quad \text{or} \quad y \in \Omega^{(\ell)}
 \quad \Rightarrow \quad \bist_{\Omega}(x,y) = \dist(x,y).
  \end{equation}

\item  Assume that $\Omega$ is midpoint convex, i.e., for every $x,y \in \Omega$ there exists   $z \in \Omega$ such that $\dist(x,z)=\dist(z,y) = \dist(x,y)/2$, and that for some $d>0$, $\lambda_0 > 0$, it holds
 \begin{equation}\label{eq:ahlfors}
   \lambda(B_r(x)) \ge \lambda_0 r^d \quad \text{for  every $r \le \ell$, and $\lambda$-a.e.\  $x \in \Omega$.}
 \end{equation}
 Let $\pi \in \mathcal{C}(\mu, \lambda)$ and $v \in \mathcal{S}(\bist_{\Omega}^p, \cY)$ be optimizers for $\WD_{\Omega}^p(\mu, \lambda)$ respectively in the primal and dual formulation, then, for $\pi$-a.e.\ $(x,y)$,
 \begin{equation}\label{eq:l-infty-lp-bound-dirichlet}
  x \in \Omega^{(2\ell)} 
  \quad \Rightarrow \quad  \lambda_0 \dist(x,y)^{p+d} \les_{p}  \WD_{\Omega}^p\bra{\mu, \lambda}.
 \end{equation}
 \end{enumerate}
 \end{lemma}


 \begin{proof}
 We set for brevity $u:= \sfQ_{\bist^p_{\Omega}}(v)$.
To prove \emph{(1)}, we first  recall that by definition of the $\bist^p_{\Omega}$-sub-differential, we have for $x\in \partial_{\bist_{\Omega}^p} v(y)$,
 \begin{equation}\label{eq:pi-ae}
 u (x) - v(y) = \bist_{\Omega}(x,y)^p.
\end{equation}
Thus
\begin{equation*}\begin{split}
 \bist_{\Omega}(x,y)^p & \le \sup_{x'\in \cX} u(x') - \inf_{y' \in \cY} v(y') \le \sup_{y' \in \cY}v(y') - \inf_{y' \in \cY} v(y') + \bist_{\Omega}(\cX, \cY)^p\\
 & \le \osc_{\cY}(v) + \dist(\cX, \cY)^p = \ell^p,
 \end{split}
\end{equation*}
having used \Cref{eq:sup-v-transform}. 
 Hence, inequality \Cref{eq:l-infty-bound-dirichlet} is proved and \Cref{eq:dist-omega-dist-p-l-infty-bound} follows from the  definitions of $\Omega^{(\ell)}$ and $\bist_{\Omega}$. 

We now turn to  \emph{(2)}. We first recall that  since both $\pi$ and $v$ are optimizers for their respective formulations, it follows that for $\pi$-a.e.\ $(x,y)$ we have $x\in \partial_{\bist_{\Omega}^p} v(y)$ or equivalently \eqref{eq:pi-ae}.
 The identity extends by continuity to every $(x,y) \in \supp \pi$. This is because inequality $\le$ holds for every $(x,y)$, and integrating with respect to $\pi$, we have
\begin{equation*}
 \int_{\cX \times \cY} \bra{ u (x) - v(y) - \bist_{\Omega}(x,y)^p} d \pi(x,y) = \WD_{\Omega}^p(\mu, \lambda) -  \WD_{\Omega}^p(\mu, \lambda) =0.
\end{equation*}
Let us notice that \Cref{eq:pi-ae} and the validity of the inequality $\le$ for general pairs (not in $\supp \pi$) yields the inequality
\begin{equation}\label{eq:monotonicity}
 \bist_{\Omega}(x,y)^p + \bist_{\Omega}(x',y')^p \le \bist_{\Omega}(x,y')^p + \bist_{\Omega}(x,y')^p, \quad \text{for every $(x,y), (x',y') \in \supp \pi$.}
\end{equation}
Consider $(x,y) \in \supp \pi$ with $x \in \Omega^{(2\ell)}$, so that in particular $\dist(x,y) \le \ell$. 
  Set for simplicity of notation $\alpha := \dist(x,y)/2$ and let $z$ be a midpoint  between $x$ and $y$, so that
 \begin{equation*}
 \dist(x,z) = \dist(y,z) = \alpha \le \ell/2,
 \end{equation*}
 and therefore $\dist(z, \Omega^c) \ge 2 \ell - \ell/2 \ge \ell$. Consider the ball $B_{\eps \alpha}(z)$, where $\eps = \eps(p) \in (0,1)$ is to be specified below, sufficiently small. Notice that,  for every $x' \in B_{\eps \alpha}(z)$, it holds
\begin{equation*}
\dist(x',x)  \le  \dist(x',z)+\dist(x,z)  \le (1+\eps)\alpha \le 2 \alpha \le \ell,
\end{equation*}
hence, by the triangle inequality, $x' \in \Omega^{(\ell)}$.
By \Cref{eq:monotonicity} applied with $(x',y') \in \supp \pi$ and $x' \in B_{\eps \alpha}(z)$, it follows that
\begin{equation}\label{eq:p-monotonicity}\begin{split}
\dist(x,y)^p + \dist(x',y')^p & =  \bist_{\Omega}(x, y)^p + \bist_{\Omega}(x', y')^p\\
& \le  \bist_{\Omega}(x', y) + \bist_{\Omega}(x, y')^p \\
& \le \dist(x',y)^p + \dist(x,y')^p.
\end{split}
\end{equation}
%
Then, we use the triangle inequality  to bound from above
\begin{equation*}
\dist(x',y)^p \le \bra{ \dist(x',z)+\dist(z,y) }^p \le (1+\eps)^p \alpha^p 
\end{equation*}
and similarly, with an additional use of the elementary bound $(a+b)^p \le (1+\eps)a^p + c\eps^{-(p-1)} b^p$ where $c = c(p)<\infty$,
\begin{equation*}\begin{split}
 \dist(x,y')^p &\le \bra{ \dist(x,z)+ \dist(z,x') + \dist(x',y')}^p \le \bra{ (1+\eps) \alpha^p + \dist(x',y') }^p  \\
 & \le  (1+\eps)^{p+1} \alpha^p + \frac{c}{\eps^{p-1}}\dist(x',y')^p.
 \end{split}
\end{equation*}
Using these bounds in \Cref{eq:p-monotonicity} and recalling that $\dist(x,y)^p = 2^p \alpha^p$, we find the inequality
\begin{equation*}
 2^p \alpha^p + \dist(x',y')^p \le (2+\eps) (1+\eps)^p \alpha^p + \frac{c(p)}{\eps^{p-1}}\dist(x',y')^p.
\end{equation*}
Choosing $\eps = \eps(p)\in (0,1)$ sufficiently small, e.g.\ such that
\begin{equation*}
 2^p  -(2+\eps) (1+\eps)^p \ge 2^{p-1}-1,
\end{equation*}
we deduce that
\begin{equation*}
\alpha^p  \les_p \dist(x',y')^p.
\end{equation*}
Recalling that $\dist(x',y') = \bist_{\Omega}(x', y')$, and integrating with respect to $(x',y') \in \supp \pi$ such that $x' \in B_{\eps \alpha}(z)$ -- and using the assumption \Cref{eq:ahlfors} -- we deduce that
\begin{equation*}
 \lambda_0 \alpha^{p+d} \les_p \int_{\cX \times \cX} \bist_{\Omega}(x', y')^p d \pi(x', y')= \WD_{\Omega}^p(\mu, \lambda),
\end{equation*}
hence \Cref{eq:l-infty-lp-bound-dirichlet} follows.
 \end{proof}

\subsection{The Euclidean case}
We specialize to the Euclidean setting, i.e., assume that $\cX \subseteq \R^d$ is compact and convex (e.g.\ is a sufficiently large cube) and $\dist(x,y) = |y-x|$ is the Euclidean distance. We momentarily keep $p$ general, although our main result holds for $p=2$ only. Given $\Omega \subseteq \cX$ Borel, since the function $\bist_{\Omega}$ is Lipschitz continuous with Lipschitz constant bounded by $p\diam(\Omega)^{p-1}$, it follows that any transform $u = \sfQ_{\bist_{\Omega}^p}(v)$ is Lipschitz and by Rademacher's theorem, Lebesgue-a.e.\ differentiable. The classical Brenier theorem \parencite{Br91} exploits this fact to characterize the optimal transport plan from an absolutely continuous measure  in terms of a map of gradient type: in the following result we extend it to the ``boundary'' case.

\begin{lemma}
 Let $p>1$, $p'=p/(p-1)$, $\lambda$, $\mu$ be measures with $\lambda(\Omega) = \mu(\Omega)$, $\lambda$ absolutely continuous with respect to the Lebesgue measure and $\mu$ supported on $\cY$. Let $\pi\in \mathcal{C}(\mu, \lambda)$ and $v \in \cS(\bist_{\Omega}^p, \cY)$ be respectively optimizers for $\WD_{\Omega}^p(\mu, \lambda)$ in the primal and dual formulation and set
 \begin{equation*}
  u:= \sfQ_{\bist_{\Omega}^p}(v).
 \end{equation*}
 Then,
 \begin{equation}\label{eq:estimnabvp}
  \int_{\Omega} |\nabla u|^{p'} \lambda \le p^{p'} \WD_\Omega^p(\mu,\lambda),
 \end{equation}
and with $\ell$ defined as in \eqref{eq:def-ell},  for $\pi$-a.e.\ $(x,y)$, it holds
  \begin{equation}\label{eq:brenier}
   x \in \Omega^{(\ell)} \Rightarrow y = x -  \bra{ \frac{ \nabla u (x)}{p}}^{\bra{p'/p}}. 
  \end{equation}
\end{lemma}

In \Cref{eq:brenier} we used conveniently the notation $x^{(\alpha)} = |x|^{\alpha-1} x$.


   \begin{remark}
  As a consequence of \Cref{eq:estimnabvp} and \Cref{lem:gagliardo} with $p'=p/(p-1)$ in place of $p$, we also collect the following bound:
  \begin{equation}\label{eq:osc-u-gagliardo}
   \osc_{\cX}(u)^{p'+d} \les   \diam(\Omega)^{d(p-1)+p'} \WD_{\Omega}^p\bra{ \mu, \lambda},
  \end{equation}
having used that $\Lip_{\cX}(u) \le p \diam(\Omega)^{p-1}$. The implicit constant depends on $p$, $d$ and  $\Omega$ but is invariant with respect to dilations.
 \end{remark}

 \begin{proof}
 Both $u$ and  $x \mapsto \dist(x,\Omega^c)$ are Lipschitz continuous hence differentiable Lebesgue-a.e.\ in $\cX$. Moreover,
 \begin{equation}\label{eq:nor-nabla-dist}
   \nor{ \nabla \dist(\cdot, \Omega^c)}_\infty  \le \Lip_{\cX}\bra{ \dist(\cdot, \Omega^c)} \le 1.
 \end{equation}
Let  $(x,y)$ be such that \Cref{eq:pi-ae} holds and $x$ is a point of differentiability for both functions. Then, we prove
 \begin{equation}\label{claimnablav}
  |\nabla u(x)|^{\frac{p}{p-1}}\le p^{\frac{p}{p-1}} \bist_\Omega^p\bra{x,y},
 \end{equation}
so that integration with respect to $\pi$ yields \Cref{eq:estimnabvp}. We consider separately the case $\bist_\Omega\bra{x,y}=|x-y|$ and the case $\bist_\Omega\bra{x,y}= \bra{ \dist(x,\Omega^c)^p+\dist(y,\Omega^c)^p}^{1/p}$. In the first case, we have by definition of $u$, for every $x'\in \cX$,
\begin{equation*}
 u(x')\le v(y)+\bist_\Omega\bra{x',y}^p\le v(y)+|x'-y|^p,
\end{equation*}
with equality at $x'=x$. Differentiating at $x'=x$ yields 
\begin{equation}\label{eq:nabla-u-x}
 \nabla u(x)=p |x-y|^{p-2}(x-y) = p (x-y)^{(p-1)}
\end{equation}
so that
\begin{equation*}
 |\nabla u(x)|^{\frac{p}{p-1}}= p^{\frac{p}{p-1}}|x-y|^p=p^{\frac{p}{p-1}}\bist_\Omega\bra{x,y}^p.
\end{equation*}
If instead $\bist_\Omega\bra{x,y}=\bra{\dist(x,\Omega^c)^p+\dist(y,\Omega^c)^p}^{1/p}$, we have, for every $x'$,
\begin{equation*}
 u(x')\le v(y)+\bist_\Omega\bra{x',y}^p\le v(y)+\dist(x',\Omega^c)^p+\dist(y,\Omega^c)^p,
\end{equation*}
with equality for $x'=x$. Differentiating again at $x'=x$ yields
\begin{equation*}
 \nabla u (x)=p \dist(x, \Omega^c)^{p-1} \nabla \dist(x, \Omega^c).
\end{equation*}
Using \Cref{eq:nor-nabla-dist}, we find
\begin{equation*}
 |\nabla u (x)|^{\frac{p}{p-1}}\le p^{\frac{p}{p-1}}\dist(x, \Omega^c)^{p} \le p^{\frac{p}{p-1}} \bist_\Omega\bra{x,y}^p,
\end{equation*}
which concludes the proof of \Cref{claimnablav}. Finally, to prove \Cref{eq:brenier}, given $x\in \Omega^{(\ell)}$,  by \Cref{eq:dist-omega-dist-p-l-infty-bound} we have $|x-y| = \bist_{\Omega}(x,y)$ and therefore identity \Cref{eq:nabla-u-x} holds. Composing both sides with the function
\begin{equation*}
 \R^d \ni z\mapsto  z^{\bra{p'/p}} = z ^{\bra{ \frac 1 {p-1} }} = |z|^{\frac{2-p}{p-1}} z
\end{equation*}
yields \Cref{eq:brenier}.
 \end{proof}

Next, we recall a general upper bound for the optimal transport cost in terms of a negative Sobolev norm, see  
e.g.\ \parencite{peyre2018comparison}, but also \parencite[Proposition 2.3]{ambrosio2019pde}, \parencite[Theorem 2]{Le17}, and \parencite{goldman2021convergence} in the context of random optimal transport problems.  

\begin{lemma}\label{lem:peyre} Let $\Omega \subseteq \R^d$ be a bounded domain with  Lipschitz boundary, $\eta: \Omega \to (0, \infty)$ be bounded and such that a weighted Poincaré inequality \Cref{eq:weighted-poincare} holds. Given finite measures $\mu= f\eta$, $\lambda= g\eta$, both absolutely continuous with respect to $\eta$, and such that $f$, $g \in L^p(\eta)$ and for some $g_0 >0$ one has $g(x) \ge g_0>0$ for every $x\in \Omega$, it holds
 \begin{equation*}
  \WN_{\Omega}^p(\mu,\lambda) 
   \les  \frac{1}{g_0^{p-1}} \nor{f-g}_{H^{-1,p}(\eta)}^p.
 \end{equation*}
 \end{lemma}

 Using \Cref{eq:weighted-poincare} in \Cref{eq:CZ-weighted} we further deduce the inequality
  \begin{equation}\label{eq:CZ-weighted}
  \WN_{\Omega}^p(\mu,\lambda)   \les \frac{c_P(\eta,p)^p}{g_0^{p-1}} \nor{f-g}_{L^{p}(\eta)}^p.
  \end{equation}

 \begin{proof}
  By the triangle inequality and \Cref{eq:convexity-wass}, we find, for every $\eps \in (0,1)$,
  \begin{equation*}\begin{split}
   \WN_{\Omega}^p(\mu,\lambda)& \le  (1+\eps) \WN_{\Omega}^p\bra{\mu, (\mu+\lambda)/2 } + \frac{c}{\eps^{p-1}} \WN_{\Omega}^p\bra{(\mu+\lambda)/2, \lambda }\\
   & \le \frac{(1+\eps)}{2} \WN_{\Omega}^p\bra{\mu, \lambda } + \frac{ c}{ \eps^{p-1}}  \WN_{\Omega}^p\bra{(\mu+\lambda)/2, \lambda },
  \end{split}\end{equation*}
 and choosing e.g.\ $\eps = 1/4$ we deduce that
 \begin{equation*}
  \WN_{\Omega}^p(\mu,\lambda) \les  \WN_{\Omega}^p\bra{(\mu+\lambda)/2, \lambda }.
 \end{equation*}
 Hence, we may assume that $f(x), g(x) \ge g_0/2$ for every $x \in \Omega$. Let now $b\in L^p(\eta)$ be such that $-\div_\eta b=f-g$ (recall the definition \Cref{defdiveta} of $\div_\eta$). Setting $\rho_t=(1-t)\mu+t\lambda$ and $j_t=b\eta$ we see that the couple $(\rho_t,j_t)$ is admissible for the Benamou-Brenier formulation of $\WN_{\Omega}^p(\mu,\lambda)$, see e.g. \parencite{Santam} thus
\begin{equation*}
  \WN_{\Omega}^p(\mu,\lambda) \le \int_{\Omega} \int_0^1 \frac{1}{\rho_t^{p-1}}|j_t|^p dt\les  \frac{1}{g_0^{p-1}} \int_{\Omega} |b|^p \eta.
\end{equation*}
Minimizing over all such $b$ we conclude the proof by \Cref{eq:dual-neg-sob}.
 \end{proof}

 We are now in a position to state and prove the main result of this section, which provides a quantitative upper bound for $\WN^p$ in terms of $\WD^p$.

\begin{theorem}\label{thm:main-deterministic}
Let $d\ge 1$, $p>1$, $p'=p/(p-1)$, $c, r >0$, and
\begin{enumerate}
 \item   $\Omega \subseteq \R^d$ be convex and bounded,  $\eta: \R^d \to [0,1]$ be $C^2$ smooth and such that $\eta(x) = 1$ if $\dist(x, \Omega^c) \le r$, $\eta(x) = 0$ if $\dist(x, \Omega^c) \ge 2r$, $\nor{ \nabla \eta }_\infty \le c/r$ and $\nor{ \nabla^2 \eta}_\infty \le c/r^2$,
  \item $\mu$ be a measure on $\Omega$ supported on $\cY \subseteq \overline{\Omega}$ with
  \begin{equation*}
  \frac 1 2 \le \frac{\mu(\Omega)}{|\Omega|} \le 2, \quad \text{and set} \quad \lambda := \frac{\mu(\Omega)}{|\Omega|} \chi_{\Omega},
  \end{equation*}
  \item  $v \in \mathcal{S}(\bist_{\Omega}^p, \cY)$ be an optimizer in the dual formulation for $\WD_{\Omega}^p(\mu,\lambda)$,  and set
  \begin{equation*}
   u:= \sfQ_{\bist_{\Omega}^p}(v),  \qquad \ell := \bra{ \osc_{\cY}(v) + \dist(\Omega, \cY)^p}^{1/p}.
  \end{equation*}
  \end{enumerate}
If $r \ge 4\ell$, then, for every $\eps \in (0,1)$ it holds
   \begin{equation}\begin{split}\label{eq:main-deterministic}
 \WN^p_{\Omega} \bra{ \mu, \lambda  } - \WD_{\Omega}^p\bra{\mu, \lambda } &  \les  \eps \WD_{\Omega}^p\bra{\mu, \lambda } +  \frac 1 {\eps^{p-1}} \Bigg[  \WN^p_{\Omega} \bra{  \frac{\int_{\Omega} \eta d \mu }{\int_{\Omega}  \eta  }  \eta , \eta \mu} \\
 & \quad +  \diam(\Omega)^{p} \bra{ \abs{ \frac{\fint_{\Omega} \eta d\mu }{ \fint_{\Omega}  \eta }-1}^p  \int_{\Omega} \eta^p + r^{-2p} \bra{ \WD_{\Omega}^p(\mu, \lambda) }^{\frac{d+2p}{d+p}}  }\\
 & \quad +  \nor{ (\nabla u)^{\bra{p'/p}} \cdot \nabla \eta - \fint_{\Omega} (\nabla u)^{\bra{p'/p}} \cdot \nabla \eta}_{H^{-1,p}(\Omega)}^p \Bigg]
 \end{split}
\end{equation}
where the implicit constant depends on $c$, $d$, $p$, $\Omega$, (but not on $r$) and is invariant with respect to dilations of $\Omega$.
\end{theorem}

\begin{remark}\label{rem:eta-cube}
 Let us notice that in the case of $\Omega = (0,L)^d$ a cube of side length $L >0$, one can easily find $\eta= \eta_{d,L,\delta}$ satisfying assumption \emph{(1)} with $r := \delta L$, for any given $\delta \in (0,1)$ and with $c = c(d)< \infty$ depending on $d$ only. Indeed, it is sufficient to argue in the case $L=1$ and $d=1$, and define  a suitable $\xi_\delta:\R \to [0,1]$ even and such that $\xi_{\delta}(x)=1$ for $|x|  \le 1/2-2\delta$, $\xi_{\delta}(x) = 0$ for $|x| \ge 1/2-\delta$   and then set, for $x= (x_i)_{i=1}^d \in (0,L)^d$,
 \begin{equation*}
  \eta_{d,L, \delta}(x) := 1-\prod_{i=1}^d \xi_{\delta}(x_i/L-1/2).
 \end{equation*}
 A more explicit choice for $\xi_{\delta}$ will be given in \Cref{sec:random}.
\end{remark}

\begin{proof}
 We split the proof into multiple steps.

\setcounter{step}{0}
 \stepcounter{step} \noindent \emph{Step \thestep \, (Notation).} Let  $\pi \in \mathcal{C}(\mu, \lambda)$ be an optimizer in the primal formulation of $\WD_{\Omega}^p(\mu, \lambda)$. Write $\bra{\pi(\cdot|x)}_{x \in \Omega}$ for the regular conditional distribution of $y$ given $x$, with respect to  $\pi$, which is a Markov kernel, that we further use to define, for  $x \in \Omega$, the function
%
%
\begin{equation*}
  \pi(\eta)(x) := \int_{\Omega} \eta(y)  \pi(dy|x).
\end{equation*}
Since the marginal distribution of $y$ with respect to $\pi$ is $\mu$, we notice that
 \begin{equation*}
  \int_{\Omega} \pi(\eta) d \lambda = \int_{\Omega \times \Omega} \eta(y) d \pi(x,y) = \int_{\Omega} \eta d \mu.
 \end{equation*}
To keep notation simple, we write
\begin{equation*}
 \kappa := \frac{\int_{\Omega} \pi(\eta) \lambda} {\int_{\Omega} \eta \lambda } = \frac{\int_{\Omega} \eta d\mu }{\int_{\Omega}  \eta \lambda },
\end{equation*}
and
\begin{equation*}
 \xi(x):= 1- \eta (x), \qquad \pi(\xi) (x):= 1- \pi(\eta) (x)= \int_{\Omega } \xi(y) \pi(dy|x),
\end{equation*}
and notice that $\xi, \pi(\xi) \in [0,1]$. Moreover, we notice that the measure $\xi(y) \pi(dx, dy)$ on $\Omega \times \Omega$ is a coupling between $\xi \mu$ (the marginal law of $y$) and $\pi(\xi)  \lambda$,
i.e.,
\begin{equation}\label{eq:coupling-xi-pi}
 \xi(y) \pi(dx, dy)\in \mathcal{C}(\pi(\xi) \lambda, \xi \mu).
\end{equation}
Indeed, for every $E \subseteq \Omega$ Borel,
\begin{equation*}
 \int_{\Omega \times \Omega} \chi_E(x) \xi(y) d\pi(x, y) =\int_{E} \int_{\Omega } \xi(y) \pi(dy|x) d\lambda(x) = \int_E \pi(\xi) d\lambda.
\end{equation*}


We introduce two auxiliary measures on $\Omega$:
\begin{equation*}
 \mu_1 :=  \eta \mu +  \pi(\xi) \lambda, \quad \text{and} \quad \mu_2  :=   \kappa \eta \lambda  +  \pi(\xi) \lambda,  
\end{equation*}
such that $\mu_2$ is absolutely continuous with respect to the Lesbesgue measure. 
Let us notice that both measures have total mass $\mu(\Omega)$.


\stepcounter{step} \noindent \emph{Step \thestep \, (Reduction to three bounds).}
We claim that the thesis follows from two applications of \Cref{eq:triangle-wass} and the following three inequalities:
\begin{eqnarray}
 \label{eq:main-01}  \WN_{\Omega}^p\bra{ \mu, \mu_1 }  & \le &  \WD_{\Omega}^p\bra{ \mu, \lambda },\\
\label{eq:main-12}    \WN_{\Omega}^p\bra{  \mu_1, \mu_2 }  &  \le  &  \WN_{\Omega}^p \bra{  \kappa  \eta \lambda, \eta \mu} 
\end{eqnarray}
and
\begin{equation}\begin{split}
 \label{eq:main-23} \WN_{\Omega}^p\bra{ \mu_2, \lambda }  & \les \diam(\Omega)^p \sqa{ |\kappa-1|^p  \int_{\Omega} \eta^p +  r^{-2p} \bra{ \WD_{\Omega}^p(\mu, \lambda) }^{\frac{d+2p}{d+p}} } \\
  & \quad +  \nor{(\nabla u)^{(p'/p)} \cdot \nabla \eta - \fint_{\Omega}(\nabla u)^{(p'/p)}\cdot \nabla \eta}_{H^{-1,p}(\Omega)}^p.
 \end{split}
\end{equation}
Indeed,  for every $\eps \in (0,1)$, it holds
\begin{equation*}\begin{split}
 \WN_{\Omega}^p\bra{ \mu, \lambda } & \le (1+\eps) \WN_{\Omega}^p\bra{ \mu, \mu_1 }  +\frac{c_1}{\eps^{p-1}}  \WN_{\Omega}^p\bra{ \mu_1, \lambda }\\
 & \le (1+\eps) \WN_{\Omega}^p\bra{ \mu_1, \mu } + \frac{c_2}{\eps^{p-1}}\sqa{ \WN_{\Omega}^p\bra{ \mu_1, \mu_2}  + \WN_{\Omega}^p\bra{ \mu_2, \lambda} },
 \end{split}
\end{equation*}
where the second inequality follows from \Cref{eq:triangle-wass} with $\eps=1/2$, and the constant $c_1, c_2<\infty$ depend on $p$  only and may be  different in both lines (and are unrelated to the constant $c$ in the assumptions). Inserting \eqref{eq:main-01}, \eqref{eq:main-12}, and  \Cref{eq:main-23} yields \Cref{eq:main-deterministic}.

We then establish separately the validity of the three claimed inequalities. 

\stepcounter{step} \noindent \emph{Step \thestep \, (Proof of \eqref{eq:main-01}).} Writing $\mu = \eta \mu + \xi \mu$,  by \Cref{eq:convexity-wass} we have
\begin{equation*}\begin{split}
 \WN_{\Omega}^p \bra{ \mu, \mu_1} & = \WN_{\Omega}^p \bra{\eta \mu + \xi \mu,  \eta \mu +  \pi(\xi) \lambda } \\
 & \le \WN_{\Omega}^p \bra{\xi \mu,  \pi(\xi) \lambda }.
 \end{split}
 \end{equation*}
 To further bound from above, we use the coupling  \Cref{eq:coupling-xi-pi} so that
 \begin{equation*}
  \WN_{\Omega}^p \bra{\xi \mu,  \pi(\xi) \lambda } \le \int_{\Omega \times \Omega} |y-x|^p \xi(y) d \pi(x,y).
 \end{equation*}
 The assumption $r \ge 4\ell$ entails in particular that $\xi(y) = 0$ if $y\notin \Omega^{(\ell)}$, therefore by
 \Cref{eq:dist-omega-dist-p-l-infty-bound} we deduce that in the integrand above we can replace $|y-x|$ with $\bist_{\Omega}(x,y)$.  Therefore,
 \begin{equation*}\begin{split}
  \int_{\Omega \times \Omega} |y-x|^p \xi(y) d \pi(x,y) & =  \int_{\Omega \times \Omega} \bist_{\Omega}(x,y)^p \xi(y) d \pi(x,y)\\
  & \le \int_{\Omega \times \Omega} \bist_{\Omega}(x,y)^p d \pi(x,y)\\
  & = \WD_{\Omega}^p(\mu, \lambda)
  \end{split}
 \end{equation*}
 by optimality of $\pi$. Thus, \eqref{eq:main-01} is proved.

\stepcounter{step} \noindent \emph{Step \thestep \, (Proof of \eqref{eq:main-12}).}  This inequality is even simpler than the previous one, because, again by \Cref{eq:convexity-wass} we have that
\begin{equation*}\begin{split}
 \WN_{\Omega}^p \bra{ \mu_2, \mu_1} &=   \WN_{\Omega}^p \bra{ \kappa  \eta \lambda + \pi(\xi) \lambda, \eta \mu  + \pi(\xi) \lambda } \\
 & \le \WN_{\Omega}^p \bra{  \kappa  \eta \lambda, \eta \mu}.
 \end{split}
 \end{equation*}

\stepcounter{step} \noindent \emph{Step \thestep \, (Proof of \Cref{eq:main-23}).}
%
 We apply \Cref{lem:peyre} on $\Omega$ with $\mu_2$ in place of $\mu$ (there) and  $g_0 = 1/2$: 
\begin{equation*}
 \WN_{\Omega}^p\bra{  \mu_2, \lambda } \les \nor{ \mu_2 - \lambda }_{H^{-1,p}(\Omega) }^p .
\end{equation*}
Writing $\lambda = \pi(\eta) \lambda + \pi(\xi)  \lambda$, we find that
\begin{equation*}
 \mu_2 - \lambda =  \kappa \eta \lambda -   \pi(\eta) \lambda = \bra{ \kappa \eta - \pi(\eta)} \frac{\mu(\Omega)}{|\Omega|}.
\end{equation*}
Since $\mu(\Omega)/|\Omega| \le 2$, we need only to show that
\begin{equation}\label{eq:main-23-only-this}\begin{split}
 \nor{ \kappa \eta  -   \pi(\eta) }_{H^{-1,p}(\Omega))}^p & \les \diam(\Omega)^p  \sqa{ |\kappa-1|^p  \int_{\Omega} \eta^p +  r^{-4} \bra{ \WD_{\Omega}^p(\mu, \lambda) }^{\frac{d+2p}{d+p}} } \\
 & \quad +  \nor{(\nabla u)^{(p'/p)}\cdot \nabla \eta - \fint_{\Omega}(\nabla u)^{(p'/p)}\cdot \nabla \eta}_{H^{-1,p}(\Omega)}^p.
 \end{split}
\end{equation}
To this aim, recall that the function
\begin{equation*}
 \R^d \ni x \mapsto u(x) = \sfQ_{\bist_\Omega^p}(v)(x)
\end{equation*}
is Lipschitz continuous, hence differentiable Lebesgue a.e.\ on $\Omega$. We use this fact to introduce the auxiliary function, defined for Lebesgue a.e.\ $x \in \Omega$,
\begin{equation*}
\rho(x) :=  \pi(\xi)(x) - \xi (x) +   \bra{ \frac {\nabla u(x)}{p}}^{\bra{p'/p}}  \cdot \nabla \xi (x)  = \eta (x) -  \pi(\eta)(x) - \bra{ \frac {\nabla u(x)}{p}}^{\bra{p'/p}} \cdot \nabla  \eta(x)
\end{equation*}
where the second identity follows from the definition of $\xi$. We collect in particular the identity
\begin{equation*}
 \pi(\eta) =  \eta - \bra{ \frac {\nabla u}{p}}^{\bra{p'/p}} \cdot \nabla \eta  - \rho.
\end{equation*}
Hence, by the triangle inequality, we find
\begin{equation*}\begin{split}
 \nor{ \kappa \eta  -   \pi(\eta) }_{H^{-1,p}(\Omega)}^p  & = \nor{( \kappa -1)\eta  + \bra{ \frac {\nabla u}{p}}^{\bra{p'/p}} \cdot \nabla \eta  + \rho}_{H^{-1,p}(\Omega)}^p \\
 & \les   \nor{(\kappa-1) \eta   - \fint_{\Omega}(\kappa-1)\eta  }_{H^{-1,p}(\Omega)}^p  \\
 & \quad + \nor{ { \nabla u }^{\bra{p'/p}} \cdot \nabla \eta   - \fint_{\Omega} { \nabla u }^{\bra{p'/p}} \cdot \nabla \eta } _{H^{-1,p}(\Omega) }^p \\
 &  \quad + \nor{ \rho -  \fint_{\Omega} \rho} _{H^{-1,p}(\Omega) }^p
  \end{split}
\end{equation*}
and we bound separately the three terms.  For the first one, we use \Cref{eq:poincare-negative}, so that
\begin{equation*}
 \nor{(\kappa-1) \eta    - \fint_{\Omega}(\kappa-1)\eta  }_{H^{-1,p}(\Omega)}^p \les \diam(\Omega)^p |\kappa-1|^p  \int_{\Omega} \eta^p
\end{equation*}
which we recognize as a contribution in  \Cref{eq:main-23-only-this}.   The second term already appears in \Cref{eq:main-23-only-this}, hence we focus on the third term.

We use again \Cref{eq:poincare-negative}, so that
\begin{equation*}
 \nor{ \rho -  \fint_{\Omega} \rho } _{H^{-1,p}(\Omega) }^p \les \diam(\Omega)^p \int_{\Omega} \rho^p.
\end{equation*}
We see that the proof of \Cref{eq:main-23-only-this}, hence \Cref{eq:main-23} and therefore of the entire result is completed if we argue that
\begin{equation}\label{eq:final-step-main-deterministic}
 \int_{\Omega} \rho^p \les r^{-2p} \bra{ \WD_{\Omega}^p(\mu, \lambda) }^{\frac{d+2p}{d+p}}.
\end{equation}

\stepcounter{step} \noindent \emph{Step \thestep \, (Proof of \Cref{eq:final-step-main-deterministic})}.
We notice that $\rho(x) = 0$ for Lebesgue a.e.\ $x \notin \Omega^{(r/2)}$. Indeed, by construction we have $\xi(x) = 0$ and $\nabla \xi (x) = 0$ for such $x$'s, and furthermore since $\xi(y) \le \chi_{\Omega^{(r)}}(y)$, we have that
\begin{equation*}
\begin{split}  \int_{ \Omega\setminus \Omega^{(r/2)}} \pi(\xi) (x)  d\lambda (x) &  = \int_{\Omega \times \Omega} \chi_{\Omega\setminus \Omega^{(r/2)}}(x)  \xi(y) d\pi(x,y) \\
& \le  \int_{\Omega \times \Omega} \chi_{\Omega\setminus  \Omega^{(r/2)}}(x) \chi_{\Omega^{(r)} }(y) d\pi(x,y) =0,
\end{split}
\end{equation*}
because for $\pi$-a.e.\ $(x,y)$, if $y \in \Omega^{(r)}$  by \Cref{eq:dist-omega-dist-p-l-infty-bound} we have
\begin{equation*}
 |x-y| \le  \ell \le \frac{r}{4},
\end{equation*}
and therefore, by triangle inequality,
\begin{equation*}
 \dist(x, \Omega^c) \ge \dist(y, \Omega^c) - |x-y| \ge r - \frac{r}{4} > \frac r 2.
\end{equation*}

Moreover, we have again by \Cref{eq:brenier}  that for $\pi$-a.e.\ $(x,y)$, if $x \in \Omega^{(r/2)}$, the conditional distribution $\pi(\cdot|x)$ reduces to a Dirac measure at $y = x-   \bra{ \frac{\nabla u(x)}{p}}^{(p'/p)}$, hence
\begin{equation*}
 \pi(\xi) (x) = \xi \bra{ x - \bra{ \frac{\nabla u(x)}{p}}^{(p'/p) }}.
\end{equation*}
Therefore, for such $x$'s we recognize that  $\rho(x)$ is the remainder of a first order expansion of $\xi$ centered at $x$ and evaluated at $y$ and therefore we can bound
\begin{equation*}
  \abs{\rho(x) } \les \nor{ \nabla^2 \xi }_\infty |y-x|^2 \les |y-x|^2r^{-2}
\end{equation*}
having used that $\nabla^2 \xi = - \nabla^2 \eta$. Finally, since $r/2 \ge 2\ell$, using \Cref{eq:l-infty-lp-bound-dirichlet} we can bound from above
\begin{equation*}
 |y-x|^p \les  \bra{ \WD_{\Omega}^p\bra{\mu, \lambda}}^{p/(d+p)}.
 \end{equation*}
Thus,
\begin{equation*}
\begin{split}
 \int_{\Omega } \rho^p d\lambda & = \int_{\Omega^{(r/2)}}  \rho^p  d\lambda \\
 & \les r^{-2p}  \bra{ \WD_{\Omega}^p\bra{\mu, \lambda}}^{p/(d+p)} \int_{\Omega^{(r/2)}\times \Omega} |y-x|^p d\pi(x,y) \\
  & \les r^{-2p}  \bra{ \WD_{\Omega}^p\bra{\mu, \lambda}}^{p/(d+p)} \int_{\Omega^{(r/2)}\times \Omega} \bist_{\Omega}^p (x,y) d\pi(x,y) \\
 & \le r^{-2p} \bra{ \WD_{\Omega}^p(\mu, \lambda) }^{\frac{d+2p}{d+p}},
 \end{split}
\end{equation*}
and the proof of \Cref{eq:final-step-main-deterministic} is completed.
\end{proof}

\subsection{Stability and the quadratic case}

 Moving towards the end of this section, we specialize to the case $p=2$, so that $p'/p =1$ and we are in a position to apply \Cref{lem:magic-identity} and bound from above the last term in \Cref{eq:main-deterministic} by
 \begin{equation}\label{H-12magic}
  \nor{ \nabla u \cdot \nabla \eta - \fint_{\Omega} \nabla u \cdot \nabla \eta}_{H^{-1,2}(\Omega)}^2  \les \bra{ r^{-2}  + \diam(\Omega)^2 r^{-4} } \Var_{\Omega^{(r)}}(u).
\end{equation}
  As already stated in the introduction, the problem is now  to bound from above $\Var_{\Omega^{(r)}}(u)$. For this we rely on recent progress on the uniform convexity properties of the Kantorovich functional from \cite{delalande2023quantitative} (in the form presented in \parencite{delalande2022quantitative}). While a similar stability result probably holds also for $\WD^2_{\Omega}$, see e.g. \parencite{mischler2024,chizat2024sharper}, we use instead that thanks to \Cref{lem:uniform-upper-bounds}, the optimal transport plan for $\WD^2_{\Omega}(\mu,\lambda)$, is also optimal for the usual quadratic cost when restricted to $\Omega^{(r)}$ so that \parencite{delalande2022quantitative} applies.


\begin{theorem}\label{thm:delalande}
Let $\Omega \subseteq \R^d$ be a bounded convex set and   $\mu$ be a measure supported on $\cY \subseteq \overline{\Omega}$ with
\begin{equation}\label{eq:equivalence-mu-omega}
   \frac 1 2  \le \frac{\mu(\Omega)}{|\Omega|} \le 2, \quad \text{and set} \quad
  \lambda := \frac{ \mu(\Omega)}{|\Omega|}\chi_{\Omega}.
 \end{equation}
 Then, for every $v \in \cS(\bist_{\Omega}^2, \cY)$ optimizer for the dual formulation of $\WD_{\Omega}^2(\mu, \lambda)$, and every $\tilde v \in \cS(\bist_\Omega^2, \cY)$, it holds
 \begin{equation}\label{eq:taylor-entropic-delalande-mu-vt-main-body}
  \Var_{  \Omega^{(3\ell)} } (u - \tilde u) \les \diam(\Omega)^2  \sqa{\WD^{2} _{\Omega}(\mu, \lambda ) - \bra{ \int_{\Omega } \tilde u d \lambda - \int_{\cY} \tilde v d \mu }} 
 \end{equation}
 where $u := \sfQ_{\bist_{\Omega}^2}(v)$, $\tilde u := \sfQ_{\bist_{\Omega}^2}(\tilde v)$ and
 \begin{equation}\label{defell}
    \ell := \bra{ \max\cur{\osc_{\cY}\bra{v}, \osc_{\cY}\bra{\tilde v}} + \dist(\cX, \cY)^2 }^{1/2}.
 \end{equation}
\end{theorem}
\begin{proof}
Up to multiplying $\mu$ and $\lambda$ by $|\Omega|/\mu(\Omega)$ we may assume that $\mu(\Omega)=|\Omega|$ and thus $\lambda=\chi_{\Omega}$. We start by noting that since $v,\tilde v \in \cS(\bist_\Omega^2, \cY)$, we may apply \Cref{lem:uniform-upper-bounds}
to both functions and get that for $v'\in \{v,\tilde{v}\}$ and $u':= \sfQ_{\bist_{\Omega}^2}(v')$ we have for $x\in \Omega^{(\ell)}$
\[
 u'(x)=\inf_{y\in \cY\cap \Omega} v'(y)+ \bist_\Omega^2(x,y)=\inf_{y\in \cY\cap \Omega^{(\ell)}} v'(y)+ |x-y|^2
\]
so that $u'$ is $\dist^2-$concave in $\Omega^{(\ell)}$. Similarly, $v'$ is also $\dist^2-$concave in $\Omega^{(\ell)}$ and $u'=   \sfQ_{\dist^2}(v')$ in $\Omega^{(2\ell)}$. Let then $\pi$ be an optimizer for $\WD^{2} _{\Omega}(\mu, \lambda )$ and define $\hat{\pi}$ by
\[
 \int \xi(x,y) d\hat{\pi}(x,y)=\int \xi(x,y) 1_{\Omega^{(3\ell)}}(x) d\pi(x,y).
\]
Notice that this is well-defined since the first marginal of $\pi$ is absolutely continuous with respect to the Lebesgue measure. By definition the first marginal of $\hat{\pi}$ is $1_{\Omega^{(3\ell)}}$. We let $\hat{\mu}$ be its second marginal. Since for $(x,y)\in \supp \hat{\pi}\subset \supp \pi$ we have $x\in \partial_{\bist^2_\Omega} v(y) \cap \Omega^{(3\ell)}\subset \partial v(y)$, we get that $\hat{\pi}$ is an optimizer for $\WN^2(\hat{\mu},1_{\Omega^{(3\ell)}})$. Moreover using once more that on $ \supp \hat{\pi}$ we have $x\in  \partial v(y)$ we also have that $v$ is an optimizer for the dual version of $\WN^2(\hat{\mu},1_{\Omega^{(3\ell)}})$.

Thanks to \cite[Corollary 1.31]{delalande2022quantitative}, with
\begin{multline*}
 \phi^0(x):= \frac{1}{2} (|x|^2-u(x)), \quad \psi^0(y):=\frac{1}{2}(|y|^2+v(y)), \qquad \textrm{and} \\  \phi^1(x):= \frac{1}{2} (|x|^2-\tilde{u}(x)), \quad \psi^1(y):=\frac{1}{2}(|y|^2+\tilde{v}(y))
\end{multline*}
we have
\begin{equation*}\begin{split}
 \Var_{\Omega^{(3\ell)} }(u -  \tilde u) & =\Var_{\Omega^{(3\ell)} }(\phi^0 -  \phi^1)\\ & \les \diam(\Omega)^2  \sqa{ \int_{\Omega^{(3\ell)}}  \phi^1 dx -\int_{\Omega^{(3\ell)}}  \phi^0 dx  + \langle \psi^1-\psi^0, (\nabla \phi^0)\# 1_{\Omega^{(3\ell)}}\rangle}.
 \end{split}
\end{equation*}
By \Cref{eq:brenier}, we have $(\nabla \phi^0)_\sharp 1_{\Omega^{(3\ell)}}=\hat{\mu}$ (where $\sharp$ denotes the push-forward) and thus
\begin{equation*}
 \begin{split}
   \int_{\Omega^{(3\ell)}}  \phi^1 dx -\int_{\Omega^{(3\ell)}}  \phi^0 dx  + \langle \psi^1-\psi^0, (\nabla \phi^0)_\sharp 1_{\Omega^{(3\ell)}}& =\int \lt[\phi^1(x)+\phi^0(x) + \psi^1(y)-\psi^0(y)\rt] d\hat{\pi}(x,y)\\
& =\int \lt[u(x)-\tilde u(x) + \tilde v(y)-v(y)\rt] d\hat{\pi}(x,y).
 \end{split}
\end{equation*}
Since $\hat{\pi}\le \pi$ and since on $\supp \pi$,
\[
 u(x)-\tilde u(x) + \tilde v(y)-v(y)=\bist_\Omega^2(x,y) -(\tilde u(x)-\tilde v(y))\ge 0
\]
we finally obtain
\begin{equation*}
 \begin{split}
   \Var_{\Omega^{(3\ell)} }(u -  \tilde u) &\les \diam(\Omega)^2 \int \lt[\bist_\Omega^2(x,y) -(\tilde u(x)-\tilde v(y))\rt] d\pi(x,y)\\
 & =\diam(\Omega)^2\sqa{\WD^{2} _{\Omega}(\mu, \lambda ) - \bra{ \int_{\Omega } \tilde u d \lambda - \int_{\cY} \tilde v d \mu }}.
 \end{split}
\end{equation*}
This concludes the proof of \Cref{eq:taylor-entropic-delalande-mu-vt-main-body}.
\end{proof}

We apply the above result to strengthen the super-additivity inequality \Cref{eq:Dirichlet-sub-additive}.

\begin{corollary}\label{cor:stability-superadd-1}
 Let  $\Omega \subseteq \R^d$ be a bounded convex set and  $\mu$ be a measure supported on $\cY \subseteq \overline{ \Omega}$  satisfying \Cref{eq:equivalence-mu-omega}. Consider a finite partition $\Omega = \bigcup_{k} \Omega_k$ such that $\mu(\Omega_k)>0$ for every $k$ 
and let $v \in \cS(\bist_{\Omega}^2, \cY)$ be an optimizer for the dual formulation of $\WD_{\Omega}^2(\mu, \mu(\Omega)/|\Omega|)$. Then,  for every $\eps \in (0,1)$, it holds

 \begin{equation}\label{eq:taylor-delalande-sub-additivity}
\begin{split}
   \Var_{  \Omega^{(3\ell)} }\bra{ \sfQ_{\bist_{\Omega}^2}(v) }& \les  \diam(\Omega)^2 \Bigg[ \WD_{\Omega}^2\bra{ \mu, \frac{\mu(\Omega)}{|\Omega|}} -  (1+\eps)^{-1} \sum_k  \WD_{\Omega_k}^2\bra{\mu,   \frac{\mu(\Omega_k)}{|\Omega_k|} } \\
  & \quad + \eps^{-1}\WN^2\bra{ \sum_{k} \frac{\mu(\Omega_k)}{|\Omega_k|} \chi_{\Omega_k}, \frac{\mu(\Omega)}{|\Omega|}}  \Bigg]  +   |\Omega| \max_k\diam(\Omega_k)^4,
   \end{split}
 \end{equation}
 with
 \begin{equation}
  \ell :=2 \max_{k} \diam(\Omega_k)
 \end{equation}
\end{corollary}

\begin{proof}
We apply \Cref{thm:delalande} on $\cX = \overline{\Omega}$, with $\tilde v$ as in \Cref{eq:tilde-v-superadd} and consequently $\tilde u:= \sfQ_{\bist_{\Omega}^2}(\tilde v)$ so that  in particular \Cref{eq:norm-tilde-v} holds. Since $\mu(\Omega_k)>0$ for every $k$, we deduce that $\cY \cap \Omega_k \neq \emptyset$ for every $k$, hence
\begin{equation}\label{osctildeu}
 \dist(\Omega, \cY) \le \max_k \diam(\Omega_k) \quad \text{and} \quad \max\cur{\osc_{\cY}(\tilde v), \osc_{\Omega}(\tilde u) }\le 3\max_{k} \diam(\Omega_k)^2.
 \end{equation}
Using \Cref{eq:taylor-entropic-delalande-mu-vt-main-body} and \Cref{eq:potential-subadditive-new} with $\lambda = \mu(\Omega)/|\Omega|$ (recall that we identify absolutely continuous measures with their densities) we find
\begin{equation*}\begin{split}
 \Var_{\Omega^{(3\ell)}}(u-\tilde u) & \les \diam(\Omega)^2 \sqa{ \WD_{\Omega}^2(\mu, \lambda) - \sqa{\int_{\Omega}\sfQ_{\bist_{\Omega}^2}( \tilde v) d \lambda -  \int_{\cY} \tilde v d \mu } }\\
 & \les \diam(\Omega)^2 \Bigg[ \WD_{\Omega}^2(\mu, \lambda ) - (1+\eps)^{-1} \sum_k  \WD_{\Omega_k}^2\bra{\mu,   \frac{\mu(\Omega_k)}{\lambda(\Omega_k)} \lambda} \\
 & \quad \quad + \eps^{-1}\WN^2\bra{  \sum_k\frac{\mu(\Omega_k)}{\lambda(\Omega_k)} \chi_{\Omega_k}\lambda, \lambda} \Bigg].
 \end{split}
\end{equation*}
The thesis then follows from the triangle inequality for the variance and \Cref{eq:variance-osc}:
\begin{equation*}\begin{split}
 \Var_{\Omega^{(3\ell)}}(u) & \les \Var_{\Omega^{(3\ell)}}(u-\tilde u)  + \Var_{\Omega^{(3\ell)}}(\tilde u)  \les \Var_{\Omega^{(3\ell)}}(u-\tilde u) + |\Omega|\osc_{\Omega}\bra{\tilde u}^2\\
 & \stackrel{\Cref{osctildeu}}{\les} \Var_{\Omega^{(3\ell)}}(u-\tilde u) + |\Omega| \max_{k}\diam(\Omega_k)^4. \qedhere
 \end{split}
\end{equation*}
\end{proof}

\section{Application to random optimal transport}\label{sec:random}

The aim of this section is to prove \Cref{thm:main-poisson}. We actually establish a more general version   which applies to random Borel measures, under the following assumptions.

\begin{assumption}\label{ass:mu}Let $\mu$ denote a random Borel measure on $\R^d$ such that
\begin{enumerate}[i)]
\item (stationarity) for every $v \in \R^d$,  it holds $\theta_v \mu  =  \mu$ in law, where $\theta_v \mu(A) = \mu(A+v)$ for every $A$ Borel,
\item  (integrability) for every bounded Borel $A \subseteq \R^d$, $\mu(A)$ is integrable.
\item (concentration) there exists $\alpha\in [0, d)$ such that for every cube $Q \subseteq \R^d$ with $\diam(Q) \ge 1$ and bounded Borel $\eta: Q \to [0,\infty)$, it holds for every $q \ge 1$,
\begin{equation}\label{eq:concentration-abstract}
\nor{ \int_Q \eta d\mu - \EE\sqa{ \int_Q \eta d\mu }}_{L^q(\PP)} \les_q  \diam(Q)^{(d+\alpha)/2}\nor{\eta}_\infty.
\end{equation}
\end{enumerate}
\end{assumption}

Notice that, by stationarity and integrability, the function $A \mapsto \EE\sqa{ \mu(A)}$ is a translation invariant ($\sigma$-finite) measure, hence for some constant $c \in [0, \infty)$ it holds
\begin{equation*}
\EE\sqa{ \mu(A) } = c | A|.
\end{equation*}
Without loss of generality, we assume in what follows that $c=1$, so that $\EE\sqa{ \mu(A)} = |A|$ (if $c=0$ the statements become trivial).

Notable examples of random measures $\mu$ satisfying \Cref{ass:mu} include:
\begin{enumerate}[a)]
 \item The Poisson point process, that can be defined as the limit $\mu = \lim_{r \to \infty} \mu_r$, where
 \begin{equation*}
  \mu_r := \sum_{i=1}^{N_r}\delta_{X_i},
 \end{equation*}
 where $N_r$ denotes a Poisson random variable with  $\EE\sqa{N_r} = |\B(r)|$ and $(X_i)_{i=1}^\infty$ are i.i.d.\ (and independent from $N_r$) uniformly distributed on the ball $\B(r)$. In this case, \Cref{eq:concentration-abstract} holds with $\alpha = 0$.
 \item The Brownian interlacement occupation measure \parencite{mariani2023wasserstein}, that can be defined similarly as a limit $\mu_\infty = \lim_{r \to \infty} \mu_r$, where in this case
 \begin{equation*}
  \mu_r := \sum_{i=1}^{N_r} \int_0^\infty \delta_{B^i_s +X_i} ds
 \end{equation*}
 and $N_r$ denotes a Poisson random variable with $\EE\sqa{N_r} = \Cap(\B(r))$ (the Newtonian capacity) $(X_i)_{i=1}^\infty$ are uniformly distributed on the sphere $\partial \B(r)$ and $(B^i)_{i=1}^\infty$ are independent standard Brownian motions on $\R^d$ -- and all the listed random variables are independent. In this case \Cref{eq:concentration-abstract} holds with $\alpha = 2$ as proved in \Cref{sec:ass-mu} below.
\end{enumerate}

Generalizing  ideas and tools from \parencite{BaBo, goldman2021convergence, ambrosio2022quadratic, goldman2022optimal}, in \parencite[Proposition A.2]{mariani2023wasserstein} the following existence result is established.

 \begin{proposition}
Let $\mu$ be a random measure on $\R^d$ satisfying i), ii), iii) in \Cref{ass:mu}  with $d>2+\alpha$. Then, for every $p \ge 1$, the following limit exists:
\begin{equation*}
\lim_{L \to \infty}  \frac{1}{L^d} \EE\sqa{ \WN_{(0,L)^d}^p\bra{ \mu , \frac{ \mu( (0,L)^d)}{L^d} } } \in [0, \infty).
\end{equation*}
 \end{proposition}


In particular, for every $p \ge 1$ it holds
\begin{equation*}
 \EE\sqa{ \WN_{(0,L)^d}^p\bra{ \mu , \frac{ \mu( (0,L)^d)}{L^d} } } \les L^d.
\end{equation*}

For $\WD$ a similar statement can be obtained with a much simpler proof, see   \parencite[Proposition A.3]{pieroni-gaussian} in the case of the Poisson point process. For $p\ge 1$, the following limit exists:
\begin{equation*}
\lim_{L \to \infty}  \frac{1}{L^d} \EE\sqa{ \WD_{(0,L)^d}^p\bra{ \mu , \frac{ \mu( (0,L)^d)}{L^d} } } \in [0, \infty).
\end{equation*}
Our main result is that  if \Cref{ass:mu} holds then the two limits coincide when $p=2$.

\begin{theorem}\label{thm:main-general}
 Let $\mu$ be a random measure on $\R^d$ satisfying \Cref{ass:mu} with
 \begin{equation*}
  d > 2 + \alpha.
 \end{equation*}
Then, the following two limits (exist and) coincide:
 \begin{equation}\label{eq:main-general}
  \lim_{L \to \infty} \frac{1}{L^d} \EE\sqa{ \WN_{(0,L)^d}^2\bra{ \mu, \frac{ \mu((0,L)^d)}{L^d}  }} = \lim_{L \to \infty} \frac{1}{L^d} \EE\sqa{ \WD_{(0,L)^d}^2\bra{ \mu, \frac{ \mu((0,L)^d)}{L^d} }}.
 \end{equation}
\end{theorem}

The  proof is an application of the deterministic estimate \Cref{eq:main-deterministic} from \Cref{thm:main-deterministic} combined with \Cref{H-12magic} and \Cref{eq:taylor-delalande-sub-additivity} from \Cref{cor:stability-superadd-1} on a suitable event. As a stochastic ingredient, we will use the crucial fact that by convergence of the right-hand side term in \Cref{eq:main-general}, the first line on the right-hand side of \Cref{eq:taylor-delalande-sub-additivity} is small and thus $ \Var_{  \Omega^{(3\ell)} }\bra{ \sfQ_{\bist_{\Omega}^2}(v) }$ is also small. We will also need the following two lemmas to bound from above the third right-hand side term in \Cref{eq:taylor-delalande-sub-additivity} and the second right-hand side  term  in \Cref{eq:main-deterministic}. The first lemma reads as follows.

\begin{lemma}\label{lemma:upper-bound-mu-intermediate-cubes}
 Let \Cref{ass:mu} hold with $d>2+\alpha$. For every $L_0 \ge 1$  and dyadic $m = 2^h \ge 2$, letting $L := m L_0$, it holds
 \begin{equation*}
 \frac 1 {L^d} \EE\sqa{\WN^{p}_{(0,L)^d}\bra{\frac{\mu((0,L)^d)}{L^d}, \sum_{k=1}^{m^d} \frac{\mu(Q_{L_0}^k)}{L_0^k} \chi_{Q_{L_0}^k}  } } \les L^{(2+\alpha-d)p/2}_0.
 \end{equation*}
 where $(0,L)^d = \bigcup_{k=1}^{m^d} Q_{L_0}^k$ is a partition into cubes of side length $L_0$. 
\end{lemma}

Before we state the second lemma, we make an explicit choice for the cut-off $\eta$, following  \Cref{rem:eta-cube}.
Let first
\begin{equation*}
 f(x)=10x^3-15x^4+6x^5=1+(x-1)^3(6x^2+3x+1)
\end{equation*}
so that 
\begin{equation}\label{condf}
 f(0)=f'(0)=f''(0)=f'(1)=f''(1)=0 \qquad \textrm{and} \qquad  f(1)=1.
\end{equation}
For $\delta>0$, we then consider $\xi_\delta:\R\to [0,1]$ even and defined on $[0,\infty)$ as
\begin{equation*}
 \xi_\delta(x)=\begin{cases}
                  1 & \text{for $x \in [0,1/2-2\delta]$,}\\
                  1- f\bra{ \frac{ x-(\frac{1}{2}-2\delta)}{\delta}} & \text{for $x \in [1/2-2\delta, 1/2-\delta]$,}\\
                  0 & \text{for $x\in [1/2-\delta,\infty)$.}
             \end{cases}
\end{equation*}
By \Cref{condf} we see that $\xi_\delta\in C^2(\R)$. Moreover,  $\delta \nor{ \xi'_\delta}_\infty+\delta^2 \nor{ \xi''_\delta}_\infty\les 1$. We then let, for $L>0$,
\begin{equation}\label{eq:eta-d-L-delta}
 \eta_{d,L,\delta}(x)=1-\prod_{i=1}^d \xi_\delta(x_i/L - 1/2),
\end{equation}
so that $\nor{\nabla \eta_{d,L,\delta}}_\infty \les (\delta L)^{-1}$,  $\nor{\nabla^2 \eta_{d,L,\delta}}_\infty \les (\delta L)^{-2}$, where the implicit constant depends on $d$ only. Notice for later use that, writing $\Gamma = \cur{x \in (0,L)^d \, : \, \eta_{d,L,\delta} (x) =0}$, then $\Gamma$ is a cube of side length $(1-2\delta)L$ and it holds for every $x \in (0,L)^d$,
 \begin{equation}\label{etasimdist}
  \bra{ \frac{\dist(x,\Gamma)}{\delta L}}^3 \les \eta_{d,L,\delta}(x) \les  \bra{ \frac{\dist(x,\Gamma)}{\delta L}}^3.
\end{equation}

%

\begin{lemma}\label{lemma:upper-bound-eta}
 Let \Cref{ass:mu} hold with $d>2+\alpha$. For every dyadic $\delta = 2^{-h_0}>0$, there exists $c(\delta)<\infty$ such that, for every $L \ge 1$,
 \begin{equation}\label{claimmatching}
 \frac{1}{L^d}\EE\lt[\WN_{(0,L)^d}^p\bra{\eta \mu,  \frac{ \int_{(0,L)^d}\eta d \mu}{\int_{(0,L)^d} \eta } \eta }  \rt]\les \delta + c(\delta) L^{(2+\alpha-d)p/2}.
 \end{equation}
 where  $\eta:= \eta_{d,\delta,L}$ as in \Cref{eq:eta-d-L-delta}.
\end{lemma}

We show first how to deduce \Cref{thm:main-general} using \Cref{lemma:upper-bound-mu-intermediate-cubes} and \Cref{lemma:upper-bound-eta}.

\begin{proof}[Proof of \Cref{thm:main-general}]\setcounter{step}{0}
We split the proof into several steps.

\stepcounter{step} \noindent \emph{Step \thestep \, (Notation).}
We write for brevity
\begin{equation*}
 \cost(L) := \frac{1}{L^d} \EE\sqa{ \WN_{Q_L}^2\bra{ \mu, \frac{ \mu(Q_L)}{L^d}  }}, \quad \costb(L) := \frac{1}{L^d} \EE\sqa{ \WD_{Q_L}^2\bra{ \mu, \frac{ \mu(Q_L)}{L^d}  }},
\end{equation*}
and let $\cost(\infty) = \lim_{L \to \infty} \cost(L)$, $\costb(\infty) = \lim_{L \to \infty} \cost(L)$. Since $\cost(L) \ge \costb(L)$ for every $L$, it is sufficient to establish the inequality $\cost(\infty) \le \costb(\infty)$. We introduce the quantity
\begin{equation*}
 \omega(L):= \sup_{L' \ge L} \abs{ \costb(\infty) - \costb(L') } \stackrel{L\to \infty}{\to} 0.
\end{equation*}

Without loss of generality we may assume that $L = 2^{\bar h}$ is dyadic. We write $\lambda$ for the constant density $\lambda = \mu(Q_L)/L^d \chi_{Q_L}$ on $Q_L:=(0,L)^d$.  Let $\gamma\in(0,1)$ to be chosen below and set $h:=\gamma \bar h$, $L_0:=2^h=L^\gamma$ and   $m:= L/L_0=2^{\bar h-h}=2^{(1-\gamma)\bar h}\in\N$. Notice that  $L_0 \to \infty$ as $L\to \infty$. With this notation,  $Q_L$ is partitioned up to Lebesgue negligible (hence $\PP$-a.s.\ $\mu$ negligible sets) into $m^d$ disjoint sub-cubes $(Q_{L_0}^k)_{k=1, \ldots, m^d}$, each of side length $L_0$.

Furthermore, we introduce auxiliary parameters $\delta = 2^{-h_0}$ for $h_0\in \N$, $\eps \in (0,1)$ and let $\eta:= \eta_{d,L, \delta}$ be as in \Cref{eq:eta-d-L-delta}, so that assumption \emph{(1)} in \Cref{thm:main-deterministic} holds on $\Omega = Q_L$ with $r:=\delta L$  and $c = c(d)<\infty$. We will let eventually $\delta \to 0$ and $\eps\to0$, but only after sending $L \to \infty$.

\stepcounter{step} \noindent \emph{Step \thestep \, (Reduction to a good event).}
We introduce the event
\begin{equation*}
  A:=  \cur{  \WD_{Q_L}^2\bra{ \mu, \lambda } \le L^{d+\eps}} \cap \cur{ \abs{\mu(Q_L) - L^d} \le L^{(1+\eps)(d+\alpha)/2}} \cap  \cur{ \mu(Q_{L_0}^k) >0, \, \forall k= 1, \ldots, m^d}. 
\end{equation*}
To estimate $\PP\bra{A^c}$ we use a union bound. First,  by \Cref{eq:concentration-abstract} we get for any $q>0$
\begin{equation}\label{eq:420}
\PP\bra{ \abs{\mu(Q_L) - L^d} \ge L^{(1+\eps)(d+\alpha)/2}} \les_q L^{-q\eps \bra{d+\alpha}/{2}}
\end{equation}
and
\begin{equation}\label{change422}
\sum_{k=1}^{m^d}\PP\bra{\mu(Q_{L_0}^k) = 0 }\les_q m^d L_0^{-q\bra{d-\alpha}/{2}}.
\end{equation}
For $L$ large enough, in the event $\cur{ \abs{\mu(Q_L) - L^d} \le L^{(1+\eps)(d+\alpha)/2}}$ (and in particular on $A$),  we have
\begin{equation*}
 \frac 1 2 \le \frac{\mu(Q_L)}{L^d} \le 2,
 \end{equation*}
so that by Markov inequality and \Cref{eq:wass-different-p}, we have for every $q'\ge2$,
\begin{equation*}\begin{split}
\PP\bra{ \cur{ \WD_{Q_L}^2\bra{ \mu, \lambda } \ge L^{d+\eps}} \cap \cur{ \abs{\mu(Q_L) - L^d} \le L^{(1+\eps)(d+\alpha)/2}} } &\les_{q'} \EE\sqa{ \WD_{Q_L}^{q'}\bra{ \mu, \lambda } } L^{-(d+\eps) q'/2} \\
& \les_{q'} L^{d-(d+\eps)q'/2}.
\end{split}
\end{equation*}
Combining the above with  \Cref{eq:420} and \Cref{change422}, this gives
  \begin{equation}\label{change426}
  P(A^c) \les_{q,q'} L^{d-(d+\eps) q'/2} + L^{-q\eps (d+\alpha)/2}+ m^d L_0^{-q(d-\alpha)/2},
 \end{equation}
 where the last term is then equal to
 \begin{equation*}
  m^d L_0^{-q(d-\alpha)/2}=2^{\bar h d -h(d+q\bra{d-\alpha}/2)}=2^{\bar h( d -\gamma(d+q\bra{d-\alpha}/{2}))}.
 \end{equation*}

On $A^c$, we trivially bound from above
 \begin{equation*}
 \WN_{Q_L}^2\bra{\mu, \lambda } \les_d \mu(Q_L)L^2
 \end{equation*}
 hence we have, by Cauchy-Schwarz inequality,
  \begin{equation}\label{change429}\begin{split}
  \frac{1}{L^d} \EE\sqa{ \chi_{A^c} \WN_{Q_L}^2\bra{\mu, \lambda }} & \les L^{2-d} \EE\sqa{\mu(Q_L)^2}^{1/2}P(A^c)^{1/2}\\
  & \les L^{2}\bra{ L^{d/2-(d+\eps) q' /4} + L^{-q\eps (d+\alpha)/4} + m^{d/2} L_0^{-q(d-\alpha)/4}  }=:\err_{\gamma}(L).
  \end{split}
 \end{equation}
 Notice that
 \begin{equation*}
  L^{2}m^{d/2} L_0^{-q(d-\alpha)/4}=2^{2 \bar h}2^{\bar h \frac{d}{2} -h(d+q\bra{d-\alpha}/{4})}=2^{\bar h\sqa{\bra{ 2+\frac{d}{2}}-\gamma \bra{ d+q\frac{d-\alpha}{4}}}}.
 \end{equation*}
 We now choose for every $\gamma,\eps\in(0,1)$ fixed, $q'$ and $q$  sufficiently large, such that
\begin{equation*}
  2 +d/2-(d+\eps) q' /4< 0, \quad  2-  q\eps \bra{d+\alpha}/4 <0,\quad \text{and} \quad \bra{2+\frac{d}{2}}-\gamma \bra{d+q\frac{d-\alpha}{4}}<0.
\end{equation*}
With this choice we get $\lim_{L\to \infty} \err_{\gamma}(L)=0$ and thus
\begin{equation*}
\limsup_{L\to \infty} \frac{1}{L^d} \EE\sqa{ \chi_{A^c} \WN_{Q_L}^2\bra{\mu, \lambda }}=0.
\end{equation*}

The thesis is reduced to
\begin{equation*}
 \lim_{L\to \infty}\frac{1}{L^d} \EE\sqa{ \chi_{A} \WN_{Q_L}^2\bra{\mu, \lambda }} \le \costb(\infty).
\end{equation*}
Let us notice for later use that, arguing similarly on each $Q_{L_0}^k$ and using \Cref{eq:concentration-abstract} we get
\begin{equation*}\begin{split}
 \frac{1}{L_0^d} \EE\sqa{ \chi_{A^c} \WD_{Q_{L_0}^k }^2\bra{\mu, \frac{\mu(Q_{L_0}^k)}{|L_0^d| } }} & \les L_0^2 L_0^{-d} \EE\sqa{\mu(Q_{L_0}^k)^2}^{1/2}P(A^c)^{1/2}\\
 & \stackrel{\Cref{change426}}{\les} L_0^{2} \bra{ L^{d/2-(d+\eps) q' /4} + L^{-q\eps (d+\alpha)/4} + m^{d/2} L_0^{-q(d-\alpha)/4}  }\\
 & \stackrel{L_0 \le L}{\le}\err_{\gamma}(L).
 \end{split}
\end{equation*}
 Hence,
 \begin{equation*}
  \frac{1}{L^d}\sum_{k=1}^{m^d} \EE\sqa{ \chi_{A^c} \WD_{Q_{L_0}^k }^2\bra{\mu, \frac{\mu(Q_{L_0}^k)}{|L_0^d| } }}=\frac{1}{m^d} \sum_{k=1}^{m^d} \frac{1}{L^d_0}\EE\sqa{ \chi_{A^c} \WD_{Q_{L_0}^k }^2\bra{\mu, \frac{\mu(Q_{L_0}^k)}{|L_0^d| } }} \les  \err_{\gamma}(L).
 \end{equation*}
Up to multiplying $\err_{\gamma}(L)$ by a small (universal) constant, we thus find
   \begin{equation}\label{eq:for-later-use-main}
    \frac 1 {L^d} \sum_{k=1}^{m^d} \EE\sqa{ \chi_A \WD_{Q_{L_0}^k}^2\bra{\mu, \frac{\mu(Q_{L_0}^k)}{L_0^k}} } \ge  \costb(L_0)- \err_{\gamma}(L). 
   \end{equation}

\stepcounter{step} \noindent \emph{Step \thestep \, (Variance bound).} On the event $A$, we apply \Cref{cor:stability-superadd-1} with $\cY =\operatorname{supp}(\mu) \cap \overline{\Omega}$ and obtain $v \in \cS(\bist_{Q_L}^2, \cY  )$ optimizer for the dual formulation of $\WD_{Q_L}^2(\mu, \lambda)$.  Writing $u := \sfQ_{\bist_{Q_L}^2}(v)$, we have by  \Cref{eq:osc-u-gagliardo} with $p=2$ (hence $p'=2$) and $\Omega = Q_L$, that
\begin{equation*}
   \osc_{Q_L}(u) \les_{d}   L^{1+\frac{ d+\eps}{d+2}},
  \end{equation*}
and since $\mu(Q_{L_0}^k)>0$ for every $k$,
\begin{equation*}
 \dist(Q_L, \cY) \les_d L_0.
\end{equation*}
By \Cref{oschatQ} we have
\begin{equation*}
 \ell=(\osc_{\cY}(v)+\dist(Q_L, \cY)^2)^{1/2}\le (\osc_{Q_L}(u)+2\dist(Q_L, \cY)^2)^{1/2}\les \max\cur{ L^{ \bra{ 1+  \frac{ d+\eps}{d+2}}/2 } , L_0 }.
\end{equation*}
Therefore, since $L_0 \ll L$ (which is guaranteed by our choice $L_0=L^\gamma$ for $\gamma\in (0,1)$) we have, for every given  $\eps, \delta \in (0,1)$,  that
\begin{equation}\label{rggl}
 r = \delta L \gg \ell.
\end{equation}
By \Cref{eq:taylor-delalande-sub-additivity} and \eqref{eq:var-tilde-lambda} on $\Omega = Q = \bigcup_{k} Q_{L_0}^k$, we obtain that for every $\eps' \in (0,1)$,
\begin{equation*}\begin{split}
  \Var_{  Q_L^{(r)} }(u)  & \les L^{2} \Bigg[ \WD^{2} _{Q_L}(\mu, \lambda ) - (1+\eps')^{-1}\sum_{k} \WD_{Q_{L_0}^k}^2\bra{\mu, \frac{\mu(Q_{L_0}^k)}{L_0^k}}  \\
  & \quad +  \frac 1 {\eps'}  \WN^{2} _{Q_L}\bra{\lambda, \sum_{k} \frac{\mu(Q_{L_0}^k)}{L_0^k} \chi_{Q_{L_0}^k}  }\Bigg]  +   L^d L_0^4.
  \end{split}
   \end{equation*}
Taking expectation on $A$, we find
   \begin{equation*}\begin{split}
    \frac {1}{L^d} \EE\sqa{ \chi_{A} \Var_{  Q_L^{(r)} }(u) } & \les  L^2 \bra{ \costb(L) - (1+\eps')^{-1} (\costb(L_0) -\err_\gamma(L)) + \frac 1 {\eps'} L_0^{2+\alpha-d} } + L_0^4 \\
    & \les L^2 \bra{ \omega(L_0) + \eps' + \err_\gamma(L) + \frac 1 {\eps'} L_0^{2+\alpha-d}  } + L_0^4
    \end{split}
   \end{equation*}
    having used  \Cref{eq:for-later-use-main} and \Cref{lemma:upper-bound-mu-intermediate-cubes} with $p=2$.
   The choice $\eps' = L_0^{(2+\alpha-d)/2}$ finally yields
    \begin{equation}\label{eq:bound-variance-A}
     \frac {1}{L^d} \EE\sqa{ \chi_{A} \Var_{  Q_L^{(r)} }(u) } \les L^2 \bra{ \omega(L_0) + \err_\gamma(L)+ L_0^{(2+\alpha-d)/2} } + L_0^4.
    \end{equation}

   \stepcounter{step} \noindent \emph{Step \thestep \, (Application of \Cref{thm:main-deterministic})}. In the event $A$,  thanks to \Cref{rggl} we may apply \Cref{thm:main-deterministic} together with Lemma \ref{lem:magic-identity} on $\Omega = Q_L$ to obtain that for every $\eps'' \in (0,1)$, 
   \begin{equation}\begin{split}\label{eq:main-proof-wn-wb}
 \WN_{Q_L}^2\bra{ \mu, \lambda} - \WD_{Q_L}^2\bra{\mu, \lambda }&  \les  \eps'' \WD_{Q_L}^2\bra{\mu, \lambda } +  \frac 1 {\eps''} \Bigg[  \WN_{Q_L}^2 \bra{  \frac{\int_{Q_L} \eta d \mu }{\int_{Q_L}  \eta  }  \eta , \eta \mu} \\
 & \quad +  L^{2} \bra{ \abs{ \frac{\int_{Q_L} \eta d\mu }{\int_{Q_L}  \eta  }-1}^2  \int_{Q_L} \eta^2 + r^{-4} \bra{ \WD_{Q_L}^2(\mu, \lambda) }^{\frac{d+4}{d+2}}  }\\
 & \quad +  \bra{ r^{-2}  +  L^{2}  r^{-4} } \Var_{Q_L^{(r)}}\bra{ u}  \Bigg].
 \end{split}
\end{equation}
After taking expectation (on $A$) and dividing both sides by $L^d$, we estimate  all the terms in the right-hand side, starting from the third line, which is bounded using \Cref{eq:bound-variance-A}:
\begin{equation}\label{mainestimterm1}
 \frac{1}{L^d}\bra{ r^{-2}  +  L^{2}  r^{-4} }  \EE\sqa{ \chi_{A}  \Var_{Q_L^{(r)}}\bra{ u} } \les \delta^{-4} \bra{ \omega(L_0) + \err_\gamma(L) + L_0^{(2+\alpha-d)/2} + L^{-2} L_0^4}.
\end{equation}
For fixed $\delta$ the right-hand side goes to zero provided
\begin{equation}\label{eq:second-condition-mL}
 \lim_{L \to \infty} L^{-2} L_0^4  = 0  \qquad \textrm{or equivalently } \qquad \gamma\in(0,1/2).
\end{equation}
For the second line in \Cref{eq:main-proof-wn-wb},  we have two contributions: the term
\begin{equation}\label{mainestimterm2}
   L^{2-d}  r^{-4} \bra{ \WD_{Q_L}^2(\mu, \lambda) }^{\frac{d+4}{d+2}}   \les \delta^{-4} L^{\frac{2d +\eps(d+4)}{d+2}-2},
\end{equation}
which is infinitesimal as $L \to \infty$, for $\eps$ small enough, and
\begin{equation*}\begin{split}
  L^{2-d}  \abs{ \frac{\int_{Q_L} \eta d\mu }{\int_{Q_L}  \eta }-1}^2  \int_{Q_L} \eta^2 &  = L^{2-d}  \abs{ \int_{Q_L} \eta d\mu - \int_{Q_L}  \eta }^2 \frac{ \int_{Q_L} \eta^2 }{\bra{ \int_{Q_L} \eta}^2 }\\
  & \les \delta^{-1} L^{2-2d} \abs{\int_{Q_L} \eta d\mu  - \int_{Q_L}\eta  }^2
  \end{split}
\end{equation*}
where the  inequality follows by the properties of $\eta$. We further bound in expectation, using \Cref{eq:concentration-abstract}:
\begin{equation}\label{mainestimterm3}\begin{split}
  \delta^{-1} L^{2-2d}\EE\sqa{ \chi_A \abs{\int_{Q_L} \eta d\mu  - \int_{Q_L}\eta }^2} & \le \delta^{-1} L^{2-2d}\Var_\PP\bra{ \int_{Q_L} \eta d\mu}\\
& \les \delta^{-1} L^{2+\alpha-d}
  \end{split}
\end{equation}
that is infinitesimal as $L \to \infty$, by the assumption $d>2+\alpha$. For the terms in the first line of the right-hand side of \Cref{eq:main-proof-wn-wb}, we have first that
\begin{equation}\label{mainestimterm4}
 \frac{\eps''}{L^d} \EE\sqa{ \chi_{A} \WD_{Q_L}^2\bra{\mu, \lambda } } \les \eps'' \costb(L)\les \eps'',
\end{equation}
hence we are  left with only one contribution, i.e., the second term in the first line, that we bound in expectation using \Cref{lemma:upper-bound-eta} with $p=2$:
 \begin{equation}\label{mainestimterm5}
  \frac 1 {L^d} \EE\sqa{ \WN_{Q_L}^2 \bra{\frac{\int_{Q_L} \eta d \mu }{\int_{Q_L}  \eta } \eta,  \eta \mu} } \les \delta + c(\delta)L^{2+\alpha-d}.
 \end{equation}
Plugging \Cref{mainestimterm1},\Cref{mainestimterm2},\Cref{mainestimterm3},\Cref{mainestimterm4} and \Cref{mainestimterm5} in \Cref{eq:main-proof-wn-wb} yields in combination with \Cref{change429},
\begin{equation}\label{conclusionmainestim}
 \frac 1 {L^d} \EE\sqa{\WN_{Q_L}^2\bra{ \mu, \lambda}} - \frac 1 {L^d} \EE\sqa{\WD_{Q_L}^2\bra{\mu, \lambda }}  \les  \eps'' +\frac{1}{\eps''}\lt[\delta + c(\delta)\bra{o_L(1) + o_{L_0}(1)+  \err_\gamma(L) + L^{-2} L_0^4}\rt].
\end{equation}

  \stepcounter{step} \noindent \emph{Step \thestep \, (Conclusion)}.
  We now choose $\gamma\in(0,1/2)$ arbitrarily. In particular, \Cref{eq:second-condition-mL} holds. Sending first $L\to \infty$ we get from \Cref{conclusionmainestim},
  \begin{equation*}
  \cost(\infty) - \costb(\infty)= \lim_{L\to \infty} \bra{ \frac 1 {L^d} \EE\sqa{\WN_{Q_L}^2\bra{ \mu, \lambda}} - \frac 1 {L^d} \EE\sqa{\WD_{Q_L}^2\bra{\mu, \lambda }}}  \les  \eps'' +\frac{\delta}{\eps''}.
  \end{equation*}
Sending first $\delta\to 0$ and then $\eps''\to 0$ (or directly choosing $\eps''=(\delta)^{1/2}$)  we obtain
\begin{equation*}
  \cost(\infty) - \costb(\infty)\le 0
\end{equation*}
which was the original claim.
\end{proof}

\subsection{Proof of \Cref{lemma:upper-bound-mu-intermediate-cubes}}

We write for brevity $\kappa(\Omega):= \mu(\Omega)/|\Omega|$. Writing $Q_L := (0,L)^d$, we  first reduce ourselves to the event
\begin{equation*}
 A = \cur{ \kappa(Q_L) \ge \frac 1 2 }.
\end{equation*}
Indeed, by Markov inequality and \Cref{eq:concentration-abstract}, we easily find, for every $q \ge 1$,
\begin{equation*}
 \PP\bra{ A^c} \les_{q}  L^{-q  },
\end{equation*}
hence, by trivially estimating on $A^c$,
\begin{equation*}
 \WN^{p}_{Q_L}\bra{\kappa(Q_L), \sum_{k=1}^{m^d}\kappa(Q^k_{L_0}) \chi_{Q^k_{L_0}} } \les L^p   \mu(Q_L),
\end{equation*}
we find, by Cauchy-Schwarz inequality,
\begin{equation*}
\frac 1 {L^d} \EE\sqa{ \chi_{A^c}  \WN^{p}_{Q_L}\bra{\kappa(Q_L), \sum_{k=1}^{m^d}\kappa(Q_{L_0}^k)} }\les_q L^{p -q/2} \les  L^{2-d+\alpha}
\end{equation*}
provided that we choose $q$ large enough.

On the event $A$, we use instead \Cref{lem:peyre} with uniform weight $\eta = \chi_{Q_L}$, so that $g_0 = \kappa(Q_L) \ge 1/2$, and we obtain
 \begin{equation*}
   \frac 1 {L^d} \EE\sqa{\chi_A \WN^{p} _{Q_L}\bra{\kappa(Q_L), \sum_{k=1}^{m^d}\kappa(Q_{L_0}^k) \chi_{Q_{L_0}^k}  } }   \les \frac{1}{L^d} \EE\sqa{ \nor{\kappa(Q_L)- \sum_{k=1}^{m^d}\kappa(Q_{L_0}^k) \chi_{Q_{L_0}^k} }_{H^{-1,p}(Q_L)}^p}.
   \end{equation*}
The thesis follows therefore by arguing that
\begin{equation}\label{eq:bound-neg-sob-first-lemma}
 \frac 1 {L^d} \EE\sqa{ \nor{\kappa(Q_L)- \sum_{k=1}^{m^d}\kappa(Q_{L_0}^k) \chi_{Q_{L_0}^k} }_{H^{-1,p}(Q_L)}^p} \les L^{(2+\alpha-d)p/2}_0.
\end{equation}
To this aim, we remark that estimating the $H^{-1,p}$ norm with the $L^p$ norm as in \Cref{eq:poincare-negative} would result in a cruder bound, depending on $m$, hence we need to iterate over dyadic decompositions. For every cube $Q$ of side length $\ell \ge 2$, partitioned into $2^d$ subcubes $(Q^i)_{i=1, \ldots, 2^d}$ each with side length $\ell/2 \ge 1$, we notice that it holds
\begin{equation*}
 \EE\sqa{ \nor{\kappa(Q)- \sum_{i=1}^{2^d}\kappa(Q^i) \chi_{Q^i} }_{H^{-1,p}(Q)}^p } \les \ell^{2+\alpha}.
\end{equation*}
Indeed, it is  sufficient to apply  \Cref{eq:poincare-negative} with $c_P(Q, p) \les \ell$ to obtain
\begin{equation*}
 \begin{split}
  \EE\sqa{ \nor{\kappa(Q)- \sum_{i=1}^{2^d}\kappa(Q^i) \chi_{Q^i} }_{H^{-1,p}(Q)}^p } & \les \ell^p  \EE\sqa{\int_{Q}\sum_{i=1}^{2^d}\abs{\kappa(Q) - \kappa(Q^i)}^p\chi_{Q^i} } \\
  & \les \ell^p \bra{   \ell^d \EE\sqa{\abs{ \kappa(Q)-1}^p} + \sum_{i=1}^{2^d} \bra{\frac{\ell}{2}}^d\EE\sqa{\abs{\kappa(Q^i)-1}^p}}\\
   & \stackrel{\Cref{eq:concentration-abstract}}{\les} \ell^{d+p\lt(1-\frac{d-\alpha}{2}\rt)}. 
 \end{split}
\end{equation*}
%
Starting from $Q = Q_L$, we use the sub-additivity property \Cref{eq:sub-add-sobolev-neg} of $H^{-1,p}$ (with uniform weight)  by further iterating $h$ times the decomposition into cubes of halved side lengths, so that at the $i$-th  iteration the side lengths of the cubes are $\ell_i := L/2^i = L_0 2^{h-i}$  and moreover we can specify a suitable  $\eps = \eps_i \in (0,1)$. At the $h$-th step, we obtain exactly the partition into the $m^d$ sub-cubes  $Q^k_{L_0}$, and we deduce that
\begin{equation}\label{eq:iterated-subadditive}\begin{split}
\frac 1 {L^d} \EE\sqa{ \nor{ \kappa(Q_L) -  \sum_{k=1}^{m^d} \kappa(Q_{L_0}^k) \chi_{Q^k_{L_0}} }_{H^{-1,p}(Q_L)}^p}  &\les  L^{-d} \sum_{i=1}^{h} 2^{id} \lt(L2^{-i}\rt)^{d+p\lt(1-\frac{d-\alpha}{2}\rt)} \eps_i^{-(p-1)} \prod_{j=1}^{i-1} (1+\eps_j) \\
&=\sum_{i=1}^{h} \lt(L2^{-i}\rt)^{p\lt(1-\frac{d-\alpha}{2}\rt)} \eps_i^{-(p-1)} \prod_{j=1}^{i-1} (1+\eps_j)\\
&=L_0^{p\lt(1-\frac{d-\alpha}{2}\rt)}\sum_{i=1}^{h} \lt(2^{h-i}\rt)^{p\lt(1-\frac{d-\alpha}{2}\rt)} \eps_i^{-(p-1)} \prod_{j=1}^{i-1} (1+\eps_j).
\end{split}
\end{equation}
We specify $\eps_i = (2^{h-i})^{-(d-2-\alpha)/2}$, so that the product terms are uniformly bounded from above:
 \begin{equation*}\begin{split}
 \prod_{j=1}^{i-1}(1+\eps_j) & \le \prod_{j=1}^h \bra{ 1+(2^{h-j})^{-\frac{d-2-\alpha}{2}} } \le \exp\bra{  \sum_{j=1}^h 2^{-\frac{h-j}{2}(d-2-\alpha)} }\\
 & \les 1.
 \end{split}
\end{equation*}
Therefore
\begin{equation*}\begin{split}
\sum_{i=1}^{h} \lt(2^{h-i}\rt)^{p\lt(1-\frac{d-\alpha}{2}\rt)} \eps_i^{-(p-1)} \prod_{j=1}^{i-1} (1+\eps_j)& \les  \sum_{i=1}^h \eps_i\les 1.
 \end{split}
\end{equation*}
Combining this with \Cref{eq:iterated-subadditive} this concludes the proof of \Cref{eq:bound-neg-sob-first-lemma}.
%

\subsection{Proof of \Cref{lemma:upper-bound-eta}}

In this case, we fix $m = 2^{h}$ such that $L_0:= L/m \in [1,2)$. Moreover, we write for brevity
\begin{equation*}
 \kappa(\Omega) := \frac{ \int_{\Omega} \eta d\mu}{\int_{\Omega} \eta},
\end{equation*}
with the convention that $\kappa(\Omega) = 0$ if $\int_{\Omega} \eta = 0$.

We preliminarily reduce ourselves to the event
\begin{equation*}
 A := \cur{  \kappa(Q_L) \ge \frac 1 4 }.
\end{equation*}
By the properties of $\eta$, in particular the fact that $\eta$ is identically $1$ if $\dist(x,Q_L^c) \le \delta L$ and $0$ if $\dist(x, Q_L^c) \ge 2 \delta L$, we see that
\begin{equation*}
\cur{ \mu(\dist(\cdot ,Q_L^c) \le \delta L) \ge \delta L^d/2 } \subseteq A,
\end{equation*}
hence by Markov inequality  and \Cref{eq:concentration-abstract} we find, for every $q\ge 1$,
\begin{equation*}
 \PP(A^c) \les_{\delta,q} L^{-q}.
\end{equation*}
On the event $A$, we partition $Q_L$ into $m^d$ disjoint subcubes $Q_L = \bigcup_{k=1}^{m^d} Q_{L_0}^k$, each of side length $L_0$, and use \Cref{eq:Neumann-exact-sub-additive} with $\eta \mu$ instead of $\mu$, $\eta \int_{Q_L} \eta \mu/ \int_{Q_L}\eta$ instead of $\lambda$, and $\eps=1/2$.
 We find
\begin{equation}\label{eq:prop-michael-first-subadditive}
 \WN^{p}_{Q_L}\bra{\eta \mu, \kappa(Q_L) \eta  } \les \sum_{k}^{m^d} \WN^{p}_{Q_{L_0}^k}\bra{ \eta \mu,\kappa(Q_{L_0}^k)  \eta  } +  \WN^p_{Q_L}\bra{ \sum_{k=1 }^{m^d} \kappa(Q_{L_0}^k)  \chi_{Q^k_{L_0} } \eta,\kappa(Q_L) \eta  }.
\end{equation}
Taking expectation (on $A$) and using that $L_0 \le 2$, we find
\begin{equation*}
 \EE\sqa{\WN^{p}_{Q_{L_0}^k}\bra{ \eta \mu, \kappa(Q_{L_0}^k)  \eta  } } \les \EE\sqa{ \int_{Q^k_{L_0}} \eta \mu} = \int_{Q^k_{L_0}} \eta
\end{equation*}
so that the first contribution in the right-hand side of \Cref{eq:prop-michael-first-subadditive} can be bounded as
\begin{equation*}
 \EE\sqa{ \sum_{k}^{m^d} \WN^{p}_{Q_{L_0}^k}\bra{ \eta \mu, \kappa(Q_{L_0}^k)  \eta  } } \les \int_{Q_L} \eta \les L^d \delta,
\end{equation*}
by the properties of $\eta$. After dividing by $L^d$, this yields a first contribution in  \Cref{claimmatching}. For the second term, we apply \Cref{lem:peyre}, so that
\begin{equation*}
 \EE\sqa{ \chi_A \WN^p_{Q_L}\bra{ \sum_{k=1 }^{m^d} \kappa(Q_{L_0}^k)  \chi_{Q^k_{L_0} } \eta,\kappa(Q_L) \eta  } } \les \EE\sqa{ \nor{ \sum_{k=1 }^{m^d} \kappa(Q_{L_0}^k)  \chi_{Q^k_{L_0} } - \kappa(Q_L) }_{H^{-1,p}(\eta)}^p} .
\end{equation*}
Let us notice that, as we argue on $A$, we can use the lower bound on $\kappa(Q_L)$ yielding $g_0 =1/4$ in \Cref{lem:peyre}, and moreover $c_P(\eta, p)< \infty$ by \Cref{rem:finite-weighted-poincare}.

To complete the proof, we argue that
\begin{equation}\label{weightedH-1pmatching}
 \frac 1 {L^d} \EE\sqa{ \nor{ \sum_{k=1 }^{m^d} \kappa(Q_{L_0}^k)  \chi_{Q^k_{L_0} } - \kappa(Q_L) }_{H^{-1,p}(\eta)}^p} \les \delta +c(\delta) L^{(2-d+\alpha)p/2}.
\end{equation}
Let us notice that this can be regarded as a weighted analogue of \Cref{eq:bound-neg-sob-first-lemma}, and indeed, as in that case, we iterate over dyadic decompositions. First, we notice that for any given cube $Q$ of side length $\ell \ge 2$, partitioned into $2^d$ subcubes $(Q^i)_{i=1, \ldots, 2^d}$ each with side length $\ell/2 \ge 1$,  it holds
\begin{equation*}
 \EE\sqa{ \nor{\kappa(Q)- \sum_{i=1}^{2^d}\kappa(Q^i) \chi_{Q^i} }_{H^{-1,p}(\chi_Q \eta )}^p } \les c_P(\chi_Q \eta, p)^p \ell^{d}\ell^{-\frac{p}{2}(d-\alpha)}.
\end{equation*}
Indeed, it is  sufficient to apply   \Cref{eq:poincare-negative} and obtain
\begin{equation*}
 \begin{split}
  \EE\sqa{ \nor{\kappa(Q)- \sum_{i=1}^{2^d}\kappa(Q^i) \chi_{Q^i} }_{H^{-1,p}(Q)}^p } & \les c_P(\chi_Q \eta, p)^p \EE\sqa{\int_{Q}\sum_{i=1}^{2^d}\abs{\kappa(Q) - \kappa(Q^i)}^p\chi_{Q^i} \eta } \\
  & \les c_P(\chi_Q \eta, p)^p \bra{  \EE\sqa{\abs{ \kappa(Q)-1}^p} \int_Q \eta  + \sum_{i=1}^{2^d}  \EE\sqa{\abs{\kappa(Q^i)-1}^p}\int_{Q^i} \eta }\\
   & \stackrel{\Cref{eq:concentration-abstract}}{\les} c_P(\chi_Q \eta, p)^p \ell^{d}\ell^{-\frac{p}{2}(d-\alpha)}. 
 \end{split}
\end{equation*}
 The crux of the argument is to obtain suitable bounds for the possibly different constants $c_P(\chi_Q \eta, p)$ that appear in the subsequent dyadic decompositions of $Q_L$. We notice that, if $\eta$ is identically null on $Q$, we can trivially estimate
\begin{equation*}
  \EE\sqa{ \nor{\kappa(Q)- \sum_{i=1}^{2^d}\kappa(Q^i) \chi_{Q^i} }_{H^{-1,p}(Q)}^p } = 0,
\end{equation*}
hence it is sufficient to consider only the contributions of the cubes $Q$ where $\eta$ is not identically null. For every $i=1, \ldots, h$, let us denote with $\cQ_{i}$ the collection of cubes of side length $L2^{-i+1}$ that are obtained by iterated dyadic decomposition of $Q_L$ and which intersect the support of $\eta$. Clearly, we have $\sharp \cQ_i \le 2^{id}$ and moreover we can always bound from above (using that $\eta$ is defined as in \Cref{eq:eta-d-L-delta})
\begin{equation*}
 \sup_{Q \in \cQ_i} \frac{ c_P(\chi_Q \eta, p)}{|Q|^{1/d}} \les_\delta 1.
\end{equation*}
Let us recall that $\delta=2^{-h_0}$. For $i\ge h_0$, we have  $\sharp \cQ_i \le \delta 2^{id}$ and we claim that
\begin{equation}\label{eq:bound-poincare}
  \sup_{Q \in \cQ_i} \frac{ c_P(\chi_Q \eta, p)}{|Q|^{1/d}} \les 1,
\end{equation}
i.e., the (rescaled) Poincaré constant becomes uniformly controlled with respect to $\delta$. Let us postpone the proof of \Cref{eq:bound-poincare} and conclude the proof of \Cref{weightedH-1pmatching}. Using \Cref{eq:bound-poincare} and arguing as in \Cref{eq:iterated-subadditive}, we collect the inequality
\begin{equation*}\begin{split}
   \EE\sqa{ \nor{ \kappa(Q_L) -  \sum_{k=1}^{m^d} \kappa(Q_{L_0}^k) \chi_{Q^k_{L_0}} }_{H^{-1,p}(Q_L)}^p}  &\le  c(\delta) \sum_{  i=1 }^{h_0-1} \lt(L2^{-i}\rt)^{p\lt(1-\frac{d-\alpha}{2}\rt)} \eps_i^{-(p-1)} \prod_{j=1}^{i-1} (1+\eps_j)\\
   & \quad +   \delta \sum_{ i = h_0  }^h \lt(L2^{-i}\rt)^{p\lt(1-\frac{d-\alpha}{2}\rt)} \eps_i^{-(p-1)} \prod_{j=1}^{i-1} (1+\eps_j). \\
   \end{split}
 \end{equation*}
For $i=1,\cdots, h_0-1$, we choose $\eps_i=1/2$ so that (recall that $h_0$ depends only on $\delta$)
\begin{equation*}
 \sum_{  i=1 }^{h_0-1} \lt(L2^{-i}\rt)^{p\lt(1-\frac{d-\alpha}{2}\rt)} \eps_i^{-(p-1)} \prod_{j=1}^{i-1} (1+\eps_j) \les_\delta L^{ (2+\alpha-d)p/2}.
\end{equation*}
For $i=h_0,\cdots, h$, we choose instead  $\eps_i = (L/2^i)^{-(d-2-\alpha)/2}$, so that arguing as above we find
\begin{equation*}
 \sum_{ i = h_0  }^h \lt(L2^{-i}\rt)^{p\lt(1-\frac{d-\alpha}{2}\rt)} \eps_i^{-(p-1)}\prod_{j=1}^{i-1} (1+\eps_j) \les 1.
\end{equation*}
%

This would prove \Cref{weightedH-1pmatching} and we are only left with the proof of \Cref{eq:bound-poincare}. For this we write $\Gamma := \cur{x \in Q_L: \eta(x) = 0}$, that is a cube of side length $(1-2\delta)L$. Given a cube $Q\in  \cQ_i$ with side length $\ell = 2^{-i}L \le \delta L$, we distinguish between two cases. If $\partial Q \cap \Gamma =\emptyset$, then by construction (in particular using that $\delta = 2^{-h_0}$ is also dyadic) we have $d(Q,\Gamma)\ges \ell$. In particular, by \Cref{etasimdist}, for every $x,y\in Q$, we have
\begin{equation*}
 \frac{\eta(x)}{\eta(y)}\sim \lt(\frac{d(x,\Gamma)}{d(y,\Gamma)}\rt)^3\les \lt(\frac{\ell +d(Q,\Gamma)}{d(Q,\Gamma)}\rt)^3\les 1.
 \end{equation*}
By \Cref{comparePoincare}, up to multiplication by a constant, translation and scaling, \Cref{eq:bound-poincare} then follows from the standard Poincar\'e inequality in $(0,1)^d$. If otherwise $\partial Q\cap  \Gamma \neq \emptyset$, set for $k\in \cur{1,\cdots, d}$, $\Gamma_k=\{0\}^k\times(0,1)^{d-k}$ and then for $x\in (0,1)^d$, $\eta_k(x)=d^3(x,\Gamma_k)$. By \Cref{comparePoincare} and  \Cref{etasimdist}, we see that after suitable translation, scaling and multiplication by a constant, it is enough to prove that, for every  $k\in \{1,\ldots, d\}$,
\begin{equation}\label{Poink}
 c_P(\eta_k\chi_{(0,1)^d}, p)\les 1.
\end{equation}
But this follows at once from \Cref{rem:finite-weighted-poincare}.

\begin{remark}
The proofs of \Cref{lemma:upper-bound-mu-intermediate-cubes} and \Cref{lemma:upper-bound-eta} have a similar structure, and could be combined together, but we prefer to keep them separate to avoid adding technicalities. Moreover, in the range $d> 4 +\alpha$ one could rely upon easier bounds, estimating the Wasserstein distance in terms of the total variation \Cref{eq:wass-tv} and  arguing as in \parencite{barthe2013combinatorial}. However this would leave out from the application to the transport of i.i.d.\ points in dimension $d \in \cur{3,4}$, and that of the Brownian interlacement occupation measure for $d \in \cur{5,6}$.
\end{remark}

\subsection{Bounds for the Brownian interlacement occupation measure}\label{sec:ass-mu}

The aim of this section is to provide a proof of the validity of \Cref{ass:mu}, with $\alpha=2$, for the Brownian interlacement occupation measure. As already noticed, the validity of \Cref{ass:mu} with $\alpha=0$ for a Poisson point process is rather  straightforward: however, since it may give some intuition for the analogue derivations in the case of the Brownian interlacement occupation measure, we report here the argument leading to \Cref{eq:concentration-abstract}. Indeed, given a Poisson point process $\mu$, a cube $Q$ and $\eta: Q \to [0,\infty)$ bounded, we write
\begin{equation*}
 \int_{Q} \eta d \mu =  \sum_{i=1}^{N(Q)} \eta(X_i),
\end{equation*}
as a sum of a Poisson number $N(Q)$, with mean $|Q|$, of i.i.d.\ variables $\eta(X_i)$ where and each $X_i$ is uniformly distributed on $Q$. Hence, $\EE\sqa{ \int_Q \eta d \mu} =  \int_Q \eta$ and we can write
\begin{equation*}
 \int_{Q} \eta d \mu - \EE\sqa{ \int_Q \eta d \mu}   = \sum_{i=1}^{N(Q)} \bra{ \eta(X_i) - \fint_Q \eta} + \bra{  N(Q) - |Q|} \fint_Q \eta.
 \end{equation*}
 The variables $Y_i := \eta(X_i) - \fint_Q \eta$ are centered and uniformly bounded by $\nor{\eta}_\infty$. Thus, by Rosenthal's inequality, for any $q\ge 2$ and $n \ge 1$,
 \begin{equation*}
  \nor{\sum_{i=1}^n Y_i}_{L^q(\PP)} \les_q \bra{ n^{1/q}+n^{1/2}} \nor{\eta}_\infty \les n^{1/2} \nor{\eta}_\infty.
 \end{equation*}
 Conditioning upon $N(Q) = n$ and using that for a Poisson variable $\EE\sqa{N(Q)^{q/2}} \les_q |Q|^{q/2}$, we then find
 \begin{equation*}
  \nor{\sum_{i=1}^{N(Q)} Y_i}_{L^q(\PP)} \les_q |Q|^{1/2} \nor{\eta}_\infty.
 \end{equation*}
 Finally, by the concentration properties of the Poisson law,
 \begin{equation*}
  \nor{\bra{  N(Q) - |Q|} \fint_Q \eta}_{L^q(\PP)} \les \nor{\bra{  N(Q) - |Q|} }_{L^q(\PP)} \nor{\eta}_\infty \les_q |Q|^{1/2} \nor{\eta}_\infty
 \end{equation*}
and \Cref{eq:concentration-abstract} is settled.
%

Moving to the Brownian interlacement occupation measure case, we first collect some facts about classical potential theory and the interlacement process, see \parencite{mariani2023wasserstein, drewitz2014introduction, sznitman2013scaling, port2012brownian}. We write $p_t(x,y)$ for the transition function of Brownian motion and
\begin{equation*}
 g(x,y) = \int_0^\infty p_t(x,y) d t = \frac{c_d}{|x-y|^{d-2}}
\end{equation*}
for the Brownian Green's kernel (where $c_d$ is a constant that depends on the dimension). Given $x \in \R^d$, we write $\PP_x, \EE_x$ for the law and expectation of Brownian motion started in $x \in \R^d$, similarly $\PP_\lambda, \EE_\lambda$ for any finite measure $\lambda$, not necessarily a probability.

For a compact $K\subset \R^d$ we denote by $e_K$ its equilibrium measure, which is uniquely defined by requiring that the associated potential
\begin{equation*}
 V(e_K)(x)=\int_{\R^d} g(x,y) d e_K(y)
\end{equation*}
satisfies $V(e_K) = 1$ on $K$ (and $V(e_k)(x)$ is infinitesimal as $|x| \to \infty$).

The mass $e_K(\R^d)=\Cap(K)$ is called the capacity of $K$ and we denote by $\tilde e_K=\frac{1}{\operatorname{Cap}(K)} e_K$ the normalised measure, assuming that $\operatorname{Cap}(K)\in (0,\infty)$. Note that the characterising property of the $e_K$ implies that
\begin{equation*}
\Delta^{-1} e_K \begin{cases} =1 & \text{ on } K \\ \leq 1 & \text{ on } \R^d \end{cases}.
\end{equation*}

Given a continuous curve $\omega = (\omega_t)_{t\ge 0}$, write $\tau_K(\omega)$ for the first hitting time of $K$. For any $x\in \R^d$ we have the fundamental identity
\begin{equation}\label{eq:fundamental}
\PP_y(\tau_K<\infty)=\int_{\R^d}g(x,y)de_K(x)
\end{equation}
relating Brownian hitting probabilities, Green's function and the equilibrium measure. Denoting with $\mu^B_\infty=\int_0^\infty \delta_{B_s} ds$ the occupation measure of a Brownian path, i.e., any Borel set $A\subset \R^d$,
\begin{equation*}
 \mu^B_\infty(A) = \int_0^\infty \chi_A\bra{B_s} ds
\end{equation*}
identity \Cref{eq:fundamental} implies
\begin{equation*}
\begin{split}
\EE_{e_K}\sqa{\mu_\infty^B(A)} & = \int_{\R^d} \int_0^\infty\int_{\R^d} \chi_A(y) p_t(x,y)dy dt de_K(x) \\
& = \int_{\R^d \times \R^d} g(x,y) \chi_A(y) dy de_K(x) \\
& = \int_A \PP_y[\tau_K<\infty] dy \leq |A|.
\end{split}
\end{equation*}
Hence, we deduce the inequality between measures
\begin{equation*}
\EE_{e_K}\sqa{\mu_\infty^B(A)} \le |A|,
\end{equation*}
with equality if $A \subseteq K$. Moving from sets to functions, we have for any $f\geq 0$ that
\begin{equation}\label{eq:eKleLeb}
\EE_{e_K}\sqa{\int_0^\infty f(B_s) ds} \leq \int_{\R^d} f.
\end{equation}


Next, we recall that given any compact $K\subseteq \R^d$, the restriction of the Brownian interlacement occupation measure $\mu$ on $K$ can be represented as
\begin{equation*}
\mu\restr K =\sum_{i=1}^{N_K} \mu_\infty^{B^i} \restr K,
\end{equation*}
where $N_K$ is a Poisson random variable with mean $\EE\sqa{N_K} = \Cap(K)$ and $(B^i)_{i=1}^\infty$ are independent Brownian motions, independent of $N_K$, each with initial law $\tilde e_K$. In particular, we get for $A\subset K$,
\begin{equation*}
\EE\sqa{\mu(A)}= \Cap(K) \frac{|A|}{\Cap(K)} = |A|
\end{equation*}
and similarly, for any function $\eta$ supported in $K$, it holds
\begin{equation*}
\EE\sqa{\int_K \eta d\mu } = \int_K\eta.
\end{equation*}

We are now in a position to prove the concentration inequality \Cref{eq:concentration-abstract}.

\begin{lemma}
Let $d \ge 3$ and $\mu$ be the  Brownian interlacement occupation measure on $\R^d$. Then, for any cube $Q$ with $\diam(Q) \ge 1$ and any non-negative and bounded $\eta$ and any $q\geq 2$,  it holds
\begin{equation*}
\nor{\int_Q\eta d\mu - \EE\sqa{\int_Q\eta d\mu}}_{L^q(\PP)} \leq c\diam(Q)^{\frac{d+2}{2}} \nor{\eta}_\infty,
\end{equation*}
where $c = c(q,d)<\infty$.
\end{lemma}

\begin{proof}
By considering $\eta \chi_Q$ we can assume $\eta$ to be supported on $Q$. We write
\begin{equation*}
\int \eta d\mu = \sum_{i=1}^N Y_i
\end{equation*}
for iid random variables $Y_i= \int \eta(B^i_s) ds$ and an independent Poisson random variable $N$ with parameter $\Cap(Q)=\diam(Q)^{d-2}\Cap(Q_1)\sim \diam(Q)^{d-2}$. Since $0\leq\int \eta(B_s) ds \leq \nor{\eta}_\infty \int \chi_Q(B_s) ds$ it follows from Lemma 3.2 in \parencite[Lemma 3.2]{mariani2023wasserstein} that $Y_1$ has finite moments of all orders, in fact that
\begin{equation*}
\EE[Y_1^q] \les_q \nor{\eta}_\infty^q \diam(Q)^{2q}.
\end{equation*}
Applying Rosenthal's inequality  as we did for the Poisson point process case -- see also \parencite[eq. (2.57)]{mariani2023wasserstein} -- it follows that  $\int \eta d\mu $ has moments of all orders and that
\begin{equation*}\begin{split}
\nor{\int_Q\eta d\mu - \EE\sqa{\int_Q\eta d\mu}}_{L^q(\PP)} & \les_q \diam(Q)^{\frac{d-2}{q}}\nor{Y_1}_{L^q(\PP)} + \diam(Q)^{\frac{d-2}{2}}\nor{Y_1}_{L^2(\PP)}\\
& \les_q \nor{\eta}_\infty \diam(Q)^{\frac{d+2}{2}}. \qedhere
\end{split}
\end{equation*}
\end{proof}

\section*{Acknowledgments and funding}

Part of this work was conceived at the CIRM 2988 conference ``PDE \& Probability in interaction: functional inequalities, optimal transport and particle systems''. The authors warmly thank the conference organizers and CIRM staff for setting up a stimulating environment. Among the participants, they thank J.-C.\ Mourrat and A.\ Delalande for fruitful comments on their respective research works. The authors thank L.~Ambrosio and O.~Mischler for useful comments and remarks on a preliminary version of the work.

M.G.\ and D.T.\ acknowledge the project  G24-202 ``Variational methods for geometric and optimal matching problems'' funded by Università Italo Francese.

M.H.\ is supported by the Deutsche Forschungsgemeinschaft (DFG, German Research Foundation) under Germany’s Excellence Strategy EXC 2044 -390685587, \textit{Mathematics M\"unster: Dynamics–Geometry–Structure} and by the DFG through the SPP 2265: \textit{Random Geometric Systems}.

D.T.\ acknowledges the MUR Excellence Department Project awarded to the Department of Mathematics, University of Pisa, CUP I57G22000700001,  the HPC Italian National Centre for HPC, Big Data and Quantum Computing -  CUP I53C22000690001, the PRIN 2022 Italian grant 2022WHZ5XH - ``understanding the LEarning process of QUantum Neural networks (LeQun)'', CUP J53D23003890006, the INdAM-GNAMPA project 2024 ``Tecniche analitiche e probabilistiche in informazione quantistica''.  Research also partly funded by PNRR - M4C2 - Investimento 1.3, Partenariato Esteso PE00000013 - "FAIR - Future Artificial Intelligence Research" - Spoke 1 "Human-centered AI", funded by the European Commission under the NextGeneration EU programme.

\printbibliography

\end{document}